\documentclass[11pt]{article}
\usepackage{epsfig,tikz,float,wrapfig}
\usepackage{amssymb,amsmath}
\usepackage[utf8]{inputenc}
\usepackage[OT1,T2A]{fontenc}
\makeatletter
\newdimen\normalarrayskip            
\newdimen\minarrayskip               
\normalarrayskip\baselineskip
\minarrayskip\jot
\newif\ifold             \oldtrue            \def\new{\oldfalse}
\def\arraymode{\ifold\relax\else\displaystyle\fi}
\def\eqnumphantom{\phantom{(\theequation)}} 
\def\@arrayskip{\ifold\baselineskip\z@\lineskip\z@
     \else
     \baselineskip\minarrayskip\lineskip1\baselineskip\fi}

\def\@arrayclassz{\ifcase \@lastchclass \@acolampacol \or
\@ampacol \or \or \or \@addamp \or
   \@acolampacol \or \@firstampfalse \@acol \fi
\edef\@preamble{\@preamble
  \ifcase \@chnum
     \hfil$\relax\arraymode\@sharp$\hfil
     \or $\relax\arraymode\@sharp$\hfil
     \or \hfil$\relax\arraymode\@sharp$\fi}}

\def\@array[#1]#2{\setbox\@arstrutbox=\hbox{\vrule
     height\arraystretch \ht\strutbox
     depth\arraystretch \dp\strutbox
width\z@}\@mkpream{#2}\edef\@preamble{\halign \noexpand\@halignto
\bgroup \tabskip\z@ \@arstrut \@preamble \tabskip\z@ \cr}%
\let\@startpbox\@@startpbox \let\@endpbox\@@endpbox
  \if #1t\vtop \else \if#1b\vbox \else \vcenter \fi\fi
  \bgroup \let\par\relax
  \let\@sharp##\let\protect\relax
  \@arrayskip\@preamble}
\def\eqnarray{\stepcounter{equation}%
              \let\@currentlabel=\theequation
              \global\@eqnswtrue
              \global\@eqcnt\z@
              \tabskip\@centering              
              \let\\=\@eqncr
              $$%
            \halign to \displaywidth  \bgroup
             \eqnumphantom \@eqnsel
      \hskip\@centering                               
    $\displaystyle  \tabskip\z@ {##}$%
    &\global\@eqcnt\@ne \hskip 2\arraycolsep
         $ \displaystyle  \arraymode{##}$\hfil
    &\global\@eqcnt\tw@ \hskip 2\arraycolsep
         $\displaystyle\tabskip\z@{##}$\hfil
         \tabskip\@centering
    &{##}\tabskip\z@\cr}
\makeatother

\def\Tr{{\rm Tr}}
\def\Ad{{\rm Ad}}
\def\sign{{\rm sign}}

\def\Bf#1{\mbox{\boldmath $#1$}}

\def\bxi{{\Bf\xi}}
\def\Beta{{\Bf\eta}}

\def\rank{{\rm rank}}
\def\diag{{\rm diag}}

\def\beq{\begin{equation}}
\def\eeq{\end{equation}}
\def\be{\beq\new\begin{array}{c}}
\def\ee{\end{array}\eeq}

\newcommand{\rf}[1]{(\ref{#1})}
\begin{document}
\newtheorem{lemma}{Lemma}
\newtheorem{proposition}{Proposition}
\newtheorem{theorem}{Theorem}

\begin{flushright}
ITEP/TH-01/14\\
FIAN/TD-01/14
\end{flushright}

\begin{center}
\baselineskip20pt
{\bf \LARGE Loop groups, Clusters, Dimers and Integrable systems}
\end{center}
\bigskip
\begin{center}
\baselineskip12pt
{\large V.~V.~Fock$^{a,b}$\footnote{on leave of absence} and A.~Marshakov$^{c,b,d}$}\\
\end{center}
\bigskip

\begin{center}
$^a${\it IRMA, Universit\'e de Strasbourg, France}\\
$^b${\it Institute for Theoretical and Experimental Physics, Moscow, Russia}\\
$^c${\it Theory Department, Lebedev Physics Institute, Moscow, Russia}\\
$^d${\it Department of Mathematics, NRU HSE, Moscow, Russia}\\
\end{center}
\bigskip\medskip

\begin{center}
{\large\bf Abstract} \vspace*{.2cm}
\end{center}

\begin{quotation}
\noindent
We describe a class of integrable systems on Poisson submanifolds of the affine Poisson-Lie groups $\widehat{PGL}(N)$, which can be enumerated by cyclically irreducible elements the co-extended affine Weyl groups $(\widehat{W}\times \widehat{W})^\sharp$. Their phase spaces admit cluster coordinates, whereas the integrals of motion are cluster functions. We show, that this class of integrable systems coincides with the constructed by Goncharov and Kenyon out of dimer models on a two-dimensional torus and classified by the Newton polygons. We construct the correspondence between the Weyl group elements and polygons, demonstrating that each particular integrable model admits infinitely many realisations on the Poisson-Lie groups. We also discuss the particular examples, including the relativistic Toda chains and the Schwartz-Ovsienko-Tabachnikov pentagram map.
\end{quotation}

\tableofcontents

\section{Introduction}

The main idea of the paper is to demonstrate the equivalence of two \emph{a priori} different methods of construction and description of a wide class of integrable models, and thus -- to propose the unified approach of their investigation. In the first well-known method \cite{Sembook} the phase space is taken as a quotient of double Bruhat cells of a Kac-Moody Lie group, with the Poisson structure defined by a classical $r$-matrix, and the integrals of motion are just the Ad-invariant functions. The second method was suggested recently by A.~Goncharov and R.~Kenyon \cite{GK}, and it grows up out of the study of dimer models of statistical physics on bipartite graphs on a two-dimensional torus. We are going to show that in fact the latter class of integrable systems is a particular case of the former, corresponding to the affine group of type $\hat{A}_{N-1}$.

The best known example of the integrable system of this class \cite{FM} is the relativistic Toda chain, discovered by S.~Ruijsenaars \cite{Ruj} and studied by Yu.~Suris \cite{Suris} and many others, which gives the common Toda chain in a certain limit, corresponding in our terms to passing to Lie algebras from the Lie groups. Another known example is the pentagram map -- a discrete integrable system on the space of broken lines in a projective plane, discovered by R.~Schwartz \cite{Schwartz} and studied by him with V.~Ovsienko and S.~Tabachnikov \cite{OT}-\cite{OT2}, and recently by M.~Glick \cite{Glick}, B.~Khesin and F.Soloviev \cite{Sol}\cite{KheSol}, M.~Shapiro, M.~Gekhtman, A.~Vainshtein and S.Tabachnikov \cite{GeShaLast}. On the other hand, V.~Ovsienko and S.~Tabachnikov have shown in \cite{OT}, that the discrete flow of the pentagram map gives the Boussinesq flow in the continuum limit. This observation, generalised in \cite{KheSol} for other dimensions, suggests that the technique, proposed in the paper, can be also applied to study of continuous integrable systems such as $n$-KdV hierarchies.

\subsection{Integrable systems on Poisson-Lie groups}

Our starting observation is that on a Poisson-Lie group, with the Poisson bracket defined by classical $r$-matrix, the Ad-invariant functions Poisson commute with each other. For a finite dimensional simple group the number of independent Ad-invariant functions is equal to the rank of the group, and thus the  corresponding integrable system can arise on a symplectic leaf of dimension not more than twice the rank \cite{FM}. However, for affine groups the number of independent Ad-invariant functions is infinite though all Poisson submanifolds are still finite dimensional, and thus the set of integrable systems one gets in this way is much larger. Such integrable systems can be constructed on any affine Poisson-Lie group $\hat{G}$, but in the paper we will need only the systems on the groups of type $\hat{A}$ with trivial center. This group can be realised as a group of matrix-valued Laurent polynomials $\mathcal{A}(\lambda)$ of a signle variable with nonzero constant determinant and considered up to multiplication by a nonzero constant. For a given $\mathcal{A}(\lambda)$ the set $\{(\lambda,\mu)|\det (\mathcal{A}(\lambda)-\mu)=\sum_{ij}\mathcal{H}_{ij}\lambda^i\mu^j=0\}$ is an algebraic curve in $(\mathbb{C}^\times)^2$, called \textit{spectral curve}, embedded into the torus $(\mathbb{C}^\times)^2$. This curve, considered up to the torus automorphisms, is a conjugacy class invariant. The spectral curve comes together with the line bundle, given by the kernel of $\mathcal{A}(\lambda)-\mu\cdot{\rm Id}$. The map from the group $\hat{G}$ to the space of curves is called the action map, and it is a Poisson map if we take a trivial Poisson bracket on the space of curves. The functions $\mathcal{H}_{ij}$ themselves are not well defined, since changing them by $\mathcal{H}_{ij}\to \mathcal{H}_{ij}\alpha^i\beta^j\gamma$ would correspond to the same curve, however one can use this freedom to make any three nonvanishing coefficients to be equal to unities. With this condition $\mathcal{H}_{ij}$ become well defined and do Poisson-commute. The map to the pair (curve, line bundle on it) is called the action-angle map. The flows generated by the Poisson commuting integrals of motion or Hamiltonians amount to the constant flow of the line bundle along the Jacobian of the spectral curve. We describe such integrable systems in more detail in sect.~\ref{s:isrmatr}.

A loop group $\hat{G}$ does not have a cluster variety structure. However it is embedded as a Poisson submanifold into a central coextension $\hat{G}^\sharp$, which admits a standard decomposition into disjoint union of the Poisson submanifolds $(\hat{G}^\sharp)^u$, called double Bruhat cells \cite{ZelKog} which are already the cluster varieties. These cells are enumerated by the elements $u$ of a coextension $(\widehat{W}\times\widehat{W})^\sharp$ of the square of the Weyl group $\widehat{W}$ of $\hat{G}$ by the automorphism group of the Dynkin diagram, which is a cyclic group for the series $\hat{A}_{N-1}$.

Intersections of the Bruhat cells of $\hat{G}^\sharp$ with $\hat{G}$, quotiented by conjugation by the finite dimensional Cartan subgroup $H$, (we call them also the double Bruhat cells and denote by $\hat{G}^u$) are the phase spaces $\hat{G}^u/\Ad H$ of our integrable systems. The dimension of such phase space is (in the case of affine groups) one less, than the length of $u$. Given a presentation of $u$ as a reduced word of the standard generators, one can define the cluster coordinates $\boldsymbol{x}=\{x_f\}$, enumerated by the letters of the word (except for the generator of the coextension) and subject to $\prod x_f=1$.  For a given Bruhat cell the Laurent polynomial $\det (\mathcal{A}(\boldsymbol{x},\lambda,\mu)=\sum \mathcal{H}_{ij}(\boldsymbol{x})\lambda^i\mu^j$ has a fixed Newton polygon $\Delta$, and the number of the Poisson commuting Hamiltonians, nonvanishing on the cell, is three less than the number of integral points of $\Delta$. In order for this system to be integrable these integrals of motion must be independent, and their number should be maximal possible $(\dim \hat{G}^u+\mbox{corank }\hat{G}^u)/2$ for the given dimension of the cell and given rank of the Poisson bracket). This condition is satisfied on double Bruhat cells, corresponding to $u$ having minimal length in its conjugacy class (such elements are also called \textit{cyclically irreducible}).

In sect.~\ref{s:loop} following \cite{Drinfeld} we introduce cluster coordinates (cluster seeds) on the double Bruhat cells of the group $\hat{G}^\sharp$ and on their quotients by conjugation by the Cartan subgroup (desribing them first for the simple groups in sect.~\ref{ss:clustcell} and  generalising then to the affine case). Given a set of coordinates $\boldsymbol{x}$, corresponding to a reduced word, we construct the corresponding matrix polynomial $\mathcal{A}(\boldsymbol{x},\lambda)$ as a product of elementary matrices each of which is either a constant or depends on just a signle coordinate $x_f$. Then we formulate our integrable systems in terms of these cluster coordinates. Coordinates corresponding to different decomposition of the word $u$ are related by cluster transformations.

\subsection{Goncharov-Kenyon integrable systems}

Recall that a cluster variety is an algebraic variety, covered by charts isomorphic to algebraic tori $(\mathbb{C}^\times)^N$ with transition function being compositions of special birational transformations called \textit{mutations} (see Appendix \ref{ap:cluster}). The Goncharov-Kenyon (GK) approach proposes an integrable system, associated with a dimer model on a graph on two-dimensional torus. The integrable system structure turns out to be compatible
with the structure of cluster variety -- implying duality, discrete group action, quantisation, tropical limit and many other attractive features.  Moreover, the mutually Poisson-commuting integrals of motion can be chosen to be the cluster functions. Many such systems admit an abelian group of discrete cluster transformations, commuting with the integrable flows.

The scheme of the construction of the GK integrable system, described in detail in sect.~\ref{s:dimers}, is roughly as follows. The starting point is a bipartite graph on a closed surface $\Sigma$, satisfying certain minimality and non-triviality conditions (following \cite{GK} we  consider the surface $\Sigma$ to be of genus one). Many aspects of the construction can be generalised for surfaces of higher genus, but we postpone this generalisation for a forthcoming publication.

Consider the space of discrete connections on the graph $\Gamma$ with values in the multiplicative group (to be specific we assume it to be the multiplicative group of non-zero complex numbers $\mathbb{C}^\times$). Since every edge of a bipartite graph can be canonically oriented, say from white to black vertex, we can interpret this space as the multiplicative group of one cochains $C^1(\Gamma)$, which is just the product of $\mathbb{C}^\times$ factors, corresponding to each edge of $\Gamma$.
The quotient of this space by discrete gauge transformations can be interpreted as the cohomology group $H^1(\Gamma)$. This group is an extension of the group of coboundaries $B^2(\Sigma)$ by the cohomology group of torus $H^1(\Sigma)=\mathbb{C}^\times\times\mathbb{C}^\times$. The elements of $B^2(\Sigma)$ are collections of numbers on faces of $\Sigma$, with their product equal to unity. The differential $\partial H^1(\Gamma)\to B^2(\Sigma)$ (given by monodromies around the faces) defines on $H^1(\Gamma)$ the structure of a principal $H^1(\Gamma)$-bundle over $B^2(\Sigma)$.

For every discrete connection $\boldsymbol{A}=\{A_e\}\in C^1(\Gamma)$ and for every choice of a spin structure on torus $\Sigma$ we define in sect.~\ref{ss:dirac} a discrete Dirac operator  $\mathfrak{D}(\boldsymbol{A}):\mathbb{C}^B\to\mathbb{C}^W$ from the functions on black vertices to the functions on white vertices. This operator degenerates on a subvariety of  $C^1(\Gamma)$ which is the vanishing locus of the determinant $\det \mathfrak{D}(\boldsymbol{A})$. This variety is gauge invariant and thus it defines a subvariety in $H^1(\Gamma)$. Intersection of this locus with a fibre over a point $\boldsymbol{x}=\{x_f\}\in B^2(\Sigma)$ gives an algebraic curve with line bundles on it given by the kernel of $\mathfrak{D}(\boldsymbol{A})$. 
Observe, that the determinant of $\mathfrak{D}(\boldsymbol{A})$ is a sum of monomials over the perfect matchings of the white and a black vertices, and this is how the dimer configurations on the graph $\Gamma$ come into play. The spin structure permits to control signs of the monomials.
This construction  defines therefore the action map of the phase space of our integrable system $B^2(\Sigma)$ to the space of plane algebraic curves and the action-angle map to the pairs (plane curve, line bundle on it) which is a birational isomorphism.

In coordinates the action map reads as follows. Choose a trivialisation of the bundle $\partial$ which amounts to an isomorphism $H^1(\Gamma)=B^2(\Sigma)\times H^1(\Sigma)$ and then choose a lift $H^1(\Gamma)$ to $C^1(\Gamma)$. Under these identifications we can associate a connection $\boldsymbol{A}(\boldsymbol{x},\boldsymbol{\lambda})$ to any  $\boldsymbol{x}\in B^1(\Sigma)$ and $\boldsymbol{\lambda}=(\lambda,\mu)\in H^1(\Sigma).$ The spectral curve is defined by the equation $\det \mathfrak{D}(\boldsymbol{A}(\boldsymbol{x},\boldsymbol{\lambda}))=\sum \mathcal{H}_{ij}(\boldsymbol{x})\lambda^i\mu^j=0$ and it does not depend on the choices made, if we consider the curve up to automorphisms of $H^1(\Sigma)=\mathbb{C}^\times\times\mathbb{C}^\times$. The coefficients $\mathcal{H}_{ij}$ are defined up to a transformation $\mathcal{H}_{ij}\to \mathcal{H}_{ij}\alpha^i\beta^j\gamma$ and we can use this degree of freedom to make three of them to be equal to unities. The remaining coefficients give a collection of the Poisson commuting functions w.r.t. a
natural Poisson structure on $B^2(\Sigma)$, which is maximal if the graph satisfies a certain minimality condition.


The space $B^2(\Sigma)$ possesses a canonical log-constant Poisson bracket introduced in \cite{GK} as follows. Embedding the graph $\Gamma$ into the surface $\Sigma$ induces cyclic order of ends of the edges at every vertex (a fat graph structure). Consider a surface $\tilde{\Sigma}$ corresponding to the same bipartite graph $\Gamma$, but with the cyclic order changed to the opposite in white vertices and kept unchanged in the black ones. Since the graph $\Gamma$ is embedded into $\tilde{\Sigma}$ we have the map $H^1(\tilde{\Sigma})\to B^2(\Sigma)$, which is a composition of the standard embedding with the coboundary operator. The space $H^1(\tilde{\Sigma})$ has a canonical Poisson structure, coming from the intersection index on $\tilde{\Sigma}$ and the map to $B^2(\Sigma)$ induces the Poisson bracket on the latter.

This bracket can be extended by multiplication invariance to the space $B^2(\Sigma)$ of collections of nonzero numbers attached to the faces of the graph and defines a cluster seed with the skew-symmetric exchange matrix defined by the Poisson bracket. In \cite{GK} it is observed that graphs admit elementary transformations called \textit{spider moves} such that integrable systems corresponding to them are isomorphic provided the phase spaces are related by a cluster mutation. Equivalence classes of integrable systems under such transformations are enumerated by Newton polygons of $\det \mathfrak{D}(\boldsymbol{x},\boldsymbol{\lambda})$ and the number of independent integrals of motion is just the number of integral points strictly inside these Newton polygons.

\subsection{Relations between two approaches}

We claim in sect.~\ref{ss:dimerstoloop} that the GK integrable systems coincide with the integrable systems on the Poisson-Lie loop groups $\widehat{PGL}(N)$. The isomorphism identifies not only their phase spaces and commuting flows, but also the discrete group action and the canonical cluster coordinates.

In both constructions the spectral curve of an integrable system is given by degeneracy condition of some matrix operator ($\mathcal{A}(\boldsymbol{x},\lambda)-\mu\cdot{\rm Id}$ in the group-theory approach and the Dirac operator $\mathfrak{D}(\boldsymbol{A}(\boldsymbol{x},\lambda,\mu))$ respectively). Though the matrices do not coincide, their determinants do --
roughly the correspondence goes as follows. The determinant of any matrix $\mathfrak{D}(\boldsymbol{A})$ can be written as a Grassman integral:
$$
S_\Gamma(\boldsymbol{A}) = \det \mathfrak{D}(\boldsymbol{A})=\int\exp\left(\sum \mathfrak{D}(\boldsymbol{A})^b_w\xi_b\eta^w\right)\prod_b d\xi_b\prod_w d\eta^w
$$
Therefore $\det\mathfrak{D}(\boldsymbol{A})$ can be interpreted as a partition function $S_\Gamma(\boldsymbol{A})$ of some lattice fermions in the background gauge field $\boldsymbol{A}$. Cutting torus into a cylinder corresponds to rewriting this partition function as a trace of the evolution operator from one boundary circle to another. This evolution operator is given by the matrix $\mathcal{A}(\lambda)$ acting in the external algebra of the $N$-dimensional space. Moreover cutting further this cylinder into a set smaller cylinders, one can present the evolution operator as a product of elementary steps, each depending on no more than one variable $x_f$ and exactly coinciding with elementary matrices, used to parameterise the double Bruhat cells, thus establishing the coincidence of spectral curves.

The first part of this program is establishing correspondence between words in generators of $(\hat{W}\times \hat{W})^\sharp$, enumerating cluster coordinate systems in the first approach, with the bipartite graphs, enumerating coordinates on the second approach. Moreover this correspondence should identify the letters of the word with the faces of the corresponding bipartite graphs, drawn on torus $\Sigma$, since these both sets correspond to the cluster coordinates in corresponding cases.

In order to do this we use the third combinatorial object, suggested by Dylan Thurston in unpublished paper \cite{Thurston} (and already used in \cite{GK} in our context), which we call the \textit{Thurston diagrams} and describe in detail in Appendix~\ref{ap:thurston}. A Thurston diagram is an isotopy class of a collection of curves on a surface, either closed or connecting two boundary points with only triple intersection points and such, that the connected components of the complement (faces) are colored in white and grey with any two faces sharing a segment of a curve having different colors (chessboard coloring). Such diagrams admit elementary modifications called Thurston moves. As it was already observed by D.~Thurston and A.~Henriques, every Thurston diagram defines a cluster seed (a chart on a cluster manifold) with cluster variables attached to the white faces. Thurston moves correspond to mutations (passing from one chart to another). Having Thurston diagrams on open surfaces one can glue together boundary components respecting their coloring, and thus obtain a new surface with a Thurston diagram.

In order to construct a Thurston diagram out of a reduced decomposition of an element $u\in (\widehat{W}\times \widehat{W})^\sharp$, we first associate a Thurston diagram on a cylinder with a single triple point and with $N$ grey (and white) segments on every boundary circle to every generator of $(\widehat{W}\times \widehat{W})^\sharp$ (except the cocentral one). Then we glue the cylinders together according to the order of the generators in the reduced decomposition, and finally we glue both ends of the resulting cylinder together with a twist, given by the power of the cocentral generator $\Lambda$.

In order to construct a bipartite graph out of a Thurston diagram we put a black vertex at every triple point and a white vertex at every grey face. Then we draw three edges from each black vertex inside the three grey sectors, meeting at this vertex, to the respective white vertices. It is easy to see, that the set of letters of the reduced word is in a canonical bijection with the set of white faces of the Thurston diagram and the latter are in bijection with faces of the bipartite graph. Next observation, almost as simple, is that this bijection induces a bijection between the cluster seeds, i.e. the Poisson bracket between the coordinates coincide.

Finally, we need to show that equations $\det \left(\mathcal{A}(\boldsymbol{x},\lambda)-\mu\right)=0$ and $\mathfrak{D}(\boldsymbol{A}(\boldsymbol{x},\lambda,\mu))=0$ define the same curve. 
For this purpose in sect.~\ref{s:dimers} we extend the lattice fermion partition functions on a bipartite graph to surfaces with boundary. Graphs on such surfaces are allowed to have vertices of the third type, terminating on the boundary and which can be connected to both white and black vertices, but not to each other. Denote the set of such vertices by $T$. The Dirac operator now acts as $\mathfrak{D}(\boldsymbol{A}):\mathbb{C}^{B\cup T}\to \mathbb{C}^{W\cup T}$, and for extra Grassmann variables $\boldsymbol{\zeta}=\{\zeta_t|t\in T\}$ we define
$$
S(\boldsymbol{A},\boldsymbol{\zeta})=\int\exp\left(\sum \mathfrak{D}(\boldsymbol{A})^b_w\xi_b\eta^w+\mathfrak{D}(\boldsymbol{A})^b_t\xi_b\zeta^t+
\mathfrak{D}(\boldsymbol{A})^t_w\zeta_t\eta^w\right)\prod_b d\xi_b\prod_w d\eta^w
$$
Gluing two boundary components of $\Sigma$ with a bipartite graph $\Gamma$ on it in a way, that terminal vertices are glued to terminal vertices, one gets a bipartite graph $\overline{\Gamma}$ on the glued surface $\overline{\Sigma}$. The connection $\boldsymbol{A}$ on $\Gamma$ induces a connection $\overline{\boldsymbol{A}}$ on $\overline{\Gamma}$: we just multiply the numbers of two halves of a glued edges. If $S_\Gamma(\boldsymbol{A},\boldsymbol{\zeta})$ is a partition function for $\Gamma$ then the partition function for $\overline{\Gamma}$ is given by:
$$S_{\overline{\Gamma}}(\overline{\boldsymbol{\zeta}},\overline{\boldsymbol{A}})=\int S_\Gamma(\boldsymbol{\zeta},\boldsymbol{A})e^{\sum \zeta_t\zeta_{\sigma(t)}}\prod d\zeta_t d\zeta_{\sigma(t)},$$
where the index $t$ runs over terminal edges on one side, $\sigma$ is a map sending a terminal vertex to the one it is glued to and $\overline{\boldsymbol{\lambda}}$ is $\boldsymbol{\lambda}$ with entries corresponding to glued vertices removed.

On the other hand observe that for any $N\times N$ matrix $M$ one can associate a function of $2N$ Grassmann variables $\boldsymbol{\xi}=\{\xi_i\}$ and $\boldsymbol{\eta}=\{\eta_i\}$ given by $S_M(\boldsymbol{\xi},\boldsymbol{\eta})=\exp M^i_j\xi_i\eta^j$. Matrix product corresponds to convolution of the corresponding functions
$$
S_{M_1M_2}(\boldsymbol{\xi},\boldsymbol{\eta})=\int S_{M_1}(\boldsymbol{\xi},\boldsymbol{\eta}') S_{M_2}(\boldsymbol{\xi}',\boldsymbol{\eta})e^{\sum\xi'_i\eta'_i}\prod d\xi'_id\eta'_i
$$
Thus if a partition function on a graph on a cylinder coincides with a partition function of a matrix, the partition function of several cylinders glued together corresponds to product of the matrices. Therefore in order to show the coincidence of the curves we need to cut the torus into small cylinders and verify for each of them the coincidence of partition functions. We complete the proof of the correspondence in sect.~\ref{ss:dimerstoloop}, and formulate there our main result.

The discrete flows in our integrable systems are considered in sect.~\ref{ss:automorph}.
In sect.~\ref{ss:examples} we discuss several examples of our integrable systems. In particular in  sect.~\ref{ss:pentagram} we show, that the discrete integrable system on the space of polygons in the projective plane discovered by R.~Schwartz \cite{Schwartz} can be realised as a particular case of the scheme described in the paper. More particular examples of the systems of this class have been already considered in \cite{AMJGP}.

\setcounter{equation}0
\section{Integrable systems and $r$-matrices
\label{s:isrmatr}}

Recall the standard construction of integrable systems related to the classical $r$-matrices on simple or affine Lie groups (see e.g. \cite{Sembook}), which we assume to be complex and with vanishing center. The phase spaces for these integrable systems are certain Poisson submanifolds of the Poisson-Lie groups, and the mutually commuting Hamiltonians or integrals of motion are given by the conjugation-invariant functions. Indeed, let $G$ be a Lie group, $\mathfrak{g}={\rm Lie}(G)$ be the corresponding Lie algebra and $r\in \mathfrak{g}\otimes \mathfrak{g}$ - a solution of the Yang-Baxter equation $[r_{12},r_{13}]+[r_{13},r_{23}]+[r_{12},r_{23}]=0.$
Such $r$-matrix defines a Poisson bracket on the group $G$
\begin{equation}\label{rbra}
\{ g \mathop{,}^{\otimes}g \} = -\frac{1}{2} [r,g\otimes g]
\end{equation}
and this bracket is compatible with the group structure in the sense that group multiplication $G\times G\to G$ and the inversion $G\to G$ are the Poisson maps.

It is easy to see directly from \rf{rbra}, that any two $\Ad$-invariant functions on $G$ do Poisson-commute with each other. Indeed, presenting the $r$-matrix as $r=\sum_Iv^{(1)}_I\otimes v^{(2)}_I$, where all $v_I\in \mathfrak g$, and denoting by $L_v$ (resp. $R_v$) the left (resp. right) vector field corresponding to $v$, the Poisson bracket for any two functions is given by $\{\mathcal{H}_1,\mathcal{H}_2\}=\sum_I\left(L_{v^{(1)}_I}\mathcal{H}_1L_{v^{(2)}_I}\mathcal{H}_2-
R_{v^{(1)}_I}\mathcal{H}_1R_{v^{(2)}_I}\mathcal{H}_2\right)$.
Since any $\Ad$-invariant function $\mathcal{H}$ satisfies $L_v \mathcal{H} = -R_v \mathcal{H}$,
the bracket of two such functions vanishes. Observe also, that since this argument is local, i.e. the bracket vanishes even if the functions are defined not on the whole $G$, but on any Poisson $\Ad$-invariant subvariety of $G$.

We shall restrict ourselves to the case of simple or affine Lie group $\hat{G}$ where there exists a canonical Drinfeld-Jimbo solution of the Yang-Baxter equation:
\begin{equation}\label{rDD}
r =  \sum_{\alpha\in\Delta _+} d^\alpha e_\alpha \otimes e_{\bar{\alpha}} + \frac{1}{2}\sum_{i\in \Pi} d^i h_i\otimes h_i
\end{equation}
where $\Delta_+$ is the set of positive roots, $d^\alpha= (\alpha,\alpha)/2$,  $\Pi$ is the set of positive simple roots and $\bar{\alpha}$ is just another notation for $-\alpha$, $e_\alpha$ and $h_i$ constitute the standard Cartan-Weyl basis of $\mathfrak{g}$.
To simplify the presentation we will assume in what follows that the group is simply-laced, i.e. $(\alpha,\alpha)=2$ and $d^\alpha=1$ for all roots $\alpha$.

On a simple group obviously there exists $\mbox{rank}G$ independent $\Ad$-invariant functions: a possible choice of these functions is the set $\{\mathcal{H}_i\}$, where $i \in \Pi$, $\mathcal{H}_i(g) = \Tr\ \pi_{\mu^i}(g)$ and $\pi_{\mu^i}$ be the $i$-th fundamental representation of $G$ with the highest weight $(\mu^i,\alpha_j)=\delta^i_j$ dual to $\alpha_i$, $i\in\Pi$. These function define integrable systems on Poisson submanifolds of $G$ of rank at most $2\cdot \mbox{rank}G$ (see e.g. \cite{FM}).
For a loop group $\hat{G}$, which we understand below as a group of Laurent polynomials with values in a simple group $G$, the number of independent $\Ad$-invariant function is infinite since every coefficient of $\Tr\ \pi_{\mu^i}(g)$ is now an $\Ad$-invariant function, and thus a loop group gives much larger set of integrable models.

On the space $G/\Ad H$ one can define, following \cite{HKKR}, an action of a discrete birational Poisson transformation $\tau:G/\mbox{Ad}H\to G/\mbox{Ad}H$, preserving the double Bruhat cells and the functions $\mathcal{H}_i$. Namely let $g=g_+g_-$ be the Gauss decomposition of $g \in G/\mbox{Ad}H$. Define $\tau(g)$ as the product $g_-g_+$. (The Gauss decomposition is ambiguously defined on $g$ since the Cartan part can be equally well attached to the upper or lower-triangular one or just split between the two. But on the quotient $G/\mbox{Ad}H$ the action of $\tau$ is nevertheless well defined). In sect.~\ref{ss:automorph} we show how this transformation can be generalised for the loop groups.


Recall now the classification of symplectic leaves of $G$ (see for example \cite{HKKR,ZelKog}). The group $G$ can be decomposed $G=\coprod_{u\in W\times W}G^u$ into the double Bruhat cells, enumerated by elements of the group $W\times W$, where $W$ is the Weyl group of $G$. Each double Bruhat cell is \textit{isotypic}, i.e. it is birationally equivalent to the product of symplectic manifold and a manifold with trivial Poisson bracket. The dimension of a cell $G^u$ is given by $\dim G^u = l(u)+\rank\ G$, where $l(u)$ is the length of $u$.

One can modify the Poisson manifold $G$ in order to make all constructions to be a little bit more symmetric. Namely, consider the action of the Cartan subgroup $H\in G$ on $G$ by conjugation. Since $H$ is a Poisson subgroup of $G$ with trivial Poisson structure, the quotient $G/\mbox{Ad}H$ inherits the Poisson structure and the collection of Poisson-commuting functions $\{\mathcal{H}_i\}$ as well as a decomposition into Poisson submanifolds $G/\mbox{Ad}H=\coprod_{u\in W\times W}G^u/\mbox{Ad}H$. The dimensions of the corresponding Poisson submanifolds are now just $\dim G^u/\mbox{Ad}H= l(u)$.

\section{Cluster parametrisation of double Bruhat cells. \\ Simple groups
\label{ss:clustcell}}

Following \cite{Drinfeld} we describe here how to introduce the structure of a cluster variety on $G^u$ and $G^u/\mbox{Ad}H$. Namely, starting from a decomposition of $u$ into reduced product $s_{i_1}\cdots s_{i_l}$ of standard generators of $W\times W$ we define a cluster seed (see appendix~\ref{ap:cluster}) - a split algebraic torus, provided with log-constant Poisson structure, and a Zariski open its Poisson embedding into $G^u$. Similarly, a cluster seed for the space $G^u/\Ad H$ is constructed from the same data. Finally, we show that for the group $G=PGL(N)$ the cluster seeds are isomorphic to those, corresponding to the Thurston diagrams constructed out of the decomposition $u=s_{i_1}\cdots s_{i_l}$ on a disk.

\subsection{Cartan-Weyl generators of a simple group}

To describe cluster coordinates on $G^u$ and $G^u/\Ad H$ we need to introduce first a set of generators of the Lie group $G$ analogous to the Cartan-Weyl generators of the corresponding Lie algebra $\mathfrak{g}$.
The set of standard generators of a Weyl group $W$ is in canonical bijection with the set $\Pi$ of simple roots of $G$, and we shall not distinguish between these two sets. The set of generators of the second copy of $W$ will be identified with the set of negative simple roots ${\bar\Pi}$.

Recall, that given a Cartan matrix $C_{ij}$ the corresponding Lie algebra $\mathfrak{g}$ is generated by $\{h_i|i \in \Pi\}$ and $\{e_i|i \in \Pi\cup{\bar\Pi}\}$, satisfying the relations\footnote{For simplicity we denote by $\bar{i}$ the root, opposite to the root $i$, and extend $h$ and $C$ to the negative roots, assuming that $h_i = h_{\bar{i}}$, that $C_{\bar{i},\bar{j}}=C_{ij}$ and that $C_{ij}=0$ if $i$ and $j$ have different signs.}
\begin{equation}\label{rela}
\begin{array}{l}
\ [h_i,h_j]=0,\\
\ [h_i,e_{j}] = \sign(j)C_{ij} e_j,\\
\ [e_i,e_{\bar{i}}] =\sign(i) h_i,\\
\ (\mbox{ad}\, e_i)^{1-C_{ij}}e_j = 0 \mbox{ for } i+j\neq 0
\end{array}
  \end{equation}
For the same Lie algebra one can replace the set $\{h_i\}$ by the set $\{h^i\}$, defined by $h_i=\sum_{j\in\Pi}C_{ij}h^j$, then the relations \rf{rela} take the form:
\begin{equation}
\label{relx}
\begin{array}{l}
\ [h^i,h^j]=0,\\
\ [h^i,e_{j}] = \sign(j)\delta_i^j e_j,\\
\ [e_i,e_{\bar{i}}] =\sign(i) C_{ij}h^j,\\
\ (\mbox{ad}\, e_i)^{1-C_{ij}}e_j = 0 \mbox{ for } i+j\neq 0
\end{array}
\end{equation}
For any $i \in \Pi\cup{\bar\Pi}$ introduce the group element $E_i = \exp(e_i)$ and a one-parameter subgroup $H_i(x)=\exp(h^{i}\log x)$ which will be our set of generators of the group $G$.
The commutational relations (\ref{relx}) imply the relations between $H_i$ and $E_i$ (we list them in  appendix~\ref{ap:relx} for the simply-laced case).
As an immediate consequence of (\ref{relx}) notice, that $H_i(x)$ commutes with $E_j$, unless $i=j$, and that $E_i$ commutes with $E_j$ if $C_{ij}=0$ and $i\neq \bar{j}$.

\subsection{Construction of the cluster seeds}

Take a decomposition $u=s_{i_1}\cdots s_{i_{l}}$ of an element $u\in W\times W$.  A seed $((\mathbb{C}^\times)^l,\varepsilon)$ and a map
$ev :(\mathbb{C}^\times)^{l}\rightarrow G/\Ad\ H$ is associated to such $u$ by
\begin{equation}
\label{paramG/AdH}
x_1,\ldots,x_l \mapsto H_{i_1}(x_1)E_{i_1}\cdots H_{i_l}(x_l)E_{i_l}
\end{equation}
The image of this map is Zariski open in the double Bruhat cell $C^u$, and it is an embedding, if the word $u$ is reduced. Different reduced decompositions of the same element $u\in W\times W$ give rise to different parametrisations, related by a cluster transformation.

To construct a seed $((\mathbb{C}^\times)^{(l+r)},\varepsilon)$, where $r$ is the rank of $G$, i.e. to parametrise the cells of simple group $G$ itself, one needs just to multiply this expression from the right by an arbitrary element of the Cartan subgroup:
\begin{equation}\label{paramG}
(x_1,\ldots,x_{l+r}) \mapsto H_{i_1}(x_1)E_{i_1}\cdots H_{i_l}(x_l)E_{i_l}H_{1}(x_{l+1})\cdots H_{r}(x_{l+r})
\end{equation}
The construction of the corresponding exchange matrix $\varepsilon$ is given in the appendix~\ref{ap:exchange}.

Using the fact, that $H_i(x)E_j=E_jH_i(x)$ unless $i=j$, one can rewrite expressions (\ref{paramG/AdH}) and (\ref{paramG}) in many different ways by moving any $H_i$ until it meets $E_i$ or $E_{\bar{i}}$. Therefore every cluster variable is naturally associated to a positive simple root $i$, and for a given $i$ to a minimal segment of the word $s_{i_1}\cdots s_{i_{l}}$ delimited by $s_i$, $s_{\bar{i}}$, or ends of the word (for a cyclic word all segments are delimited just by $s_i$ or $s_{\bar{i}}$).

One can check (see \cite{Drinfeld}) that cluster seeds, corresponding to different decompositions of the same word $u$, are related by a cluster transformation. In particular, their images coincide up to codimension one. For example, if $C_{ij}=C_{ji}=-1$, the coordinates corresponding to decomposition $u=As_is_js_iB$, where $A$ and $B$ are arbitrary words, are related to the coordinates corresponding to decomposition $\tilde{u}=As_js_is_jB$ by a mutation in the variable, associated to the segment $[s_is_js_i]$. Similarly, for the relation between $As_is_{\bar{i}}B \leftrightarrow A s_{\bar{i}}s_iB$ we should make a mutation in the variable associated to the segment $[s_is_{\bar{i}}]$.  If $C_{ij}=0$ applying relation for $As_is_jB\leftrightarrow As_js_iB$ does not change parametrisation, since the corresponding matrices $E_i$ and $E_j$ commute.

If a decomposition is not reduced, the maps (\ref{paramG/AdH}) and (\ref{paramG}) are still defined, but they are not embeddings any longer. Instead such map is an embedding, corresponding to a reduced word, pre-composed with a projection along some coordinates and some mutations. Indeed, the map corresponding to the word $As_is_iB$ is a composition of the mutation in the variable associated to the segment $[s_is_i]$ with the map corresponding to $As_iB$.

Therefore we claim, that the double Bruhat cells of $G$ and of $G/\mbox{Ad}H$ are cluster varieties. In \cite{Drinfeld} it is proven that the Poisson brackets on $G^u$ (and thus on $G^u/\mbox{Ad}H$), given by the exchange matrix $\varepsilon$, and \rf{rbra} given by the Drinfeld-Jimbo $r$-matrix \rf{rDD} coincide.

\subsection{Generators and Thurston diagrams for the group $PGL(N)$}

In the case $G=PGL(N)$ we can make the construction much more explicit.
The generators $E_i$ and $H_i(x)$ in the standard representation have a particularly simple form:
\begin{equation}
\renewcommand{\arraycolsep}{1pt}
\label{genslN}
H_i(x) =
\begin{pmatrix}
x&0&&\cdots&&0\\[-5pt]
0&\ddots&&&&0\\[-5pt]
&&x&&\\[-5pt]
\vdots&&&1&\\[-5pt]
&&&&\ddots&0\\[-3pt]
0&&\cdots&&0&1
\end{pmatrix}\quad\quad
E_i = E_{\bar i}^{\rm tr} =
\begin{pmatrix}
1&0&&\cdots&&0\\[-5pt]
0&\ddots&&&&0\\[-5pt]
&&1&1&\\[-5pt]
\vdots&&&1&\\[-5pt]
&&&&\ddots&0\\[-3pt]
0&&\cdots&&0&1
\end{pmatrix}
\end{equation}
where the lowest line of $H_i(x)$ with $x$, and the line in $E_i$ containing the off-diagonal unity have number $i>0$.
For negative $i$ the corresponding matrix $E_{-i} = E_{\bar i}$ is just transposed to the matrix of the positive root.


Let us now give an alternative description of the cluster seeds, corresponding to the decomposition of $u\in W\times W$ into product of generators for the group $PGL(N)$, using the Thurston diagrams
(see appendix~\ref{ap:thurston}). Every decomposition of $u$ into the product of generators corresponds to a Thurston diagram, and the latter in its turn corresponds to a cluster seed. We claim that this alternative way gives the same seed. To verify this statement one needs, first, to compare the seeds corresponding to a single generator. This can be done by comparing the exchange graphs for such elementary Thurston diagrams containing only one triple point (see fig.\ref{fi:thurston-example}A) and the exchange graph (chord) described in appendix \ref{ap:exchange}. Then we need to verify that gluing exchange graphs corresponds to gluing Thurston diagrams. We leave these verifications as easy exercises.

For the group $G=PGL(N)$ every word $u\in W\times W$ in the generators corresponds to a Thurston diagram on a disc, which we draw as an infinite vertical strip. Every diagram consists of $2N$ curves connecting the sides of the strip and having no vertical tangent. The curves are oriented in such a way, that for generic vertical section the orientations of the curves at the intersection points with the section alternate. Generators $s_i$ (and $s_{\bar{i}}$) correspond to a triple intersection of the curves with the numbers $2i-1,2i,2i+1$ (and $2i,2i+1,2i+3$ respectively), counted from above along a section. An example of a Thurston diagram is shown on fig.~\ref{fi:thurston-example}. The properties of the correspondence between the diagrams and the words are listed in appendix~\ref{ap:thurston}. The face variables, corresponding
to the top and bottom white faces, for a simple group $G=PGL(N)$ are restricted to be equal to unity.

\begin{figure}[ht]
\begin{center}
\begin{tikzpicture}
\draw[fill=gray!20] (0,0.4)..controls(0.3,0.4)..(0.6,0)--(0,0)--cycle;
\draw[fill=gray!20] (1.2,0.4)..controls(0.9,0.4)..(0.6,0)--(1.2,0)--cycle;
\draw[fill=gray!20] (0,-0.8)--(0,-0.4)..controls(0.3,-0.4)..(0.6,0)..controls(0.9,-0.4)..(1.2,-0.4)--(1.2,-0.8)--cycle;
\draw[fill=gray!20] (0,0.8)--(1.2,0.8)--(1.2,1.2)--(0,1.2)--cycle;
\draw[fill=gray!20] (0,-1.2)--(1.2,-1.2)--(1.2,-1.6)--(0,-1.6)--cycle;
\draw[thick] (0,-2)--(0,1.6);
\draw[thick] (1.2,-2)--(1.2,1.6);
\draw[-stealth] (0.55,0.2)--(0.3,-0.1);
\draw[-stealth] (0.4,-0.15)--(0.8,-0.15);
\draw[-stealth] (0.9,-0.1) -- (0.65,0.2);
\foreach \y in {-2.4,-1.6,-0.8,0}
  {
\begin{scope}[yshift=\y cm]
\draw[-stealth,color=gray!60] (0.1,0.65)--(0.1,1.35);
\draw[stealth-,color=gray!60] (1.1,0.65)--(1.1,1.35);
\end{scope}
  }
\draw (2.1,-0.2) node {$\rightarrow$};
\draw[semithick] (3,-1) -- (4.2,-1);
\draw[semithick] (3,-0.2) -- (4.2,-0.2);
\draw[semithick] (3,0.6) -- (4.2,0.6);
\draw[semithick] (3,1.4) -- (4.2,1.4);
\draw[semithick] (3,-1.8) -- (4.2,-1.8);
\draw[ultra thick,color=gray!55,shorten >=0.3cm,-stealth] (3.6,0.6)--(3.2,-0.2);\draw[ultra thick,color=gray!55] (3.6,0.6)--(3.2,-0.2);
\draw[ultra thick,color=gray!55,shorten >=0.35cm,-stealth] (4,-0.2)--(3.6,0.6);\draw[ultra thick,color=gray!55] (4,-0.2)--(3.6,0.6);
\draw[ultra thick,color=gray!55,shorten >=0.3cm,-stealth] (3.6,-1)--(3.2,-0.2);\draw[ultra thick,color=gray!55] (3.6,-1)--(3.2,-0.2);
\draw[ultra thick,color=gray!55,shorten >=0.35cm,-stealth] (4,-0.2)--(3.6,-1);\draw[ultra thick,color=gray!55] (4,-0.2)--(3.6,-1);
\draw[ultra thick,shorten >=0.25cm,-stealth] (3.2,-0.2)--(4,-0.2);\draw[ultra thick] (4,-0.2)--(3.2,-0.2);
\fill (3.6,-1) circle (0.08);
\fill (3.2,-0.2) circle (0.08);\fill (4,-0.2) circle (0.08);
\fill (3.6,0.6) circle (0.08);
\node at (2,-2.5) {A};
\end{tikzpicture}
\qquad\qquad\qquad\qquad
\begin{tikzpicture}
\fill[color=gray!10] (-0.3,0)--(1.5,0)..controls(1.8,0)..(2.1,0.4)--(-0.3,0.4)--cycle;
\draw (-0.3,0)--(1.5,0)..controls(1.8,0)..(2.1,0.4)--(-0.3,0.4)--cycle;
\fill[color=gray!10] (-0.3,0.8)--(0.3,0.8)..controls (0.6,0.8)..(0.9,1.2)..controls (0.6,1.2)..(0.3,1.6)..controls(0,1.2)..(-0.3,1.2)--cycle;
\draw (-0.3,0.8)--(0.3,0.8)..controls (0.6,0.8)..(0.9,1.2)..controls (0.6,1.2)..(0.3,1.6)..controls(0,1.2)..(-0.3,1.2)--cycle;
\fill[color=gray!10] (-0.3,1.6)--(0.3,1.6)..controls(0,2)..(-0.3,2)--cycle;
\draw (-0.3,1.6)--(0.3,1.6)..controls(0,2)..(-0.3,2)--cycle;
\fill[color=gray!10] (0.3,1.6)..controls(0.6,2)..(0.9,2)..controls(1.2,2)..(1.5,1.6)..controls(1.2,1.6)..(0.9,1.2)..controls(0.6,1.6)..(0.3,1.6)--cycle;
\draw (0.3,1.6)..controls(0.6,2)..(0.9,2)..controls(1.2,2)..(1.5,1.6)..controls(1.2,1.6)..(0.9,1.2)..controls(0.6,1.6)..(0.3,1.6)--cycle;
\fill[color=gray!10] (3.3,1.6)--(1.5,1.6)..controls(1.8,2)..(2.1,2)--(3.3,2)--cycle;
\draw (3.3,1.6)--(1.5,1.6)..controls(1.8,2)..(2.1,2)--(3.3,2)--cycle;
\fill[color=gray!10] (2.7,0.8)..controls(2.4,1.2)..(2.1,1.2)..controls(1.8,1.2)..(1.5,1.6)..controls(1.2,1.2)..(0.9,1.2)..controls(1.2,0.8)..(1.5,0.8)..controls(1.8,0.8)..(2.1,0.4)..controls(2.4,0.8)..(2.7,0.8);
\draw (2.7,0.8)..controls(2.4,1.2)..(2.1,1.2)..controls(1.8,1.2)..(1.5,1.6)..controls(1.2,1.2)..(0.9,1.2)..controls(1.2,0.8)..(1.5,0.8)..controls(1.8,0.8)..(2.1,0.4)..controls(2.4,0.8)..(2.7,0.8);
\fill[color=gray!10] (3.3,0.4)..controls(3,0.4)..(2.7,0.8)..controls(2.4,0.4)..(2.1,0.4)..controls(2.4,0)..(2.7,0)--(3.3,0)--cycle;
\draw (3.3,0.4)..controls(3,0.4)..(2.7,0.8)..controls(2.4,0.4)..(2.1,0.4)..controls(2.4,0)..(2.7,0)--(3.3,0)--cycle;
\fill[color=gray!10] (3.3,1.2)..controls(3,1.2)..(2.7,0.8)--(3.3,0.8)--cycle;
\draw (3.3,1.2)..controls(3,1.2)..(2.7,0.8)--(3.3,0.8)--cycle;
\draw[thick,color=gray!20] (-0.3,-0.2)--(-0.3,2.2);
\draw[thick,color=gray!20] (3.3,-0.2)--(3.3,2.2);
\node at (0.3,-0.5) {$s_1$};
\node at (0.9,-0.5) {$s_{\bar{1}}$};
\node at (1.5,-0.5) {$s_1$};
\node at (2.1,-0.5) {$s_{\bar{2}}$};
\node at (2.7,-0.5) {$s_2$};
\node at (1.5,-1) {B};
\end{tikzpicture}
\end{center}
\caption{Thurston diagrams: (A) for an elementary generator $s_i$ of the Weyl group; (B) for the word $s_1s_{\bar{1}}s_1s_{\bar{2}}s_2$.\label{fi:thurston-example}}
\end{figure}
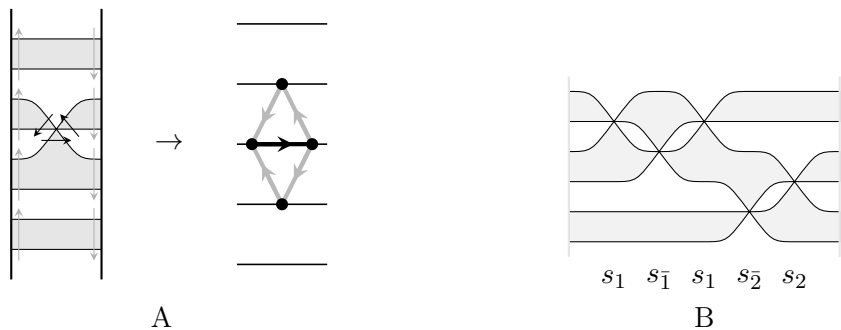
The Thurston diagrams for $G^u/\mbox{Ad}H$ are obtained from those for $G^u$ by gluing together the right and the left sides, thus getting diagrams on a cylinder instead of a strip.

\subsection{Example: Poisson submanifolds of $PGL(3)$}

\begin{figure}[H]
\begin{center}
\begin{tikzpicture}
\draw[black, fill=gray!20] (0,2)..controls(0.3,2)..(0.6,1.6)--(0,1.6);
\draw[black,fill=gray!20] (0.6,1.6)..controls(0.9,1.6)..(1.2,1.2)..controls(1.5,1.6)..(1.8,1.6)--(3,1.6)..controls(2.7,2)..(2.4,2)--(1.2,2)..controls(0.9,2)..(0.6,1.6);
\draw[black,fill=gray!20] (4.2,1.6)..controls(3.9,1.6)..(3.6,1.2)..controls(3.3,1.6)..(3,1.6)..controls(3.3,2)..(3.6,2)--(4.2,2);
\draw[black,fill=gray!20] (0,1.2)..controls(0.3,1.2)..(0.6,1.6)..controls(0.9,1.2)..(1.2,1.2)..controls(0.9,0.8)..(0.6,0.8)--(0,0.8);
\draw[black,fill=gray!20] (1.2,1.2)..controls(1.5,1.2)..(1.8,0.8)..controls(1.5,0.8)..(1.2,1.2);
\draw[black,fill=gray!20] (1.8,0.8)..controls(2.1,1.2)..(2.4,1.2)..controls(2.7,1.2)..(3,1.6)..controls(3.3,1.2)..(3.6,1.2)..controls(3.3,0.8)..(3,0.8)..controls(2.7,0.8)..(2.4,0.4)..controls(2.1,0.8)..(1.8,0.8);
\draw[black,fill=gray!20] (4.2,1.2)--(3.6,1.2)..controls(3.9,0.8)..(4.2,0.8);
\draw[black,fill=gray!20] (0,0.4)--(1.2,0.4)..controls(1.5,0.4)..(1.8,0.8)..controls(2.1,0.4)..(2.4,0.4)..controls(2.1,0)..(1.8,0)--(0,0);
\draw[black,fill=gray!20](4.2,0.4)--(2.4,0.4)..controls(2.7,0)..(3,0)--(4.2,0);
\draw[thin] (0,-0.2)--(0,2.2);
\draw[thin] (4.2,-0.2)--(4.2,2.2);
\node at (0.2,1.4) {$x_1$};
\node at (0.9,1.4) {$x_2$};
\node at (1.2,0.6) {$x_3$};
\node at (2.1,0.6) {$x_4$};
\node at (1.9,1.3) {$x_5$};
\node at (3.3,1.4) {$x_6$};
\node at (4.0,1.35) {$x_7$};
\node at (3.6,0.8) {$x_8$};
\node at (0.6,-0.5) {$s_1$};
\node at (1.2,-0.5) {$s_{\bar{1}}$};
\node at (1.8,-0.5) {$s_2$};
\node at (2.4,-0.5) {$s_{\bar{2}}$};
\node at (3,-0.5) {$s_1$};
\node at (3.6,-0.5) {$s_{\bar{1}}$};
\node at (2,-1) {A};
\end{tikzpicture}
\qquad
\begin{tikzpicture}
\draw[black, fill=gray!20] (0,2)..controls(0.3,2)..(0.6,1.6)--(0,1.6);
\draw[black,fill=gray!20] (0.6,1.6)..controls(0.9,1.6)..(1.2,1.2)..controls(1.5,1.6)..(1.8,1.6)--(3,1.6)..controls(2.7,2)..(2.4,2)--(1.2,2)..controls(0.9,2)..(0.6,1.6);
\draw[black,fill=gray!20] (3,1.6)--(1.8,1.6)..controls(1.5,1.6)..(1.2,1.2)..controls(0.9,1.6)..(0.6,1.6)..controls(0.9,2)..(1.2,2)--(3,2);
\draw[black,fill=gray!20] (0,1.2)..controls(0.3,1.2)..(0.6,1.6)..controls(0.9,1.2)..(1.2,1.2)..controls(0.9,0.8)..(0.6,0.8)--(0,0.8);
\draw[black,fill=gray!20] (1.2,1.2)..controls(1.5,1.2)..(1.8,0.8)..controls(1.5,0.8)..(1.2,1.2);
\draw[black,fill=gray!20] (3,1.2)--(2.4,1.2)..controls(2.1,1.2)..(1.8,0.8)..controls(2.1,0.8)..(2.4,0.4)..controls(2.7,0.8)..(3,0.8);
\draw[black,fill=gray!20] (0,0.4)--(1.2,0.4)..controls(1.5,0.4)..(1.8,0.8)..controls(2.1,0.4)..(2.4,0.4)..controls(2.1,0)..(1.8,0)--(0,0);
\draw[black,fill=gray!20](3,0.4)--(2.4,0.4)..controls(2.7,0)..(3,0);
\node at (0.2,1.4) {$x_1$};
\node at (0.9,1.4) {$x_2$};
\node at (1.2,0.6) {$x_3$};
\node at (2.1,0.6) {$x_4$};
\node at (0.6,-0.5) {$s_1$};
\node at (1.2,-0.5) {$s_{\bar{1}}$};
\node at (1.8,-0.5) {$s_2$};
\node at (2.4,-0.5) {$s_{\bar{2}}$};
\node at (1.5,-1) {B};
\end{tikzpicture}
\qquad
\begin{tikzpicture}
\draw[black,fill=gray!20](0,2)--(0.6,2)..controls(0.9,2)..(1.2,1.6)..controls(0.9,1.6)..(0.6,1.2)..controls(0.3,1.6)..(0,1.6);
\draw[black,fill=gray!20](0,1.2)--(0.6,1.2)..controls(0.3,0.8)..(0,0.8);
\draw[black,fill=gray!20](0,0.4)--(1.2,0.4)..controls(0.9,0)..(0.6,0)--(0,0);
\draw[black,fill=gray!20](0.6,1.2)--(1.2,0.4)..controls(1.5,0.8)..(1.8,0.8)--(1.2,1.6)..controls(0.9,1.2)..(0.6,1.2);
\draw[black,fill=gray!20](2.4,2)--(1.8,2)..controls(1.5,2)..(1.2,1.6)--(2.4,1.6);
\draw[black,fill=gray!20](2.4,1.2)..controls(2.1,1.2)..(1.8,0.8)--(2.4,0.8);
\draw[black,fill=gray!20](2.4,0.4)..controls(2.1,0.4)..(1.8,0.8)..controls(1.5,0.4)..(1.2,0.4)..controls(1.5,0)..(1.8,0)--(2.4,0);
\node at (0.2,1.4) {$x'_1$};
\node at (0.9,1.4) {$x'_2$};
\node at (0.6,0.7) {$x'_3$};
\node at (1.5,0.6) {$x'_4$};
\node at (0.6,-0.5) {$s_{\bar{1}}$};
\node at (1.2,-0.5) {$s_1s_{\bar{2}}$};
\node at (1.8,-0.5) {$s_2$};
\node at (1,-1) {C};
\end{tikzpicture}
\end{center}
\caption{Thurston diagrams for: (A) the word $1\bar{1}2\bar{2}1\bar{1}$; (B) the cyclic word $1\bar{1}2\bar{2}$; (C) the cyclic word $\bar{1}1\bar{2}2$ (for B and C the right and the left sides are glued together).\label{fi:112211}}
\end{figure}
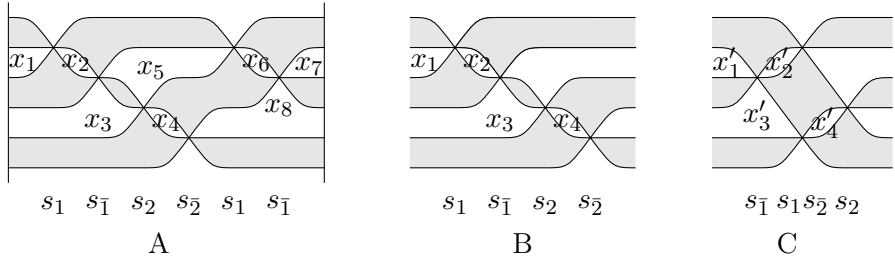

In the case of $G=PGL(3)$ there are just two simple roots, and the Cartan matrix is
\begin{equation}
\label{Csl3}
C_{ij} = 2\frac{(\alpha_i,\alpha_j)}{(\alpha_i,\alpha_i)}=(\alpha_i,\alpha_j)=\left(
\begin{array}{cc}
  2 & -1 \\
  -1 & 2
\end{array}\right)
\end{equation}
The elementary matrices \rf{genslN} are in this case
\be
\label{pgl3}
E_1 = E_{\bar{1}}^{\rm tr}=
\begin{pmatrix}
1&1&0\\[-2pt]
0&1&0\\[-2pt]
0&0&1
\end{pmatrix};\ \ \ E_2 = E_{\bar{2}}^{\rm tr}=
\begin{pmatrix}
1&0&0\\[-2pt]0&1&1\\[-2pt]0&0&1
\end{pmatrix}
\\
 H_1(x)=
\begin{pmatrix}
x&0&0\\[-2pt]
0&1&0\\[-2pt]
0&0&1
\end{pmatrix};\ \ \
H_2(x)=
\begin{pmatrix}
x&0&0\\[-2pt]
0&x&0\\[-2pt]
0&0&1
\end{pmatrix}
\ee
The big cell in $G=PGL(3)$ is parametrised by particular case of expression \rf{paramG}, corresponding to a decomposition of the longest element of $W\times W$ (here of the length $l=6$). For an element $1\bar{1}2\bar{2}1\bar{1}$ with the Thurston diagram, presented at fig.~\ref{fi:112211}A, the corresponding product is

\be
\label{bigcell}
(x_1,\ldots,x_8)\mapsto\\
\mapsto H_1(x_1)E_1 H_1(x_2)E_{\bar{1}} H_2(x_3)E_2 H_2(x_4)E_{\bar{2}} H_1(x_5)E_1 H_1(x_6)E_{\bar{1}} H_1(x_7) H_2(x_8).
\ee
The coordinates $x_1,\ldots,x_8$ are associated with the segments $]1]$, $[1\bar{1}]$, $]1\bar{1}2]$, $[2\bar{2}]$, $[\bar{2}1\bar{1}[$, $[1\bar{1}]$, $[\bar{1}[$, $[\bar{2}1\bar{1}[$, respectively.

For the Poisson submanifold of dimension $l=2\cdot\mbox{rank}G=4$  in $G/\mbox{Ad}H$, corresponding to the longest cyclically irreducible word $1\bar{1}2\bar{2}$ with the Thurston diagram, presented at fig.~\ref{fi:112211}B, the parameterisation \rf{paramG} gives
\be
\label{pretoda}
(x_1,x_2,x_3,x_4) \mapsto H_1(x_1)E_1 H_1(x_2)E_{\bar{1}} H_1(x_3)E_2 H_1(x_4)E_{\bar{2}}=
\\
=\begin{pmatrix}1&1&0\\0&1&0\\0&0&1\end{pmatrix}
\begin{pmatrix}x_1&0&0\\0&1&0\\0&0&1\end{pmatrix}
\begin{pmatrix}1&0&0\\1&1&0\\0&0&1\end{pmatrix}
\begin{pmatrix}x_2&0&0\\0&1&0\\0&0&1\end{pmatrix}\times
\\
\times
\begin{pmatrix}1&0&0\\0&1&1\\0&0&1\end{pmatrix}
\begin{pmatrix}x_3&0&0\\0&x_3&0\\0&0&1\end{pmatrix}
\begin{pmatrix}1&0&0\\0&1&0\\0&1&1\end{pmatrix}
\begin{pmatrix}x_4&0&0\\0&x_4&0\\0&0&1\end{pmatrix}=
\\
=\begin{pmatrix}x_1x_2x_3x_4+x_2x_3x_4&x_3x_4+x_4&1\\x_2x_3x_4&x_3x_4+x_4&1\\0&x_4&1\end{pmatrix}
\ee
where the matrix in the r.h.s. should be understood as an element of $PGL(3)$, i.e.
modulo multiplication by a constant.

The parametrisation for the word $\bar{1}1\bar{2}2$ (the corresponding Thurston diagram is presented at fig.~\ref{fi:112211}C) is given by
$$
(x'_1,x'_2,x'_3,x'_4) \mapsto H_1(x'_1)E_{\bar{1}} H_1(x'_2)E_1 H_1(x'_3)E_{\bar{2}} H_1(x'_4)E_2=
$$
$$
=\begin{pmatrix}x'_1x'_2x'_3x'_4&x'_1x'_3x'_4&x'_1x'_3\\
x'_1x'_2x'_3x'_4&x'_1x'_3x'_4+x'_3x'_4&x'_1x'_3+x'_3\\
0&x'_3x'_4&x'_3+1\end{pmatrix}
$$
where the coordinates $x'_1,\ldots,x'_4$ are related to the coordinates $x_1,\ldots,x_4$ by mutations in the variables $x_2$ and $x_4$:
$$
(x'_1,x'_2,x'_3,x'_4)=(x_1(1+x_2)(1+1/x_4)^{-1},1/x_2,x_3(1+1/x_2)^{-1}(1+x_4),1/x_4)
$$
The Poisson structures, corresponding to the symplectic leaves in $PGL(3)$, can be described
by exchange graphs, see appendix~\ref{ap:exchange}. For the big cell in $PGL(3)$, corresponding to the word
$u=1\bar{1}2\bar{2}1\bar{1}\in W\times W$, the corresponding graph is presented at fig.~\ref{fi:sl3-simple}, while for the Poisson submanifold in $G/\mbox{Ad}H$, corresponding to $u=1\bar{1}2\bar{2}$, the
exchange graph can be found at fig.~\ref{fi:todasl3}.

\section{Cluster parameterisation of double Bruhat cells.\\ Loop groups
\label{s:loop}}

Let us generalise this construction to loop groups, namely present it for $\widehat{PGL}(N)$. For loop groups one gets infinitely many  $\Ad$-invariant functions, and they possess finite-dimensional Poisson submanifolds, thus allowing to construct much wider class of integrable systems.

The main difference with the case of simple groups, since the Cartan matrix for affine groups is non-invertible, is that relations (\ref{rela}) and (\ref{relx}) define now non-isomorphic Lie algebras. The former defines a centrally extended loop group $\widehat{PGL}_\#(N)=\widehat{PGL}(N)\ltimes \mathbb{C}^\times$, while the latter corresponds to the co-extended one $\widehat{PGL}^\sharp(N)=\widehat{PGL}(N)\rtimes \mathbb{C}^\times$, see e.g. \cite{Kac}. For our purposes we shall use the group $\widehat{PGL}^\sharp(N)$ since it admits cluster parameterisation.
The simple roots for this group can be identified with the set $\Pi=\mathbb{Z}/N\mathbb{Z}$, with the Dynkin diagram given by a closed necklace with $N$ vertices.

\subsection{The coextended affine Weyl group, wiring and Thurston diagrams\label{ss:coextended-Weyl}}

The Weyl group $W$ of the group $\widehat{PGL}^\sharp(N)$ also admits a central co-extension $W^\sharp = W\rtimes\mathbb{Z}/N\mathbb{Z}$, and can be defined by generators, corresponding to the simple roots $\{s_i|i\in \mathbb{Z}/N\mathbb{Z}\}$, and an additional generator $\Lambda$ with relations
\be\label{Wrelations}
  s_i s_{i+1} s_i = s_{i+1}s_{i}s_{i+1},\\
  \Lambda s_i=s_{i+1}\Lambda,\\
  s_i^2=1,\ \ \ \ \Lambda^N=1
\ee
Similarly, the group $(\widehat{W}\times \widehat{W})^\sharp=(\widehat{W}\times \widehat{W})\rtimes(\mathbb{Z}/N\mathbb{Z})$ is generated by $\{s_i|i\in \Pi\cup \bar{\Pi}\}$, corresponding now to positive or negative simple roots, and $\Lambda$ subject to (\ref{Wrelations}), with one additional relation
\begin{equation}\label{WWrelations}
s_is_j=s_js_i \mbox{ if $i>0$ and $j<0$}
\end{equation}
The subgroup of $(\widehat{W}\times \widehat{W})^\sharp$, generated by $s_is_{\overline{i}}$ and $\Lambda$, is isomorphic to $\widehat{W}^\sharp$ by the obvious isomorphism $s_is_{\overline{i}}\mapsto s_i$. It will be called the \textit{diagonal subgroup} and denoted by $\widehat{W}^\sharp_+$.

Elements of the group $\widehat{W}^\sharp$ as well as their decomposition into products of generators can be visualised by \textit{wiring diagrams} similarly to the finite Weyl groups of type $A_{N-1}$, see e.g \cite{FZ}. The only difference is that for the affine case the diagrams are drawn on cylinders instead of strips.
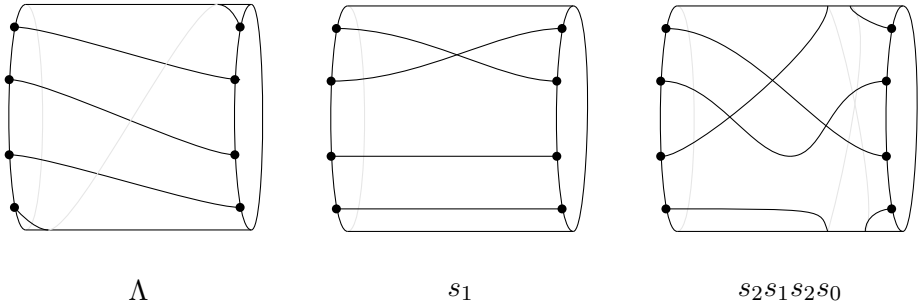
\begin{figure}[ht]
\begin{center}
\begin{tikzpicture}
\draw[color=gray!20] (0,0) ..controls (0.3,0) and (0.3,3).. (0,3);
\draw (0,0)--(3,0); \draw (0,3)--(3,3);
\draw (3,3) ..controls (2.7,3) and (2.7,0).. (3,0);
\draw (0,3) ..controls (-0.3,3) and (-0.3,0).. (0,0);
\draw (3,0) ..controls (3.2,0) and (3.2,3).. (3,3);
\fill (2.85,2.7) circle(0.06);
\fill (2.78,2) circle(0.06);
\fill (2.78,1) circle(0.06);
\fill (2.85,0.3) circle(0.06);
\fill (-0.15,2.7) circle(0.06);
\fill (-0.22,2) circle(0.06);
\fill (-0.22,1) circle(0.06);
\fill (-0.15,0.3) circle(0.06);
\draw  (-0.15,2.7)..controls(0.35,2.7) and(2.35,2)..(2.85,2) ;
\draw (-0.22,2)..controls(0.28,2) and (2.28,1)..(2.78,1);
\draw (-0.22,1)..controls(0.28,1)and(2.35,0.3)..(2.85,0.3);
\draw (2.85,2.7)..controls (2.85,2.7)and(2.75,3)..(2.55,3);
\draw[color=gray!20] (2.55,3)..controls(2.15,3)and(0.7,0)..(0.3,0);
\draw (0.3,0)..controls(0.1,0)and(-0.15,0.3)..(-0.15,0.3);
\node at (1.5, -0.8) {$\Lambda$};
\end{tikzpicture}\qquad
\begin{tikzpicture}[rotate=270]
\draw[color=gray!20] (0,0) ..controls (0,0.3) and (3,0.3).. (3,0);
\draw (0,0)--(0,3); \draw (3,0)--(3,3);
\draw (3,3) ..controls (3,2.7) and (0,2.7).. (0,3);
\draw (3,0) ..controls (3,-0.3) and (0,-0.3).. (0,0);
\draw (0,3) ..controls (0,3.3) and (3,3.2).. (3,3);
\fill (2.7,2.85) circle(0.06);
\fill (2,2.78) circle(0.06);
\fill (1,2.78) circle(0.06);
\fill (0.3,2.85) circle(0.06);
\fill (2.7,-0.15) circle(0.06);
\fill (2,-0.22) circle(0.06);
\fill (1,-0.22) circle(0.06);
\fill (0.3,-0.15) circle(0.06);
\draw (2.7,-0.15) -- (2.7,2.85);
\draw (2,2.78) -- (2,-0.22);
\draw (0.3,2.85) ..controls (0.3,2) and (1,1).. (1,-0.22);
\draw (1,2.78) ..controls (1,2) and (0.3,1).. (0.3,-0.15);
\node at (3.8,1.5) {$s_1$};
\end{tikzpicture}\qquad
\begin{tikzpicture}[rotate=270]
\draw[color=gray!20] (0,2.3) ..controls (0.2,2.7) and (2.7,2.1)..(3,2);
\draw[color=gray!20] (0,0) ..controls (0,0.3) and (3,0.3).. (3,0);
\draw[color=gray!20] (3,2.5) .. controls (2.7,2.7) and (0.3,2.2).. (0,2);
\draw (1,2.8) ..controls (1,2) and (2,2).. (2,1.5);
\draw (2,1.5) ..controls (2,1) and (1,0.5).. (1,-0.2);
\draw (0,3) ..controls (0,3.3) and (3,3.2).. (3,3);
\draw (3,3) ..controls (3,2.7) and (0,2.7).. (0,3);
\draw (3,0) ..controls (3,-0.3) and (0,-0.3).. (0,0);
\draw (0,0)--(0,3); \draw (3,0)--(3,3);
\draw (2,2.8) ..controls (2,2) and (0.3,1).. (0.3,-0.15);
\draw (2.7,2.85) ..controls (2.7,2.7) and (2.8,2.5)..(3,2.5);
\draw (0,2) .. controls (0.5,2) and (2,0.2).. (2,-0.2);
\draw (0.3,2.85) ..controls (0.3,2.7) and (0.1,2.3)..(0,2.3);
\draw (3,2) ..controls (2.7,1.9) and (2.7,2).. (2.7,-0.15);
\node at (3.8,1.5) {$s_2s_1s_2s_0$};
\fill (2.7,2.85) circle(0.06);
\fill (2,2.78) circle(0.06);
\fill (1,2.78) circle(0.06);
\fill (0.3,2.85) circle(0.06);
\fill (2.7,-0.15) circle(0.06);
\fill (2,-0.22) circle(0.06);
\fill (1,-0.22) circle(0.06);
\fill (0.3,-0.15) circle(0.06);
\end{tikzpicture}\qquad
\end{center}
\caption{Cylindrical diagrams for the elements of the coextended affine Weyl group.}
\label{fi:cylindrical}
\end{figure}

A wiring diagram is as a collection of $N$ paths $\gamma_k$ on the cylinder $\mathbb{R}/N\mathbb{Z}\times[0,1]$, connecting bijectively the integral points on one edge of the cylinder to the integral points on the other. The paths are considered up to homotopy and up to the diffeomorphisms of the torus, preserving the boundary of the cylinder point-wise -- i.e. the Dehn twists; this rule ensures the relation $\Lambda^N=1$. The group product corresponds to gluing the right side of one cylinder to the left of another. A generator $s_i$ corresponds to the diagram with single crossing of the path, connecting the $i$-th point to the $i+1$-st, with the path, connecting $i+1$-st to the $i$-th, and the rest of the paths remain horizontal. The diagram corresponding to the generator $\Lambda$ connects $i$-th point to the $i+1$-st for any $i$ (see fig.~\ref{fi:cylindrical}).

Conversely, given a wiring diagram with only pairwise crossings and with the paths, going monotonously from left to right, one can associate it with a word in $\widehat{W}^\sharp$ in the following way. Cut the cylinder into a rectangle by a line going from left to the right side of the cylinder starting on the left at some point between the points $(N-1 \mod{N})$ and $(0\mod N)$ and such that for every face of the diagram it passes through, it enters it through the leftmost point and leaves through the rightmost one. Once the cylinder is cut we can we can associate to every crossing point the generator $s_i$ just as we did for finite diagrams on a strip with additional generator $s_0$ corresponding to the intersection points occurring on the cut. The power of the coextension generator $\Lambda$ is determined by the number of the segment where the cut meets the right side of the cylinder.

The coextension homomorphism $\widehat{W}^\sharp\to \mathbb{Z}/N\mathbb{Z}$ can also be described as the intersection index of the diagram with a generator of the cylinder: orient all paths from left to right and count the number of intersection, taking into account the orientation, with the straight horizontal line from left to right, disjoint from the integral points\footnote{
Yet another way to compute this coextension is the following. To every path one can associate an integer measuring the difference of coordinates of both ends of the curve on the universal cover of the cylinder. The sum of such numbers divided by $N$ modulo $N$ is the desired homomorphism. One can easily check, that such sum is always divisible by $N$, and that it is remains unchanged modulo $N^2$ under the Dehn twist.}.

\begin{figure}[ht]
\begin{center}
\begin{tikzpicture}
\draw[thin] (4.2,0)--(4.2,3.2);
\draw[thin] (0,0)--(0,3.2);
\draw (0,0.4)--(2.4,0.4)
..controls (2.7,0.4)..(3,0)..controls (3.3,0.4)..(3.6,0.4)--(4.2,0.4);
\draw (0,1.2)--(0.6,1.2)
..controls(0.9,1.2)..(1.2,1.6)--(1.8,2.4)..controls(2.1,2.8)..(2.4,2.8)
..controls(2.7,2.8)..(3,3.2)..controls(3.3,2.8)..(3.6,2.8)--(4.2,2.8);
\draw (0,2)--(0.6,2)
..controls(0.9,2)..(1.2,1.6)..controls(1.5,1.2)..(1.8,1.2)--(3,1.2)
..controls(3.3,1.2)..(3.6,1.6)..controls(3.9,2)..(4.2,2);
\draw (0,2.8)--(1.2,2.8)
..controls(1.5,2.8)..(1.8,2.4)..controls(2.1,2)..(2.4,2)--(3,2)
..controls(3.3,2)..(3.6,1.6)..controls(3.9,1.2)..(4.2,1.2);
\node at (1.2,-0.3) {$s_{2}$};
\node at (1.8,-0.3) {$s_{1}$};
\node at (3,-0.3) {$s_{0}$};
\node at (3.6,-0.3) {$s_{2}$};
\end{tikzpicture}
\begin{tikzpicture}
\draw[thin] (4.2,0)--(4.2,3.2);
\draw[thin] (0,0)--(0,3.2);
\draw (0,0.8)--(1.8,0.8)
..controls(2.1,0.8)..(2.4,1.2)..controls(2.7,1.6)..(3,1.6)--(4.2,1.6);
\draw (0,1.6)--(1.8,1.6)
..controls(2.1,1.6)..(2.4,1.2)..controls(2.7,0.8)..(3,0.8)--(4.2,0.8);
\draw (0,2.4)
..controls(0.3,2.4)..(0.6,2.8)..controls(0.9,3.2)..(1.2,3.2)--(4.2,3.2);
\draw (0,3.2)
..controls(0.3,3.2)..(0.6,2.8)..controls(0.9,2.4)..(1.2,2.4)--(4.2,2.4);
\node at (0.6,-0.3) {$s_{\bar{0}}$};
\node at (2.4,-0.3) {$s_{\bar{2}}$};
\end{tikzpicture}
\begin{tikzpicture}
\fill[color=gray!20] (0,0)--(3,0)..controls (2.7,0.4)..(2.4,0.4)--(0,0.4)--cycle;
\draw (3,0)..controls (2.7,0.4)..(2.4,0.4)--(0,0.4);
\fill[color=gray!20] (4.2,0)--(3,0)..controls (3.3,0.4)..(3.6,0.4)--(4.2,0.4)--cycle;
\draw (3,0)..controls (3.3,0.4)..(3.6,0.4)--(4.2,0.4);
\draw[fill=gray!20] (0,0.8)--(1.8,0.8)
..controls(2.1,0.8)..(2.4,1.2)..controls(2.7,0.8)..(3,0.8)--(4.2,0.8)--(4.2,1.2)
..controls (3.9,1.2)..(3.6,1.6)..controls(3.3,1.2)..(3,1.2)--(1.8,1.2)
..controls (1.5,1.2)..(1.2,1.6)..controls (0.9,1.2)..(0.6,1.2)
--(0,1.2)--cycle;
\draw[fill=gray!20] (0,1.6)--(1.2,1.6)..controls(0.9,2)..(0.6,2)--(0,2)--cycle;
\draw[fill=gray!20] (1.2,1.6)--(1.8,1.6)
..controls(2.1,1.6)..(2.4,1.2)..controls(2.7,1.6)..(3,1.6)--(3.6,1.6)
..controls(3.3,2)..(3.0,2)--(2.4,2)..controls (2.1,2)..(1.8,2.4)--cycle;
\draw[fill=gray!20] (4.2,1.6)--(4.2,2)..controls(3.9,2)..(3.6,1.6)--cycle;
\draw[fill=gray!20] (0,2.4)..controls(0.3,2.4)..(0.6,2.8)--(0,2.8)--cycle;
\draw[fill=gray!20] (0.6,2.8)..controls(0.9,2.4)..(1.2,2.4)--(1.8,2.4)
..controls(1.5,2.8)..(1.2,2.8)--cycle;
\draw[fill=gray!20] (4.2,2.4)--(4.2,2.8)--(3.6,2.8)
..controls(3.3,2.8)..(3,3.2)..controls(2.7,2.8)..(2.4,2.8)
..controls(2.1,2.8)..(1.8,2.4)--cycle;
\fill[color=gray!20] (0,3.2)..controls(0.3,3.2)..(0.6,2.8)
..controls(0.9,3.2)..(1.2,3.2)--(4.2,3.2)--cycle;
\draw (0,3.2)..controls(0.3,3.2)..(0.6,2.8)..controls(0.9,3.2)..(1.2,3.2)--(4.2,3.2);
\draw[thin] (4.2,0)--(4.2,3.2);
\draw[thin] (0,0)--(0,3.2);
\node at (0.6,-0.3) {$s_{\bar{0}}$};
\node at (1.2,-0.3) {$s_{2}$};
\node at (1.8,-0.3) {$s_{1}$};
\node at (2.4,-0.3) {$s_{\bar{2}}$};
\node at (3,-0.3) {$s_{0}$};
\node at (3.6,-0.3) {$s_{2}$};
\end{tikzpicture}
\end{center}
\caption{A: wiring diagram; B: wiring diagram for the second copy of the group $W$; C:~Thurston diagram.}
\label{fi:wiring}
\end{figure}
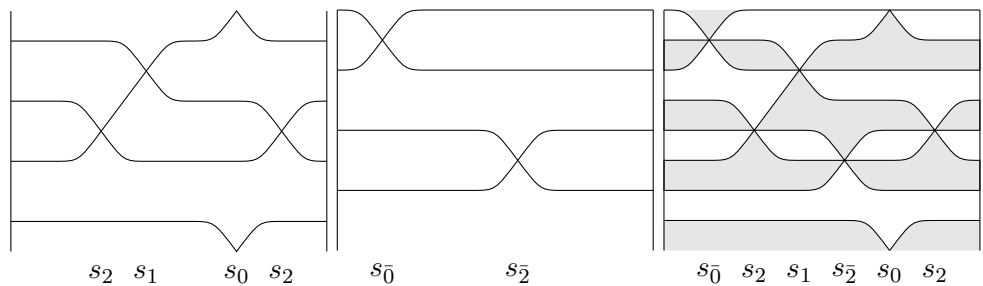

To present elements of the group $(\widehat{W}\times \widehat{W})^\sharp$  it is sufficient to give just two wiring diagrams, one for each copy of $\widehat{W}$, with the same number of lines and the same coextension. However, to present a particular decomposition of this group as a product of generators we will use Thurston diagrams (see appendix \ref{ap:thurston}), which are constructed from two wiring diagrams corresponding to two $W$-factors, drawn on the same cylinder, do that the second is slightly shifted down with respect to the first one, with two additional conditions:
(1) every vertical line intersects lines of both diagrams alternatively; (2) as a consequence all intersection points must be triple; (3) the lines corresponding to the first (resp. second diagrams) are oriented from left to right (resp. from right to left), see fig.~\ref{fi:wiring}(C). The triple intersection points are therefore of two types -- when two lines from the first diagram intersect a line from the second and visa versa; they correspond to the generators $s_i$ with $i>0$ or $i<0$, respectively.

Below we will draw cylindrical diagrams on rectangles assuming that the bottom side of the rectangle is identified with the top one.

\subsection{Realisations of the coextended loop group}

We shall use two different realisations of the group $\widehat{PGL}^\sharp(N)$: the infinite matrix realisation and the loop realisation.
In the infinite matrix realisation it can be identified with the group of infinite matrices $\{A_I^J|I,J\in \mathbb Z\}$, considered up to multiplication by a constant and satisfying the following two conditions:
\be
\label{quasi}
A_I^J = 0 \mbox{ for }|I-J|\gg 0
\\
A_{I+N}^{J+N} = xA_{I}^{J}\mbox{ for some }x\in \mathbb{C}^\times
\ee
In the Laurent realisation the group $\widehat{PGL}^\sharp(N)$ can be identified with the group of expressions $A(\lambda)T_x$, where $T_x$ is the operator of multiplicative shift by $x$:
\begin{equation}\label{shiftT}
T_x = \exp\left({\log x\lambda\frac{\partial}{\partial \lambda}}\right) = x^{\lambda \partial/\partial \lambda}
\end{equation}
and $A(\lambda)$ is a Laurent polynomial with values in $N\times N$ matrices, considered again up to a multiplicative constant. The multiplication rule of such expressions is therefore
\begin{equation}\label{AshifT}
A_1(\lambda)T_{x_1}\cdot A_2(\lambda)T_{x_2}=A_1(\lambda)A_2(x_1\lambda)T_{x_1x_2}
\end{equation}
The correspondence between the infinite matrix and the loop realizations is given by the isomorphism
\begin{equation}
A_I^J \mapsto \sum_{K\in \mathbb Z} A_{I}^{J+KN}\lambda^KT_x
\end{equation}
where in the r.h.s. $I,J\in 1,\ldots,N$ and $x$ is the quasi-periodicity factor from the condition (\ref{quasi}).


The loop group $\widehat{PGL}^\sharp(N)$ has one more triple of generators in addition to
those of $PGL(N)$, which we denote $E_0$, $E_{\bar 0}$ and $H_0(x)$.
In loop representation the matrices $E_i$ for $i\neq 0,\bar{0}$ coincide with the corresponding matrices for the corresponding finite-dimensional group $PGL(N)$, and the matrices $H_i(x)$ get multiplied by $T_x$, i.e.
\begin{equation}
\label{loop}
H_i(x) =
\begin{pmatrix}
x&0&&\cdots&&0\\[-5pt]
0&\ddots&&&&0\\[-5pt]
&&x&&\\[-5pt]
\vdots&&&1&\\[-5pt]
&&&&\ddots&0\\[-2pt]
0&&\cdots&&0&1
\end{pmatrix}T_x,\quad\quad
E_i = E_{\bar i}^{\rm tr} =
\begin{pmatrix}
1&0&&\cdots&&0\\[-5pt]
0&\ddots&&&&0\\[-5pt]
&&1&1&\\[-5pt]
\vdots&&&1&\\[-5pt]
&&&&\ddots&0\\[-2pt]
0&&\cdots&&0&1
\end{pmatrix}
\end{equation}
for $i>0$. For $i=0$ we have additionally $H_0(x)=T_x$, and
\begin{equation}
\label{loopextra}
E_0 =
\begin{pmatrix}
1&0&\cdots&0\\[-5pt]
\vdots&\ddots&&\vdots\\[-5pt]
0&&\ddots&0\\[-2pt]
\lambda&\cdots&0&1
\end{pmatrix},\qquad
E_{\bar{0}} =
\begin{pmatrix}
1&0&\cdots&\lambda^{-1}\\[-5pt]
\vdots&\ddots&&\vdots\\[-5pt]
0&&\ddots&0\\[-2pt]
0&\cdots&0&1
\end{pmatrix}
\end{equation}
It is also useful to introduce the element $\Lambda \in \widehat{PGL}(N)$, having in the infinite matrix presentation the form $\Lambda_I^J=\delta_I^{J+1}$, or in the loop representation
\begin{equation}
\label{lpshft}
\Lambda =
\begin{pmatrix}
0&1&\cdots&0\\
\vdots&\ddots&\ddots&\vdots\\
0&\cdots&0&1\\
\lambda&\cdots&0&0
\end{pmatrix}
\end{equation}
This matrix has the property
\begin{equation}
\label{ddshft}
\Lambda E_i \Lambda^{-1}= E_{i+1},\ \ \
\Lambda H_i(z) \Lambda^{-1}= H_{i+1}(z),\ \ \ i\in{\mathbb{Z}/N\mathbb{Z}}
\end{equation}
i.e. the operator $\Lambda$ acts as a unit shift along the Dynkin diagram of $\widehat{PGL(N)}$.

\subsection{Integrable systems on double Bruhat cells 
\label{ss:ouris}}

The Weyl group for affine Lie group is infinite (for example the product of all generators has  infinite order) and therefore, in contrast to the finite-dimensional case, one can consider arbitrarily long words, corresponding to the Poisson submanifolds of arbitrary large dimensions. The space of $\Ad$-invariant functions on $\widehat{PGL}^\sharp(N)$ is finitely generated, but it is infinitely generated on its subgroup $\widehat{PGL}(N)$.

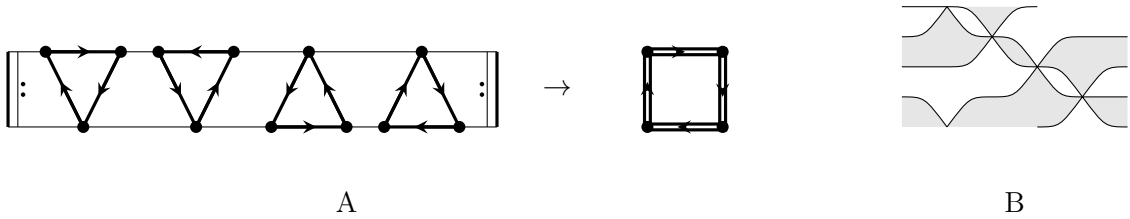
\begin{figure}[ht]
\begin{center}
\begin{tikzpicture}
 \draw[very thick,shorten >=0.4cm,-stealth] (0.5,1) -- (1.5,1);
\draw[very thick] (0.5,1) -- (1.5,1);
\draw[very thick,shorten >=0.45cm,-stealth] (1.5,1) -- (1,0);
\draw[very thick] (1.5,1) -- (1,0);
\draw[very thick,shorten >=0.45cm,-stealth] (1,0) -- (0.5,1);
\draw[very thick] (1,0) -- (0.5,1);
\fill (0.5,1) circle (0.08);
\fill (1.5,1) circle (0.08);
\fill (1,0) circle (0.08);
\draw[very thick,shorten <=0.4cm,stealth-] (2,1) -- (3,1);
\draw[very thick] (2,1) -- (3,1);
\draw[very thick,shorten <=0.45cm,stealth-] (3,1) -- (2.5,0);
\draw[very thick] (3,1) -- (2.5,0);
\draw[very thick,shorten <=0.45cm,stealth-] (2.5,0) -- (2,1);
\draw[very thick] (2.5,0) -- (2,1);
\fill (2,1) circle (0.08);
\fill (3,1) circle (0.08);
\fill (2.5,0) circle (0.08);
\draw[very thick,shorten >=0.4cm,-stealth] (3.5,0) -- (4.5,0);
\draw[very thick] (3.5,0) -- (4.5,0);
\draw[very thick,shorten >=0.45cm,-stealth] (4.5,0) -- (4,1);
\draw[very thick] (4.5,0) -- (4,1);
\draw[very thick,shorten >=0.45cm,-stealth] (4,1) -- (3.5,0);
\draw[very thick] (4,1) -- (3.5,0);
\fill (3.5,0) circle (0.08);
\fill (4.5,0) circle (0.08);
\fill (4,1) circle (0.08);
\draw[very thick,shorten <=0.4cm,stealth-] (5,0) -- (6,0);
\draw[very thick] (5,0) -- (6,0);
\draw[very thick,shorten <=0.45cm,stealth-] (6,0) -- (5.5,1);
\draw[very thick] (6,0) -- (5.5,1);
\draw[very thick,shorten <=0.45cm,stealth-] (5.5,1) -- (5,0);
\draw[very thick,] (5.5,1) -- (5,0);
\fill (5,0) circle (0.08);
\fill (6,0) circle (0.08);
\fill (5.5,1) circle (0.08);
\draw (0,0)--(6.5,0);
\draw (0,1)--(6.5,1);
\draw[very thick] (0,0)--(0,1);
\draw[very thick](6.5,0)--(6.5,1);
\draw[thin] (0.13,0)--(0.13,1);
\draw[thin](6.37,0)--(6.37,1);
\fill (0.2,0.57) circle (0.03);
\fill (0.2,0.43) circle (0.03);
\fill (6.3,0.57) circle (0.03);
\fill (6.3,0.43) circle (0.03);
\draw (7.3,0.5) node {$\rightarrow$};
\draw[very thick, double distance = 1pt] (8.5,0) -- (9.5,0)--(9.5,1)-- (8.5,1)--(8.5,0);
\draw[very thick,-stealth] (8.9,1) -- (9,1);
\draw[very thick,stealth-] (8.9,0) -- (9.0,0);
\draw[very thick,stealth-] (9.5,0.4) -- (9.5,0.5);
\draw[very thick,stealth-] (8.5,0.6) -- (8.5,0.5);
\fill (8.5,0) circle (0.08);
\fill (8.5,1) circle (0.08);
\fill (9.5,0) circle (0.08);
\fill (9.5,1) circle (0.08);
\node at (4.5,-1) {A};
\end{tikzpicture}\qquad\qquad\qquad
\begin{tikzpicture}
\draw[black,fill=gray!20](0,1.6)--(0.6,1.6)..controls(0.9,1.6)..(1.2,1.2)..controls(1.5,1.6)..(1.8,1.6);
\draw[black,fill=gray!20](0,0.8)--(0.6,0.8)..controls(0.9,0.8)..(1.2,1.2)..controls(0.9,1.2)..(0.6,1.6)..controls(0.3,1.2)..(0,1.2);
\fill[color=gray!20](0,0.4)..controls(0.3,0.4)..(0.6,0)--(0,0);
\draw (0,0.4)..controls(0.3,0.4)..(0.6,0);
\draw[black,fill=gray!20] (0.6,0)..controls(0.9,0.4)..(1.2,0.4)..controls(1.5,0.4)..(1.8,0.8)..controls(2.1,0.4)..(2.4,0.4)..controls(2.1,0)..(1.8,0);
\draw[black,fill=gray!20](1.2,1.2)..controls(1.5,1.2)..(1.8,0.8)..controls(1.5,0.8)..(1.2,1.2);
\draw[black,fill=gray!20](3,1.2)--(2.4,1.2)..controls(2.1,1.2)..(1.8,0.8)..controls(2.1,0.8)..(2.4,0.4)..controls(2.7,0.8)..(3,0.8);
\draw[black,fill=gray!20](3,0)..controls(2.7,0)..(2.4,0.4)--(3,0.4);
\node at (1.5,-1) {B};
\end{tikzpicture}
\end{center}
\caption{Exchange graph (A) and Thurston diagram (B) for the cell of $\widehat{PGL}^\sharp(2)/\Ad H$ corresponding to the word $0\bar{0}1\bar{1}$.}
\label{fi:todasl3h}
\end{figure}

The cluster parameterisation of double Bruhat cells for the loop group $\widehat{PGL}(N)^\sharp$ goes exactly along the lines we had for simple groups. The only difference on the level of exchange graphs (see appendix~\ref{ap:exchange}) is that one should add one extra line, an example of the exchange graph for $u=0\bar{0}1\bar{1}$ is shown on fig.~\ref{fi:todasl3h}(A). The Thurston diagrams for $\hat{G}^u$ are on a cylinder with left and right boundaries, which replaces the infinite strip for
the case of simple groups,  while the Thurston diagrams for $\hat{G}^u/\Ad H$ are on a torus, which
is obtained after gluing the left with the right boundary of a cylinder (i.e. on fig.~\ref{fi:todasl3h}(B) one should identify top with bottom and left with right).

Let $u=s_{i_1},\ldots,s_{i_n}\Lambda^k$ be an arbitrary reduced decomposition of $u\in (\widehat{W}\times \widehat{W})^\sharp$. The corresponding double Bruhat cell can be parametrised by $n$ variables $\boldsymbol{x}=(x_1,\ldots,x_n)$ by the formula
\be
\label{laxmap}
(x_1,\ldots,x_n) \mapsto \mathcal{A}(\lambda,\boldsymbol{x})=H_{i_1}(x_1)E_{i_1}\cdots H_{i_n}(x_n)E_{i_n}\Lambda^k
\ee
In order to get an element of $\widehat{PGL}(N)/\Ad H$ and not of $\widehat{PGL}^\sharp(N)/\Ad H$ we need to impose in \rf{laxmap} that
\begin{equation}
\label{trishift}
\prod_j x_j = 1
\end{equation}
If this condition is satisfied the matrix $\mathcal{A}(\lambda,\boldsymbol{x})$ is just a matrix with entries being Laurent polynomials, but defined up to multiplication by a constant $\alpha$ and up to a shift $\lambda \to \beta\lambda$ (which is just a result of its conjugation by $H_0(\beta)=T_\beta$).

Consider the function
\begin{equation}
\label{eq:partition-Lie}
S_{s_{i_1},\ldots,s_{i_n}\Lambda^k}(\lambda,\mu,\boldsymbol{x})=
\det(\mathcal{A}(\lambda,\boldsymbol{x})-\mu))
\end{equation}
which we shall call the generating function of the integrable model.
This function is ill-defined due to the ambiguity described above. However, let us consider it, as defined up to a transformation
$$
S_{s_{i_1},\ldots,s_{i_n}\Lambda^k}(\lambda,\mu,\boldsymbol{x})\simeq\alpha S_{s_{i_1},\ldots,s_{i_n}\Lambda^k}(\beta\lambda,\gamma\mu,\boldsymbol{x})
$$
where $\alpha$, $\beta$ and $\gamma$ are arbitrary functions of $x_1,\ldots,x_n$, then it becomes already well-defined.

Computing determinant in \rf{eq:partition-Lie}, this function is expressed as Laurent
polynomial of $\lambda$ and $\mu$:
\be
\label{S-curve}
S_{s_{i_1},\ldots,s_{i_n}\Lambda^k}(\lambda,\mu,\boldsymbol{x})=\sum_{i,j} \mathcal{H}_{ij}(\boldsymbol{x})\lambda^i\mu^j
\ee
Denote by $\Delta$ the Newton polygon of \rf{S-curve}, i.e. a convex hull on the plane of the points $(ij)$ for which the functions $\mathcal{H}_{ij}(\boldsymbol{x})\neq 0$ identically. Fix any three corners of $\Delta$, and adjust $\alpha, \beta$ and $\gamma$ to make the corresponding coefficients $\mathcal{H}_{ij}$ to be equal to unities, the resulting partition function will be called \textit{normalised}.

\begin{theorem}
The coefficients $\{\mathcal{H}_{ij}(\boldsymbol{x})\}$ of expansion \rf{S-curve} of the normalised partition functions are well defined Ad-invariant functions on the Poisson submanifold $\hat{G}^{s_{i_1},\ldots,s_{i_n}\Lambda^k}/\Ad H$, and therefore they Poisson-commute with each other.
\end{theorem}
The proof that they indeed form a set of integrals of motion of an integrable system, and the discussion of the properties of this system will be done below, after we establish the isomorphism of these systems with the integrable systems of Goncharov and Kenyon.

\section{Dimers
\label{s:dimers}}

In this section we define, following \cite{GK}, the dimer partition function for a bipartite graph with weights attached to its edges. We show that partition functions on certain graphs can be computed as minors of certain matrices of small size (related to the number of vertices). We define then the spectral variety, and describe the GK integrable systems, constructed out of certain bipartite graphs on a torus.

\subsection{Recollection about dimers}

Let $\Gamma$ be a graph, denote by $E_\Gamma$ and $V_\Gamma$ the set of its edges and vertices, respectively. A \textit{dimer configuration} on $\Gamma$ is a subset $D\subset E_\Gamma$, such that every vertex is contained in exactly one edge of $D$. Denote by $\mathcal{D}_\Gamma$ the set of all dimer configurations on $\Gamma$, which we assume to be nonempty. Fix a function ${\bf A}:e\mapsto A_e$  associating a complex number, called \emph{weight}, to any edge $e\in E_\Gamma$.

Consider the sum over all dimer configurations
\begin{equation}
\label{diparf}
S_\Gamma({\bf A}) = \sum_{D\in \mathcal{D}_\Gamma}\ \prod_{e\in D}A_e
\end{equation}
and call it the \textit{dimer partition function} of the graph $\Gamma$. This partition function is a polynomial of the weight variables $A_e$ with unit coefficients.

One can slightly generalise this construction for the graphs called \textit{open}, where one allows the \textit{terminal} edges, having one special univalent vertex, also called \textit{terminal}. Graphs without terminal edges are called \textit{closed}. One can identify two terminal vertices, and then erase resulting two-valent vertex: this procedure can be used to glue together two graphs. Conversely, one can cut an edge and declare emerging vertices to be terminal.

A dimer configuration on an open graph is a collection of edges containing every internal vertex exactly once without any condition for the terminal edges. Denote by $\mathcal{T}_{\Gamma}$ the set of terminal edges of the graph. For every subset $T\subset \mathcal{T}_{\Gamma}$ denote by $\mathcal{D}_\Gamma(T)$ the set of dimer configurations with terminal edges from $T$ occupied and all other terminal edges free. One can define the partition function with boundary
\begin{equation}\label{defpfcond}
S_\Gamma({\bf A},T) = \sum_{D\in \mathcal{D}_\Gamma(T)}\ \prod_{e\in D}A_e
\end{equation}
All such partition functions can be put together into a single generating partition function. Fix an ordering of the terminal edges and introduce odd variables $\bxi=\{\xi_t|t\in \mathcal{T}_{\Gamma}\}$. Define then
$$
S_\Gamma({\bf A},\bxi)=\sum_{T\subset \mathcal{T}_{\Gamma}} S_{\Gamma}({\bf A}, T) \prod_{t\in D\cap T} \xi_t,
$$
where the product is taken respecting the fixed order of the terminal edges.

This generating function is especially convenient, since it behaves in a nice way under gluing the terminal edges together. Namely, let the graph $\Gamma'$ be obtained from the graph $\Gamma$ by gluing a terminal edge $i$ with weight $A_i$ to the neighbor edge $j$, which follows $i$ with respect to the chosen order and carries the weight $A_j$. The resulting edge $e$ carries the weight $A_e=A_iA_j$.

\begin{lemma}\label{le:product}
$S_{\Gamma'}({\bf A}',\bxi')=\int S_{\Gamma}({\bf A},\bxi)e^{\xi_i\xi_j}d\xi_jd\xi_i$
\end{lemma}
Here $\bxi'$ denotes the collections $\bxi$ with $\xi_i$ and $\xi_j$ removed, while ${\bf a}'$ is the collection ${\bf a}$ with the weight $A_e=A_iA_j$ replacing $A_i$ and $A_j$.
The proof of this Lemma is obvious, since
$$
\int S_{\Gamma}({\bf A},\bxi)e^{\xi_i\xi_j}d\xi_jd\xi_i = \int S_{\Gamma}({\bf A},\bxi)d\xi_jd\xi_i + \int S_{\Gamma}({\bf A},\bxi)\xi_i\xi_jd\xi_jd\xi_i
$$
The first term in the r.h.s. gives the sum over the dimer configurations with both terminal edges $i$ and $j$ occupied, which are in bijection with dimer configurations for the graph $\Gamma'$ with the edge $e$ occupied. Similarly the second term gives the sum over configurations, where the edge $e$ is empty.

\subsection{Matrices and dimers on a disc}

Now we are going to show, that for any matrix one can associate a bipartite graph with terminal edges and with the weights such that its partition function is a generating functions for the minors of the matrix. Multiplication of matrices corresponds to gluing graphs.
Let $M$ be an $n\times n$ matrix, $\bxi=(\xi_1,\ldots,\xi_n)$ and $\Beta=(\eta^1\ldots\eta^n)$ be two sets of the Grassmann variables. Define the minor generating function by
\be
\label{mgf}
S(M,\bxi,\Beta)=\exp(\sum_{i,j} M^i_j\xi_i\eta^j)
\ee
This function satisfies the following relations:

\begin{lemma}\label{le:matrices}
\begin{enumerate}
 \item\label{product} $S(M_1M_2,\bxi',\Beta')=\int S(M_1,\bxi',\Beta)S(M_2,\bxi,\Beta')\prod_ie^{ -\xi_i\eta^i}{d\xi^id\eta_i}$
 \item\label{generating} $S(M,\bxi,\Beta)= 1 + \sum\limits_{k=1}^n\mathop{\sum\limits_{i_i<\cdots<i_k}}\limits_{j_1<\cdots<j_k}M^{i_1,\ldots,i_k}_{j_1,\ldots,j_k}\xi_{i_1}\cdots \xi_{i_k}\eta^{j_k}\cdots\eta^{j_1}$
 \item\label{determinant} $\int S(M,\bxi,\Beta)e^{\mu\sum \eta^i\xi_i}\prod_i{d\xi^id\eta_i}=\det(\mu-M)$
\item\label{pluecker} $\displaystyle\sum_i\left(\xi_i\frac{\partial}{\partial \tilde{\xi}_i}-\tilde{\eta}^i\frac{\partial}{\partial \eta^i}\right)S(M,\bxi,\Beta)S(M,\tilde{\bxi},\tilde{\Beta})=0$
\end{enumerate}
\end{lemma}
Here $M^{i_1,\ldots,i_k}_{j_1,\ldots,j_k}$ denotes the minor of the matrix $M$, corresponding to the rows $i_1<\cdots<i_k$ and columns $j_1<\cdots<j_k$. The property \ref{generating} shows that the function $S(M,\bxi,\Beta)$ can be considered as a generating function of all minors of the matrix $M$. The property \ref{pluecker} is equivalent to the Pl\"ucker relations for the minors of a matrix. A proof of this Lemma is given in appendix~\ref{ap:super}.

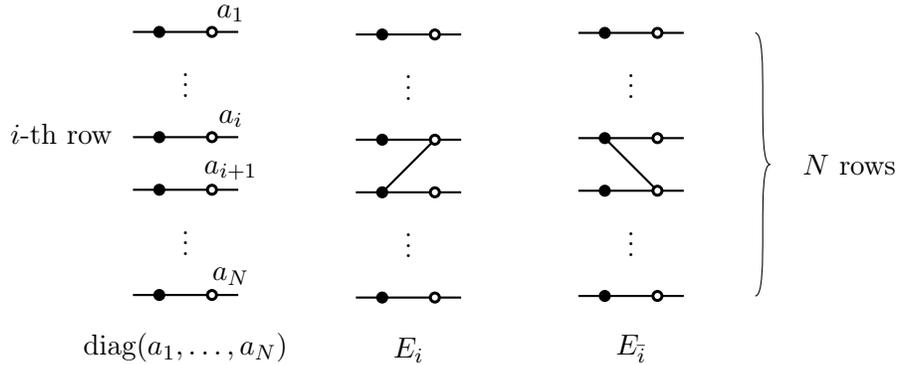
\begin{figure}[H]
\begin{center}
\begin{tikzpicture}
\foreach \y in {0,1.4,2.1,3.5}
{\draw[thick] (-0.35,\y)-- (1.05,\y);
 \fill (0,\y) circle(0.08);
 \fill (0.7,\y) circle(0.08);\fill[color=white] (0.7,\y) circle(0.04);
}
\node at (-1.3,2.15) {$i$-th row};
\node at (0.35,-0.7) {$\mbox{diag}(a_1,\ldots,a_N)$};
\node[above] at (0.95,0) {$a_N$};
\node[above] at (0.95,1.4) {$a_{i+1}$};
\node[above] at (0.95,2.1) {$a_i$};
\node[above] at (0.95,3.5) {$a_1$};
\node at (0.35,2.9) {$\vdots$};
\node at (0.35,0.8) {$\vdots$};
\end{tikzpicture}\qquad
\begin{tikzpicture}
\draw[thick] (0,1.4)--(0.7,2.1);
\foreach \y in {0,1.4,2.1,3.5}
{\draw[thick] (-0.35,\y)-- (1.05,\y);
 \fill (0,\y) circle(0.08);
 \fill (0.7,\y) circle(0.08);\fill[color=white] (0.7,\y) circle(0.04);
}
\node at (0.35,-0.7) {$E_i$};
\node at (0.35,2.9) {$\vdots$};
\node at (0.35,0.8) {$\vdots$};
\end{tikzpicture}\qquad\qquad
\begin{tikzpicture}
\draw[thick] (0,2.1)--(0.7,1.4);
\foreach \y in {0,1.4,2.1,3.5}
{\draw[thick] (-0.35,\y)-- (1.05,\y);
 \fill (0,\y) circle(0.08);
 \fill (0.7,\y) circle(0.08);\fill[color=white] (0.7,\y) circle(0.04);
}
\node at (0.35,2.9) {$\vdots$};
\node at (0.35,0.8) {$\vdots$};
\node at (0.35,-0.7) {$E_{\bar{i}}$};
\node[right] at (2.5,1.75) {$N$ rows};
\draw (2,0)..controls (2.2,0) and (2.0,1.75)..(2.2,1.75);
\draw (2,3.5)..controls (2.2,3.5) and (2.0,1.75)..(2.2,1.75);
\end{tikzpicture}
\end{center}
\caption{Graphs corresponding to the generators of $GL(N)$. The weight variables $\mathbf{a}$ are all assigned to the terminal edges (to the right ones, associated with the $\Beta$-variables in our choice). All other weights are equal to unities.}
\label{fi:graphsmatrices}
\end{figure}

Now it is easy to demonstrate, that
\begin{lemma}
For every matrix $M$ one can construct a bipartite planar graph with ordered edges $\Gamma_M$ with the weights $\boldsymbol{A}$ on edges, such that its minor generating function $S(M,\boldsymbol{\xi},\boldsymbol{\eta})$ coincides with generating partition function of the graph $S_{\Gamma_M}(\boldsymbol{A},\boldsymbol{\xi},\boldsymbol{\eta})$.
\end{lemma}
Indeed, comparing the property \ref{product} and Lemma~\ref{le:product} we see, that gluing the graphs corresponds to the matrix products. Therefore, one needs just to verify, that weighted graphs exist for a system of generators of $GL(N)$, i.e. to diagonal matrices and exponentiated simple roots \rf{genslN}. By listing all dimer configuration one can check, that graphs shown on  fig.~\ref{fi:graphsmatrices} correspond to the matrices $H_{\bf a}=\diag(a_1,\ldots,a_N)$, $E_i$ and $E_{\bar{i}}$, with the minor generating functions $S(H_{\bf a},\bxi,\Beta)=\exp(\sum_ja_j\xi_j\eta^j)$, $S(E_i,\bxi,\Beta)=\exp(\xi_i\eta^{i+1}+\sum_j\xi_j\eta^j)$ and $S(E_{\bar{i}},\bxi,\Beta)=\exp(\xi_{i+1}\eta^{i}+\sum_j\xi_j\eta^j)$, respectively. We assume, that terminal edges of the graphs are ordered counter-clockwise: first on the left side from top to the bottom ($\bxi$-variables) and then on the right from bottom to the top ($\Beta$-variables).

\subsection{Dimer partition functions with signs}

This simple construction unfortunately fails to construct in a universal way the graphs, corresponding to the generators \rf{loopextra}, \rf{lpshft} of the affine groups $\widehat{PGL}(N)$. One can do this in much simpler way by assigning the signs to dimer configurations, in particular to present an algorithm, valid
for the groups $\widehat{PGL}(N)$ both in the cases of integer and odd $N$. Moreover, this is the most
direct way to explain the appearance of terms of different signs in (\ref{eq:partition-Lie}), to be  compared below with certain dimer partition functions.
Assume from now that the graph $\Gamma$ is embedded into a surface $\Sigma$ with terminal vertices (if any) on the boundary of $\Sigma$.

Any dimer configuration can be viewed as an element of the space of chains $C_1(\Gamma,\mathbb{Z}/2\mathbb{Z})$. Fix an element $Q\in Z^1(\Sigma,\mathbb{Z}/2\mathbb{Z})\otimes Z^1(\Sigma,\partial \Sigma, \mathbb{Z}/2\mathbb{Z})$, such that its projection to homology gives the canonical element in $H^1(\Sigma,\mathbb{Z}/2\mathbb{Z})\otimes H^1(\Sigma,\partial \Sigma, \mathbb{Z}/2\mathbb{Z})$. Such $Q$ can be interpreted as a quadratic form on  $C_1(\Gamma,\mathbb{Z}/2\mathbb{Z})$ and in particular on its subspace generated by dimer configurations. The difference of any two such forms is a linear form. To construct such form one can choose a collection of cycles $\{\alpha_i\}$ representing a basis of $H^1(\Sigma,\mathbb{Z}/2\mathbb{Z})$ and a collection of cycles $\{\beta_i\}$ representing the dual basis in $H^1(\Sigma,\partial \Sigma, \mathbb{Z}/2\mathbb{Z})$. Then the form is given by $Q=\sum_i\alpha_i\otimes \beta_i$ and its value on a dimer configuration $D$ is equal to $Q(D)=\sum_i\langle\alpha_i,D\rangle\langle\beta_i,D\rangle$.

Using this quadratic form one can modify \rf{diparf} to
\be
\label{diparfs}
S_{\Gamma,Q}({\bf A}) = \sum_{D\in \mathcal{D}_\Gamma}(-1)^{Q(D)}\prod_{e\in D}A_e
\ee
and (\ref{defpfcond}) to
\be
S_{\Gamma,Q}({\bf A},T) = \sum_{D\in \mathcal{D}_\Gamma(T)}(-1)^{Q(D)}\ \prod_{e\in D}A_e
\ee
Cutting a surface $\Sigma$ provided with a form $Q$ by a closed curve induces a form $Q'$ on the resulting surface $\Sigma'$. The form $Q'$ is defined by the condition that for any dimer configuration $D$ on $\Gamma\subset\Sigma$ we have $Q(D)=Q(D')$, where $D'$ is the dimer configuration obtained by cutting the configuration $D$. In this case we say that the forms $Q$ and $Q'$ are compatible.

On the pictures we show the cocycles $\alpha_i\in Z^1(\Sigma,\mathbb{Z}/2\mathbb{Z})$ as well as from $\beta_i\in Z^1(\Sigma,\partial \Sigma, \mathbb{Z}/2\mathbb{Z})$ as corresponding cycles on the surface transversal to the graph $\Gamma$. The former ones are represented by cycles possibly going from boundary to boundary and the latter by closed cycles.

\subsection{$\widehat{PGL}(N)$ and dimers on a cylinder}

In this section we present generating functions of elements of the affine group $\widehat{PGL}(N)\subset \widehat{PGL}^\sharp(N)$ as dimer partition functions with signs. A product of generators \rf{loop}, \rf{loopextra} of $\widehat{PGL}^\sharp(N)$ belongs to $\widehat{PGL}(N)$ if the condition \rf{trishift} is satisfied.

\begin{figure}[H]
\begin{center}
\begin{tikzpicture}
\draw[color=gray!40,dashed] (-0.15,1.85)--(0.95,1.85);
\draw[dashed] (0.95,1.85)--(1.25,1.85);
\draw[dashed] (-0.1,-0.2)..controls(-0.3,-0.2)and (-0.3,3.7)..(-0.1,3.7);
\draw[color=gray!50,dashed] (-0.1,-0.2)..controls(0.1,-0.2)and (0.1,3.7)..(-0.1,3.7);
\draw (-0.3,-0.2)--(1.1,-0.2);
\draw (-0.3,3.7)--(1.1,3.7);
\draw (-0.3,-0.2)..controls(-0.5,-0.2)and (-0.5,3.7)..(-0.3,3.7);
\draw[color=gray!50] (-0.3,-0.2)..controls(-0.1,-0.2)and (-0.1,3.7)..(-0.3,3.7);
\draw (1.1,-0.2)..controls(0.9,-0.2)and (0.9,3.7)..(1.1,3.7);
\draw (1.1,-0.2)..controls(1.3,-0.2)and (1.3,3.7)..(1.1,3.7);
\foreach \y in {0,3.5}
{\draw[thick] (-0.35,\y)-- (1.05,\y);
 \fill (0,\y) circle(0.08);
 \fill (0.7,\y) circle(0.08);\fill[color=white] (0.7,\y) circle(0.04);
}
\foreach \y in {1.4,2.1}
{\draw[ultra thick, color=white] (-0.45,\y)-- (0.95,\y);
 \draw[thick] (-0.45,\y)-- (0.95,\y);
 \fill (-0.1,\y) circle(0.08);
 \fill (0.6,\y) circle(0.08);\fill[color=white] (0.6,\y) circle(0.04);
}
\node at (-1.3,2.15) {$i$-th row};
\node at (0.35,-0.7) {$\mbox{diag}(a_1,\ldots,a_N)$};
\node[above] at (0.95,0) {$a_N$};
\node[above] at (0.95,1.4) {$a_{i+1}$};
\node[above] at (0.95,2.1) {$a_i$};
\node[above] at (0.95,3.5) {$a_1$};
\node at (0.35,2.9) {$\vdots$};
\node at (0.35,0.8) {$\vdots$};
\end{tikzpicture}
\begin{tikzpicture}
\draw[color=gray!40,dashed] (-0.15,1.85)--(0.95,1.85);
\draw[dashed] (0.95,1.85)--(1.25,1.85);
\draw[dashed] (-0.1,-0.2)..controls(-0.3,-0.2)and (-0.3,3.7)..(-0.1,3.7);
\draw[color=gray!50,dashed] (-0.1,-0.2)..controls(0.1,-0.2)and (0.1,3.7)..(-0.1,3.7);
\draw (-0.3,-0.2)--(1.1,-0.2);
\draw (-0.3,3.7)--(1.1,3.7);
\draw (-0.3,-0.2)..controls(-0.5,-0.2)and (-0.5,3.7)..(-0.3,3.7);
\draw[color=gray!50] (-0.3,-0.2)..controls(-0.1,-0.2)and (-0.1,3.7)..(-0.3,3.7);
\draw (1.1,-0.2)..controls(0.9,-0.2)and (0.9,3.7)..(1.1,3.7);
\draw (1.1,-0.2)..controls(1.3,-0.2)and (1.3,3.7)..(1.1,3.7);
\draw[thick] (-0.1,1.4)--(0.6,2.1);
\foreach \y in {0,3.5}
{\draw[ultra thick, color=white] (-0.35,\y)-- (1.05,\y);
\draw[thick] (-0.35,\y)-- (1.05,\y);
 \fill (0,\y) circle(0.08);
 \fill (0.7,\y) circle(0.08);\fill[color=white] (0.7,\y) circle(0.04);
}
\foreach \y in {1.4,2.1}
{\draw[thick] (-0.45,\y)-- (0.95,\y);
 \fill (-0.1,\y) circle(0.08);
 \fill (0.6,\y) circle(0.08);\fill[color=white] (0.6,\y) circle(0.04);
}
\node at (0.35,-0.7) {$E_i$};
\node at (0.35,2.9) {$\vdots$};
\node at (0.35,0.8) {$\vdots$};
\end{tikzpicture}\qquad
\begin{tikzpicture}
\draw[color=gray!40,dashed] (-0.15,1.85)--(0.95,1.85);
\draw[dashed] (0.95,1.85)--(1.25,1.85);
\draw[dashed] (-0.1,-0.2)..controls(-0.3,-0.2)and (-0.3,3.7)..(-0.1,3.7);
\draw[color=gray!50,dashed] (-0.1,-0.2)..controls(0.1,-0.2)and (0.1,3.7)..(-0.1,3.7);
\draw (-0.3,-0.2)--(1.1,-0.2);
\draw (-0.3,3.7)--(1.1,3.7);
\draw (-0.3,-0.2)..controls(-0.5,-0.2)and (-0.5,3.7)..(-0.3,3.7);
\draw[color=gray!50] (-0.3,-0.2)..controls(-0.1,-0.2)and (-0.1,3.7)..(-0.3,3.7);
\draw (1.1,-0.2)..controls(0.9,-0.2)and (0.9,3.7)..(1.1,3.7);
\draw (1.1,-0.2)..controls(1.3,-0.2)and (1.3,3.7)..(1.1,3.7);
\draw[thick] (-0.1,2.1)--(0.6,1.4);
\foreach \y in {0,3.5}
{\draw[ultra thick, color=white] (-0.35,\y)-- (1.05,\y);
\draw[thick] (-0.35,\y)-- (1.05,\y);
 \fill (0,\y) circle(0.08);
 \fill (0.7,\y) circle(0.08);\fill[color=white] (0.7,\y) circle(0.04);
}
\foreach \y in {1.4,2.1}
{\draw[thick] (-0.45,\y)-- (0.95,\y);
 \fill (-0.1,\y) circle(0.08);
 \fill (0.6,\y) circle(0.08);\fill[color=white] (0.6,\y) circle(0.04);
}
\node at (0.35,-0.7) {$E_{\bar{i}}$};
\node at (0.35,2.9) {$\vdots$};
\node at (0.35,0.8) {$\vdots$};
\node[right] at (2.5,1.75) {$N$ rows};
\draw (2,0)..controls (2.2,0) and (2.0,1.75)..(2.2,1.75);
\draw (2,3.5)..controls (2.2,3.5) and (2.0,1.75)..(2.2,1.75);
\end{tikzpicture}
\end{center}
\caption{Graphs on a cylinder, corresponding to the common generators of $PGL(N)$ and $\widehat{PGL}(N)$. The cycle $\alpha\in H^1(\Sigma,\mathbb{Z}/2\mathbb{Z})$ is depicted by vertical dashed line around a cylinder, while the cycle $\beta\in H^1(\Sigma,\partial \Sigma, \mathbb{Z}/2\mathbb{Z})$ corresponds to the horisontal dashed line, connecting two boundaries of a cylinder on its back side.}\label{fi:pglcyl}
\end{figure}
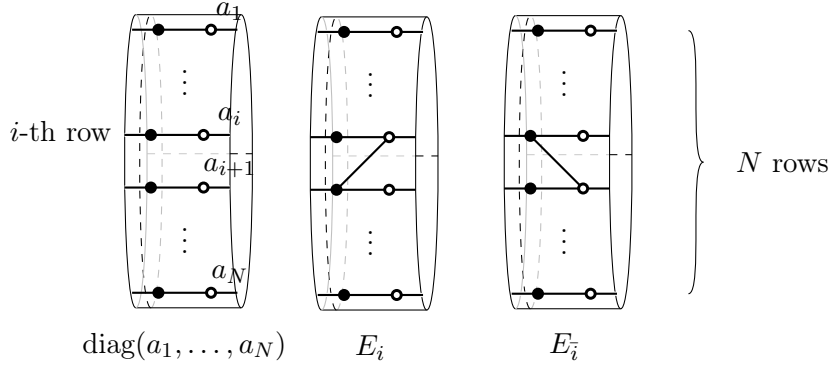

The generators \rf{loop} correspond just to the same graphs as on fig.~\ref{fi:graphsmatrices}, but which are embedded now into short horizontal cylinders, as shown on fig.~\ref{fi:pglcyl}. For these matrices the formulae (\ref{defpfcond}) and \rf{diparfs} give the same answer, if we choose the cycles $\alpha$ and $\beta$ as shown on the figure: indeed, since the cycle $\beta$ goes from boundary to boundary not  intersecting the graph $Q(D)=\langle\alpha,D\rangle\langle\beta,D\rangle=0$, and the signs of all dimer configurations are positive.

\begin{figure}[H]
\begin{center}
\begin{tikzpicture}
\draw[color=gray!40,dashed] (-0.15,1.85)--(0.95,1.85);
\draw[dashed] (0.95,1.85)--(1.25,1.85);
\draw[dashed] (-0.1,-0.2)..controls(-0.3,-0.2)and (-0.3,3.7)..(-0.1,3.7);
\draw[color=gray!50,dashed] (-0.1,-0.2)..controls(0.1,-0.2)and (0.1,3.7)..(-0.1,3.7);
\draw (-0.3,-0.2)--(1.1,-0.2);
\draw (-0.3,3.7)--(1.1,3.7);
\draw (-0.3,-0.2)..controls(-0.5,-0.2)and (-0.5,3.7)..(-0.3,3.7);
\draw[color=gray!50] (-0.3,-0.2)..controls(-0.1,-0.2)and (-0.1,3.7)..(-0.3,3.7);
\draw (1.1,-0.2)..controls(0.9,-0.2)and (0.9,3.7)..(1.1,3.7);
\draw (1.1,-0.2)..controls(1.3,-0.2)and (1.3,3.7)..(1.1,3.7);
\draw[thick] (0,3.5)..controls (0,3.5) and (0,3.7).. (0.1,3.7);
\draw[thick,color=gray!20] (0.1,3.7)..controls(0.2,3.7) and (0.5,-0.2)..(0.6,-0.2);
\draw[thick] (0.6,-0.2)..controls(0.7,-0.2) and (0.7,0)..(0.7,0);
\node[color=gray!70] at (0.1,3.0) {$-\lambda$};
\foreach \y in {0,3.5}
{\draw[ultra thick, color=white] (-0.35,\y)-- (1.05,\y);
\draw[thick] (-0.35,\y)-- (1.05,\y);
 \fill (0,\y) circle(0.08);
 \fill (0.7,\y) circle(0.08);\fill[color=white] (0.7,\y) circle(0.04);
}
\foreach \y in {1.4,2.1}
{\draw[ultra thick, color=white] (-0.45,\y)-- (0.95,\y);
 \draw[thick] (-0.45,\y)-- (0.95,\y);
 \fill (-0.1,\y) circle(0.08);
 \fill (0.6,\y) circle(0.08);\fill[color=white] (0.6,\y) circle(0.04);
}
\node at (0.35,-0.7) {$E_0$};
\node at (0.35,2.9) {$\vdots$};
\node at (0.35,0.8) {$\vdots$};
\end{tikzpicture}\qquad
\begin{tikzpicture}
\draw[color=gray!40,dotted] (-0.15,1.85)--(0.95,1.85);
\draw[dashed] (0.95,1.85)--(1.25,1.85);
\draw[dashed] (-0.1,-0.2)..controls(-0.3,-0.2)and (-0.3,3.7)..(-0.1,3.7);
\draw[color=gray!50,dotted] (-0.1,-0.2)..controls(0.1,-0.2)and (0.1,3.7)..(-0.1,3.7);
\draw (-0.3,-0.2)--(1.1,-0.2);
\draw (-0.3,3.7)--(1.1,3.7);
\draw (-0.3,-0.2)..controls(-0.5,-0.2)and (-0.5,3.7)..(-0.3,3.7);
\draw[color=gray!50] (-0.3,-0.2)..controls(-0.1,-0.2)and (-0.1,3.7)..(-0.3,3.7);
\draw (1.1,-0.2)..controls(0.9,-0.2)and (0.9,3.7)..(1.1,3.7);
\draw (1.1,-0.2)..controls(1.3,-0.2)and (1.3,3.7)..(1.1,3.7);
\draw[thick] (0.7,3.5)..controls (0.7,3.5) and (0.7,3.7).. (0.6,3.7);
\draw[thick,color=gray!20] (0.6,3.7)..controls(0.5,3.7) and (0.2,-0.2)..(0.1,-0.2);
\draw[thick] (0.1,-0.2)..controls(0,-0.2) and (0,0)..(0,0);
\node[color=gray!70] at (0.6,3.0) {$-\lambda$};
\foreach \y in {0,3.5}
{\draw[ultra thick, color=white] (-0.35,\y)-- (1.05,\y);
 \draw[thick] (-0.35,\y)-- (1.05,\y);
 \fill (0,\y) circle(0.08);
 \fill (0.7,\y) circle(0.08);\fill[color=white] (0.7,\y) circle(0.04);
}
\foreach \y in {1.4,2.1}
{\draw[ultra thick, color=white] (-0.45,\y)-- (0.95,\y);
 \draw[thick] (-0.45,\y)-- (0.95,\y);
 \fill (-0.1,\y) circle(0.08);
 \fill (0.6,\y) circle(0.08);\fill[color=white] (0.6,\y) circle(0.04);
}
\node at (0.35,-0.7) {$E_{\bar{0}}$};
\node at (0.35,2.9) {$\vdots$};
\node at (0.35,0.8) {$\vdots$};
\end{tikzpicture}\qquad
\begin{tikzpicture}
\draw[color=gray!40,dashed] (-0.15,1.85)--(0.95,1.85);
\draw[dashed] (0.95,1.85)--(1.25,1.85);
\draw[dashed] (-0.1,-0.2)..controls(-0.3,-0.2)and (-0.3,3.7)..(-0.1,3.7);
\draw[color=gray!40,dashed] (-0.1,-0.2)..controls(0.1,-0.2)and (0.1,3.7)..(-0.1,3.7);
\draw (-0.3,-0.2)--(1.1,-0.2);
\draw (-0.3,3.7)--(1.1,3.7);
\draw (-0.3,-0.2)..controls(-0.5,-0.2)and (-0.5,3.7)..(-0.3,3.7);
\draw[color=gray!50] (-0.3,-0.2)..controls(-0.1,-0.2)and (-0.1,3.7)..(-0.3,3.7);
\draw (1.1,-0.2)..controls(0.9,-0.2)and (0.9,3.7)..(1.1,3.7);
\draw (1.1,-0.2)..controls(1.3,-0.2)and (1.3,3.7)..(1.1,3.7);
\draw[thick] (0,3.5)..controls (0,3.5) and (0,3.7).. (0.1,3.7);
\draw[thick,color=gray!20] (0.1,3.7)..controls(0.2,3.7) and (0.5,-0.2)..(0.6,-0.2);
\draw[thick] (0.6,-0.2)..controls(0.7,-0.2) and (0.7,0)..(0.7,0);
%
\draw[thick] (-0.1,2.1)--(0.2,2.4);
\draw[thick] (-0.1,1.4)--(0.6,2.1);
\draw[thick] (0,0)--(0.3,0.3);
\draw[ultra thick, color=white] (0.7,3.5)--(0.4,3.2);\draw[thick] (0.7,3.5)--(0.4,3.2);\
\draw[thick] (0.6,1.4)--(0.3,1.1);
\foreach \y in {0,3.5}
{\draw[ultra thick, color=white] (0.7,\y)-- (1.05,\y);
\draw[thick] (-0.35,\y)-- (0,\y);
\draw[thick] (0.7,\y)--(1.05,\y);
 \fill (0,\y) circle(0.08);
 \fill (0.7,\y) circle(0.08);\fill[color=white] (0.7,\y) circle(0.04);
}
\foreach \y in {1.4,2.1}
{\draw[ultra thick, color=white] (0.6,\y)-- (0.95,\y);
\draw[thick] (-0.45,\y)-- (-0.1,\y);
\draw[thick] (0.6,\y)--(0.95,\y);
 \fill (-0.1,\y) circle(0.08);
 \fill (0.6,\y) circle(0.08);\fill[color=white] (0.6,\y) circle(0.04);
}

\node at (0.35,2.9) {$\vdots$};
\node at (0.35,0.8) {$\vdots$};
\node at (0.35,-0.7) {$\Lambda$};
\node[right] at (2.5,1.75) {$N$ rows};
\draw (2,0)..controls (2.2,0) and (2.0,1.75)..(2.2,1.75);
\draw (2,3.5)..controls (2.2,3.5) and (2.0,1.75)..(2.2,1.75);
\end{tikzpicture}\qquad
\end{center}
\caption{Graphs corresponding to the spectral parameter dependent generators of $\widehat{PGL}(N)$. Dashed lines again denote the cycles $\alpha$ (vertical) and $\beta$ (horisontal). Now for the $\lambda$-dependent minors $\langle\beta,D\rangle=1$ and corresponding $Q(D)=\langle\alpha,D\rangle\langle\beta,D\rangle=\langle\alpha,D\rangle$ counts the number of
the occupied terminal edges.}
\label{fi:graphaffine}
\end{figure}
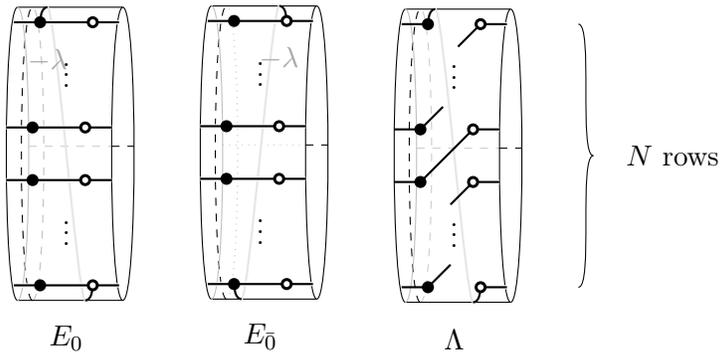

The generators \rf{loopextra}, \rf{lpshft} specific to $\widehat{PGL}(N)$ can be represented by the graphs, embedded into a cylinder, shown on fig.~\ref{fi:graphaffine}. All generators \rf{loopextra}, \rf{lpshft}
have just a single nonvanishing matrix element, equal to $\lambda^{\pm 1}$ beyond the main, next upper and next lower diagonals. Therefore, all $\lambda$-dependent minors, dependently on their size $k=1,\ldots,N$, will have positive (negative) sign for odd (even) $k$ respectively in the expansions of minor generating
functions \rf{mgf} for the matrices $E_0(\lambda)$, $E_{\bar 0}(\lambda)$ and $\Lambda(\lambda)$. The most
simple way to take into account these signs is:
\begin{itemize}
  \item to assign weight $-\lambda$ (instead of naive $+\lambda$) to every link on the graph, corresponding to a $\lambda$-dependent matrix element;
  \item to assign the signs to all dimer configurations, according to the rule \rf{diparfs}.
\end{itemize}

It is easy to see that gluing cylinders along boundary circles is compatible with the quadratic form $Q=\alpha\otimes\beta$, where $\beta$ is the horizontal generator of the cylinder and $\alpha$ is a vertical section of it in a small vicinity of one of the boundaries. Observe that such form $Q$ does not depend on two possible choices of the cycle $\alpha$ (since the number of vertices between them is even), and that this form is compatible with gluing two cylinders together.

\subsection{Dimers and the discrete Dirac operator
\label{ss:dirac}}

Let us show now, that for a closed surface the partition function with signs \rf{diparfs} coincides with determinant of a matrix, which can be called the \textit{discrete Dirac operator}. To define this operator one needs first to fix a spin structure on the surface (even with this fixing the determinant is defined up to a sign). However the dependence on the spin structure is not essential, and the equivalence class of the partition function does not depends on this choice. The construction belongs mainly to P.~Kasteleyn \cite{Kasteleyn} and it was generalised to higher genus by D.~Cimasoni and N.~Reshetikhin \cite{Cimasoni-Reshetikhin}.

Recall that a spin structure on a closed surface can be identified \cite{Johnson} with a quadratic form $Q$ on $H_1(\Sigma,\mathbb{Z}/2\mathbb{Z})$ satisfying the compatibility condition with the intersection index $\langle,\rangle$:
\begin{equation}\label{quadraticcompat}
Q(x+y)+Q(x)+Q(y)=\langle x,y\rangle\mod{2}
\end{equation}
The spin structure on the surface can be encoded by a \textit{Kasteleyn marking} of the edges of a bipartite graph. A Kasteleyn marking of a graph $\Gamma$ embedded into a surface $\Sigma$ is a marking of edges of a graph by numbers $\boldsymbol{K}=K_e\in\{\pm1\}$ such that for any face $f$ with $l(f)$ sides the product of $K_e$ over the sides be equal to $(-1)^{l(f)/2+1}$ (i.e. with monodromy $-1$ around faces with number of sides divisible by 4 and $+1$ otherwise). Such a marking can be considered as a cochain $\boldsymbol{K} \in C^1(\Gamma,\{\pm1\})$ with the condition $\partial \boldsymbol{K}=\boldsymbol{R}$ where $\boldsymbol{R}$ is the cocycle taking value $(-1)^{l/2+1}$ on a face $f$.

Given two dimer configurations $D_1$ and $D_2$ one can define their \textit{mutual parity} $\sigma(D_1,D_2)$ taking value in $\{\pm 1\}$ as follows. Every dimer configuration $D$ defines a map, denoted by the same letter $D:B\to W$ from black to white vertices. The parity $\sigma(D_1,D_2)$ is just
a parity of the map $(D_1-D_2): W\to W$, it obviously satisfies the property $\sigma(D_1,D_2)\sigma(D_2,D_3)\sigma(D_3,D_1)=1$ for any three dimer configurations $D_1,D_2$ and $D_3$.
The proof of the following lemma is essentially contained in \cite{Cimasoni-Reshetikhin}.

\begin{lemma}\label{le:Kastel}
\begin{enumerate}
\item A Kasteleyn marking exists and is unique up to a cocycle from $Z^1(\Gamma,\{\pm1\})$.
\item For any Kasteleyn marking $\boldsymbol{K}$ and for any two dimer configurations $D_0$ and $D_1$, such that $D_1-D_0 \in B_1(\Sigma, \mathbb{Z}/2\mathbb{Z})$, we have $$\sigma(D_1,D_0)\langle(D_1-D_0),\boldsymbol{K}\rangle=1$$
\item For any Kasteleyn marking and for any four dimer configurations $D_0,D_1,D_2$ and  $D_3$, such that $D_0+D_1+D_2+D_3 \in B_1(\Sigma, \mathbb{Z}/2\mathbb{Z})$, one has:
$$\sigma(D_1,D_2)\sigma(D_0,D_3)\langle (D_1-D_0),(D_2-D_0)\rangle\langle(D_0+D_1+D_2+D_3),\boldsymbol{K}\rangle=1$$
\end{enumerate}
\end{lemma}

\noindent\textbf{Corollary:} Fix a dimer configuration $D_0$ and consider a function on the set of dimer configurations $q:\mathcal{D}_\Gamma\to \{\pm1\}$ defined by $q(D)= \langle (D-D_0),\boldsymbol{K}\rangle\sigma(D,D_0)$. The statement 2 of Lemma~\ref{le:Kastel} implies that $q(D_0)=1$ and its value depends only on the homology class of $D-D_0$, and the statement 3 is equivalent to representing it in the form $q(D)=(-1)^{Q(D-D_0)}$, where $Q$ is quadratic form, representing a spin structure.

Define a linear operator $\mathfrak{D}(\boldsymbol{A}): \mathbb{C}^{W}\to \mathbb{C}^B$, by
$$\mathfrak{D}(\boldsymbol{A})=\sum_e \langle e,\boldsymbol{KA}\rangle b(e)w(e),$$
where $b(e)$ is the embedding of $\mathbb{C}$ into the factor corresponding to the black end of the edge $e$, $w(e)$ is the projection onto the factor corresponding to the white end of the edge. In other terms the matrix of the operator $\mathfrak{D}(\boldsymbol{A})$ has rows numerated by black vertices, columns numerated by white vertices and a matrix element between a white vertex $w$ and a black vertex $b$ is equal to the sum $\mathfrak{D}(\boldsymbol{A})_{bw}=\sum_{b\,\stackrel{e}{\rightarrow}\,w} K_eA_e$ over all edges $e$ connecting $b$ and $w$.

Changing $\boldsymbol{A}$ and $\boldsymbol{K}$ by a coboundary (making a discrete gauge transformation) amounts to the multiplication of $\mathfrak{D}(\boldsymbol{A})$ by a nondegenerate diagonal matrices from both sides. Thus the property of $\mathfrak{D}(\boldsymbol{A})$ being degenerate does not depends only on the class of $\boldsymbol{A}$ in $H^1(\Gamma)$.

\begin{lemma}
$\det \mathfrak{D}(\boldsymbol{A})=\pm S_\Gamma(\boldsymbol{A}).$
\end{lemma}

\noindent\textit{Remark:} Though the Dirac operator $\mathfrak{D}$ depends on the spin structure on the surface via the Kasteleyn orientation $\boldsymbol{K}$, the dependence is given by a simple rule

$$\mathfrak{D}_{\boldsymbol{K'}}(\boldsymbol{A})=\mathfrak{D}_{\boldsymbol{K}}(\boldsymbol{AK'/K}).$$
The formula makes sense since the  ratio $\boldsymbol{K'/K}$ can be considered as belonging to $H^1(\Sigma)$. For this reason we omit the subscript $\boldsymbol{K}$.

\noindent\textit{Proof.} The determinant of the operator $\mathfrak{D}(\boldsymbol{A})$ is not defined before we fix an isomorphism between the volume forms on $\mathbb{C}^{W}$ and $\mathbb{C}^B$. A choice of a bijection between $W$ and $B$ establishes this isomorphism, but different choices give different signs of the determinant and therefore it is defined up to a sign. One can fix a bijection by choosing a dimer configuration $D_0$. Once the bijection is fixed one can speak about the parity $\sign(s)$ of a map $s\colon B\to W$ and the formula for the determinant reads as $\det \mathfrak{D}(\boldsymbol{A})=\sum_{s\colon B\to W}\prod_{b\in B}(-1)^{\sign(s)}\mathfrak{D}_{bs(b)}$, where the sum is taken over bijections $s\colon B\to W$. Every dimer configuration $D$ defines a bijection which we will denote also by $D$ and the bijections which do not correspond to dimer configurations give no contribution to the determinant. Therefore the formula for the determinant can be rewritten as
\be
\det \mathfrak{D}(\boldsymbol{A})=\sum_{D\in\mathcal{D}_\Gamma}(-1)^{\sigma(D,D_0)}\langle D,\boldsymbol{KA}\rangle=
\\
= \langle D_0,\boldsymbol{K}\rangle\sum_{D\in\mathcal{D}_\Gamma}(-1)^{\sigma(D,D_0)}\langle (D-D_0),\boldsymbol{K}\rangle\langle D,\boldsymbol{A}\rangle =\pm\sum_{D\in\mathcal{D}_\Gamma}(-1)^{Q(D-D_0)}\langle D,\boldsymbol{A}\rangle
\ee
where the last equality holds due the corollary of the Lemma \ref{le:Kastel}.

\subsection{Spectral submanifold and face partition function}

Consider now a closed surface $\Sigma$, a bipartite graph $\Gamma$ on it, and the corresponding dimer partition function with signs. For some set of values of the edge parameters $\boldsymbol{A}$ the partition function vanishes and, since it is determinant of the operator $\mathfrak{D}(\boldsymbol{A})$, the kernel of $\mathfrak{D}(\boldsymbol{A})$ forms a line bundle over this set. However, there is a large group of gauge symmetry acting on the space of parameters $\boldsymbol{A}$ leaving the vanishing set invariant. The spectral submanifold is just the quotient of the vanishing set by this group.

We proceed, first, to the description of the quotient of all edge parameters by the gauge group in terms of elementary algebraic topology.

Since the edges of a bipartite graph have canonical orientation (from the white vertex to the black one) any dimer configuration can be interpreted as an element of the space $C_1(\Gamma,\mathbb{Z})$ of 1-chains  on $\Gamma$ with integer coefficients. The boundary of any dimer configuration is just a 0-chain, which is a sum of black vertices minus the sum of the white vertices, which we shall denote by $\boldsymbol{v}=\boldsymbol{b}-\boldsymbol{w}$. Therefore the difference of any two dimer configurations is a cycle.

The weights $\boldsymbol{A}=\{A_e\}$ on the edges can be seen as elements of the space $C^1(\Gamma,\mathbb{C}^\times)$ of 1-cochains of $\Gamma$ with coefficients in the multiplicative group $\mathbb{C}^\times$, and the product $\prod_{e\in D}A_e=\langle \boldsymbol{A}, D\rangle$ is just the canonical pairing between 1-chains and 1-cochains. The expression for the partition function \rf{diparfs} can be  written down as
$$
S_\Gamma(\boldsymbol{A})=\sum_{D\in \mathcal{D}_\Gamma}(-1)^{Q(D)}\langle  \boldsymbol{A},D\rangle
$$
Below, unless stated otherwise, we assume all groups of chains with integer coefficients written additively, since they are just free finitely generated Abelian groups. On the other hand, all groups of cochains will be considered with coefficients in $\mathbb{C^\times}$ written multiplicatively, since they are algebraic tori.

The function $S_\Gamma(\boldsymbol{A})$ is homogeneous with respect to the action of the gauge group $C^0(\Gamma)$. Namely
$S_\Gamma(\boldsymbol{A}\partial\boldsymbol{g})=S_\Gamma(\boldsymbol{A})\langle \boldsymbol{g},\boldsymbol{w}-\boldsymbol{b}\rangle,$
where $\boldsymbol{b},\boldsymbol{w}\in C_0(\Gamma)$ are the sums of the black and the white vertices, respectively. Indeed
\be
S_\Gamma(\boldsymbol{A}\partial\boldsymbol{g})=\sum_{D\in \mathcal{D}_\Gamma}(-1)^{Q(D)}\langle  \boldsymbol{A}\partial\boldsymbol{g},D\rangle=\sum_{D\in \mathcal{D}_\Gamma}(-1)^{Q(D)}\langle  \boldsymbol{A},D\rangle\langle\partial\boldsymbol{g},D\rangle=
\\
=\sum_{D\in \mathcal{D}_\Gamma}(-1)^{Q(D)}\langle
\boldsymbol{A},D\rangle\langle\boldsymbol{g},\partial D\rangle=S_\Gamma(\boldsymbol{A})\langle \boldsymbol{g},\boldsymbol{w}-\boldsymbol{b}\rangle
\ee
It implies that the partition function $S_\Gamma(\cdot)$ as well as the pairings $\langle \cdot,D_0\rangle$ for any dimer configuration $D_0$ can be interpreted as sections of a line bundle over the quotient $H^1(\Gamma)=C^1(\Gamma)/\partial C^0(\Gamma)$. The ratio $S_\Gamma(\cdot)/\langle \cdot,D_0\rangle$ is just a function on $H^1(\Gamma)$. The vanishing locus of this section is our spectral variety $\mathcal{L}$.

Observe that since our graph is embedded into a closed surface $\Sigma$ we have the exact sequence
\begin{equation}\label{H1Gamma}
 1\to H^1(\Sigma)\to H^1(\Gamma)\stackrel{\partial}{\to} B^2(\Sigma)\to 1,
\end{equation}
where $B^2(\Sigma)$ is the group of 2-coboundaries, and in our case it is just the group of assignments of nonzero numbers with product one to the faces of the surface. This exact sequence defines on $H^1(\Gamma)$ the structure of a principal $H^1(\Sigma)$-bundle over $B^2(\Sigma)$.

Assume now, that $\Sigma$ is a torus, and let $\mathcal{M}$ be the moduli space of the plane curves in $H^1(\Sigma)$, which is just a two-dimensional algebraic torus, considered up to a (multiplicative) shift by elements of $H^1(\Sigma)$. Let \textit{universal Jacobian of planar curves} $\mathcal{J}$ be the space of pairs (curve, line bundle on it) also up to the action of $H^1(\Sigma)$.

The points of $\mathcal{M}$ correspond to the polynomials $P(\boldsymbol{\lambda})=\sum_{\boldsymbol{i}\in H_1(\Sigma)}\mathcal{H}_{\boldsymbol{i}}\langle\boldsymbol{i},\boldsymbol \lambda\rangle$, considered up to a transformation $P(\boldsymbol{\lambda})\mapsto \alpha P(\boldsymbol{\beta}\boldsymbol{\lambda})$, where $\alpha\in \mathbb{C}^\times$ and $\boldsymbol{\beta}\in H^1(\Sigma)$.
The spaces $\mathcal{M}$ and $\mathcal{J}$ are stratified by the finite dimensional subspaces $\mathcal{M}_\Delta$ and $\mathcal{J}_\Delta$ respectively, with curves having a given Newton polygon $\Delta \subset H_1(\Sigma)$.  The polygon $\Delta$ is just the convex hull of the points $\boldsymbol{i} \in H_1(\Sigma)$ for which $\mathcal{H}_{\boldsymbol{i}}\neq 0$.

The coefficients $\mathcal{H}_{\boldsymbol{i}}$ with $\boldsymbol{i}$ in the corners of the Newton polygon $\Delta$ must be nonvanishing. Chosen any three corners, one can reduce the polynomial $F(\boldsymbol{\lambda})=\sum_{\boldsymbol{i}\in \Delta}\mathcal{H}_{\boldsymbol{i}}(\boldsymbol{x})\langle\boldsymbol{i},\boldsymbol \lambda\rangle$ to a normal form (we call it \textit{normalised} equation of a curve) - with unit coefficients in these corners. The remaining coefficients are coordinates on the space $\mathcal{M}_\Delta$.

For any point $\boldsymbol{x}\in B^2(\Sigma)$ the intersection of the spectral variety to the fiber over $\boldsymbol{x}$ is a curve, and the kernel of the Dirac operator defines a line bundle on it provided the curve is nonsingular. The fiber is isomorphic to $H^1(\Sigma)$, but this isomorphism is defined up to a shift. Therefore the map $I:B^2(\Sigma)\to\mathcal{M}$ to the moduli of curves called the \textit{action map} and a rational map $J:B^2(\Sigma)\to\mathcal{J}$ to the universal Jacobian called the \textit{action-angle map} are well defined.

The exact sequence implies that noncanonically $H^1(\Gamma)$ is isomorphic to the product $H^1(\Sigma)\times B^2(\Sigma)$. Fixing such an isomorphism and choosing a representative of $H^1(\Gamma)$ in discrete connections we get a map $\boldsymbol{A}: B^2(\Sigma)\times H^1(\Gamma)\to C^1(\Gamma)$.
The \textit{face partition function} useful to compute the action map in coordinates is given by
$$
S^f_\Gamma(\boldsymbol{\lambda},\boldsymbol{x})=
S_\Gamma(\boldsymbol{A}(\boldsymbol{\lambda},\boldsymbol{x}))
$$
Choosing another isomorphism, the face partition function changes
\begin{equation}\label{normalisation}S^f_\Gamma(\boldsymbol{\lambda},\boldsymbol{x})\mapsto \alpha(\boldsymbol{x})S^f_\Gamma(\boldsymbol{\beta}(\boldsymbol{x})\boldsymbol{\lambda},\boldsymbol{x}),
\end{equation}
where $\alpha:B^2(\Sigma)\to \mathbb{C}^\times$ and $\boldsymbol{\beta}:B^2(\Sigma)\to H^1(\Sigma)$ are arbitrary homomorphisms. The face partition function
\be
\label{fpfc}
S^f_\Gamma(\boldsymbol{\lambda},\boldsymbol{x})=\sum_{\boldsymbol{i}\in \Delta}\mathcal{H}_{\boldsymbol{i}}(\boldsymbol{x})\langle\boldsymbol{\lambda},\boldsymbol{i}\rangle
\ee
thus defines a point of the moduli space $\mathcal{M}$ of plane curves for any value of $\boldsymbol{x}$. We can always choose the splitting making three corner coefficients to be unit, so that the face partition function becomes normalised.

Observe that the normalised face partition function does not depend on the spin structure. Indeed, the difference between two quadratic forms on $H_1(\Sigma,\mathbb{Z}/2\mathbb{Z})$ satisfying (\ref{quadraticcompat}) is a linear form $L\in H^1(\Sigma,\mathbb{Z}/2\mathbb{Z})$.  The class $L$ can be also viewed as an element of $H^1(\Sigma)$ if we identify $\mathbb{Z}/2\mathbb{Z}$ with the group $\{\pm 1\}\subset \mathbb{C}^\times$. Therefore the face partition function changes as $S^f_\Gamma(\boldsymbol{\lambda},\boldsymbol{x})\mapsto S^f_\Gamma(L\boldsymbol{\lambda},\boldsymbol{x})$ which is a particular case of the transformation (\ref{normalization}).
The line bundle on the spectral curve \rf{fpfc} will be given by the kernel of the discrete Dirac operator $\mathfrak{D}(\boldsymbol{A}(\boldsymbol{\lambda},\boldsymbol{x}))$.

\subsection{Dual surface, double partition function and the Poisson bracket 
}

Recall that embedding of a graph into an orientable surface induces a fat graph structure, namely a cyclic order of ends of edges, incident to every vertex. Conversely, for a given fat graph structure there exists the only surface and the only embedding, up to a diffeomorphism provided the embedding induces surjections of the fundamental groups.
Following \cite{GK}, we introduce the notion of \emph{dual embeddings} of a bipartite graph into two surfaces. Two embeddings $\Gamma\hookrightarrow\Sigma$ and $\Gamma\hookrightarrow\tilde{\Sigma}$ are called \emph{dual} if both induce surjections of the fundamental groups, and the fat graph structures induced by both embeddings coincide in the white vertices and are opposite in the black ones.

Dual surface $\tilde{\Sigma}$ can be also imagined as follows. A \textit{zig-zag path} on $\Gamma$ is a path turning left in every white vertex and right in a white one. Glue a disk to the graph $\Gamma$ along every zig-zag path. Since every edge of $\Gamma$ belongs to exactly two zig-zag paths, we obtain a closed smooth surface $\tilde{\Sigma}$. A zig-zag path is a boundary of the face of the dual surface $\tilde{\Sigma}$, but it may be nontrivial in homology of $\Sigma$. Conversely, a boundary of a face of $\Sigma$ maybe noncontractible in $\tilde{\Sigma}$. The sum in homology $H^1(\Sigma)$ of the zig-zag paths obviously vanishes and therefore they are sides of a convex polygon $\Delta\in H^1(\Sigma,\mathbb{R})$ with vertices in integer points. The polygon $\Delta$ is uniquely defined up to a shift.

Observe, that then we have obvious maps on the level of homology:
\be
i:H_1(\Gamma)\to H_1(\Sigma),\quad\tilde \imath: H_1(\Gamma)\to H_1(\tilde{\Sigma})
\\
i^*:H^1(\Sigma)\to H^1(\Gamma),\quad\tilde \imath^*: H^1(\Gamma)\to H^1(\tilde{\Sigma})
\ee
As it is shown in \cite{GK}, if the graph $\Gamma$ is minimal, the cokernel $Z$ of the map $i+\tilde\imath$ is finite and therefore the kernel $Z^*$ of the map $i^*\times \tilde\imath^*$ is finite.

Denote by $\mathcal{C}_1\subset H_1(\Gamma)$ the kernel of $i+\tilde\imath$ and by $\mathcal{C}^1$ the cokernel of $i^*\times \tilde\imath^*.$
Fixing a splitting of the exact sequence
\begin{equation}\label{H1Dual}
1\to Z^*\to H^1(\Sigma)\times H^1(\tilde\Sigma)\to H^1(\Gamma)\to \mathcal{C}^1\to 1
\end{equation}
one can construct the map
$\bar{\boldsymbol{A}}:H^1(\Sigma)\times H^1(\tilde\Sigma)\times \mathcal{C}^1\to C^1(\Gamma)$ and thus yet another version of the partition function, which we will call the \textit{double partition function}
$$
S^d_\Gamma(\boldsymbol{\lambda},\tilde{\boldsymbol{\lambda}},\boldsymbol{c})=
S_\Gamma(\bar{\boldsymbol{A}}(\boldsymbol{\lambda},\tilde{\boldsymbol{\lambda}},\boldsymbol{c}))
$$
The double partition function depends on the splitting and therefore is defined up to a transformation
$S^d_\Gamma(\boldsymbol{\lambda},\tilde{\boldsymbol{\lambda}},\boldsymbol{c})\mapsto \alpha(\boldsymbol{c})S^f_\Gamma(\boldsymbol{\beta}(\boldsymbol{c})\boldsymbol{\lambda},\tilde{\boldsymbol{\beta}}(\boldsymbol{c})\tilde{\boldsymbol{\lambda}},\boldsymbol{c}),$ where $\alpha:\mathcal{C}^1\to \mathbb{C}^\times$, $\boldsymbol{\beta}:\mathcal{C}^1\to H^1(\Sigma)$ and $\tilde{\boldsymbol{\beta}}:\mathcal{C}^1\to H^1(\tilde{\Sigma})$ are arbitrary homomorphims.




Observe that the sequence (\ref{H1Dual}) implies that the space $B^2(\Sigma)$ is isomorphic to $(H^1(\tilde{\Sigma})\times \mathcal{C}^1)/Z^*$. Indeed, this isomorphism is given by the differential $\partial:H^1(\tilde{\Sigma})\times \mathcal{C}^1\to B^2(\Sigma).$

Recall that the space $H^1(\tilde{\Sigma})$ has a canonical Poisson bracket, originated form the intersection index $\langle,\rangle_{\tilde{\Sigma}}$ on $H_1(\tilde{\Sigma})$. Indeed, for two functions $\langle\cdot,x\rangle$ and $\langle\cdot,y\rangle$ on $H^1(\tilde{\Sigma})$ defined by integral homology classes $x, y\in H_1(\tilde{\Sigma})$ the bracket is given by
\be
\label{GKbra}
\{\langle\cdot,x\rangle,\langle\cdot,y\rangle\}=\langle x,y\rangle \langle\cdot,x+y\rangle
\ee
where the pairing between $x$ and $y$ is the intersection index~\footnote{The functions on $H^1(\tilde{\Sigma})$ dual to the canonical $A$- and $B$- cycles in $H_1(\tilde{\Sigma})$, for some
 choice $\langle A_i,B_j\rangle = \delta_{ij}$ of the latter, can play the role of the (exponentiated) co-ordinates and momenta. It is easy to see, that this is literally true for the relativistic Toda systems, see sect.~\ref{ss:examples} and \cite{AMJGP}.}. This bracket is obviously invariant under the multiplicative shifts and in particular under the action of the finite group $Z^*$.
The product $H^1(\tilde{\Sigma})\times \mathcal{C}^1$ has thus a Poisson structure trivial on the second factor and the map $\partial$ induces one on $B^2(\Sigma)$.


Recall now, that the space $B^2(\Sigma)$ has natural coordinates, given by associating numbers to the faces of $\Gamma$, embedded into $\Sigma$. The bracket \rf{GKbra} in these coordinates is given by the combinatorial formula (cf. with the bracket \rf{clubra} in Appendix~\ref{ap:cluster})
\be
\label{GKbraG}
\{x_i,x_j\}=\varepsilon^{ij}x_ix_j
\ee
Here $i,j\in B^2(\Sigma)$ are two faces, and $\varepsilon^{ij}$ is the number of common edges counted with sign plus, if going from $i$ to $j$ we leave the white vertex on the right, and with the sign minus otherwise. In other terms the exchange graph of the Poisson bracket \rf{GKbra}, \rf{GKbraG} is just dual to the graph $\Gamma$.

\section{Dimers and integrable systems for the loop groups
\label{ss:dimerstoloop}}

The aim of this section is to formulate our main result - to show, that the class of the GK integrable systems coincides with the class of naturally defined integrable systems (sect.~\ref{ss:ouris}) on the Poisson submanifold of loop groups $\widehat{PGL}(N)$.

It was proven in \cite{GK}, that the action map $I_\Gamma:B^2(\Sigma)\to \mathcal{M}$ is integrable if the bipartite graph $\Gamma$ is minimal. In other words, integrability means that the set of coefficients $\mathcal{H}_{\boldsymbol{i}}$ of the normalised face partition function $S^f_\Gamma(\boldsymbol{\lambda},\boldsymbol{x}) = \sum_{\boldsymbol{i}}\mathcal{H}_{\boldsymbol{i}}(\boldsymbol{x})\langle \boldsymbol{i},\boldsymbol{\lambda}\rangle$ Poisson-commute with each other, and - when we forget the corner coefficients of the normalised partition function - they form a maximal commuting set of functionally independent integrals of motion.

One can easily convince himself that the number of integrals of motion, given by this construction, is maximal possible. Assume for simplicity, that all curves of the Thurston diagram have no self-intersections. In this case the number of triple crossings and thus the number of white faces is equal to $2\mathsf{S}$ where $\mathsf{S}$ is the area\footnote{This can be proven using convenient triangulation of the Newton polygon and the fact that the vector product of two oriented sides coincides with the intersection index of the corresponding curves on the Thurston diagram.} of the Newton polygon $\Delta$. Therefore the dimension of the phase space is equal to $2\mathsf{S}-1$.  The number of nontrivial Hamiltonians is equal to the number $\mathsf{N}$ of integral points belonging to $\Delta$. The number of Casimirs is equal to $\mathsf{B}-3$ where $\mathsf{B}$ is the number of the boundary points of $\Delta$. The maximal number of integrals is thus $(2\mathsf{S}-1-(\mathsf{B}-3))/2= \mathsf{S}-\mathsf{B}/2+1=\mathsf{N}$, according to the Pick's theorem.

Moreover in \cite{GK} it is shown that
\begin{enumerate}
\item If $\boldsymbol{i}$ is a corner of the boundary of the Newton polygon $\Delta$, then $\mathcal{H}_{\boldsymbol{i}}=c_i(\boldsymbol{x})$ is a Casimir function (monomial in $\{x_i\}$), i.e. $\{\mathcal{H}_{\boldsymbol{i}},f\}=0$ for any function $f=f(\boldsymbol{x})$ on $B^2(\Sigma)$.
\item Two minimal graphs $\Gamma_1$ and $\Gamma_2$ can be obtained from one another by spider moves (see sect.~\ref{ss:equgraphs}) if and only if the corresponding polygons $\Delta_1$ and $\Delta_2$ can be mapped one onto another by shift or $SL(2,\mathbb{Z})$-action. Such graphs are called equivalent.
\item Two integrable systems corresponding to the same polygon are isomorphic by a cluster transformation. 
\end{enumerate}

Hence, the GK integrable systems are enumerated by the Newton polygons $\Delta$, and the coordinate systems on their phase spaces - by the bipartite graphs $\Gamma$ on torus $\Sigma$. The integrable systems on the Poisson submanifolds of $\widehat{PGL}(N)$, constructed in sect.~\ref{ss:ouris}, are enumerated by cyclically irreducible elements $u\in (\widehat{W}\times \widehat{W})^\sharp$ in the coextended double Weyl groups, and the coordinates on their phase spaces - by decompositions of these elements into a reduced word of generators $u=s_{i_1}\cdots s_{i_l}\Lambda^k$. We need to establish the correspondence between these two sets of data.

Encode, first, both sets of combinatorial data with the Thurston diagrams. Indeed, given a reduced decomposition $u=s_{i_1}\cdots s_{i_l}\Lambda^k$ we can draw a corresponding Thurston diagram on a cylinder, and then gluing together the right and the left sides of this cylinder we get a diagram on the torus. The corresponding Newton polygon is (unique up to a shift) convex polygon in $H_1(\Sigma, \mathbb{R})$ with the sides, being classes of curves on the Thurston diagram. On the other hand, Thurston diagrams are in one-to one correspondence with the bipartite graphs $\Gamma$, see Appendix~\ref{ap:thurston}.

Now, for any Newton polygon $\Delta$ there exists a reduced word $s_{i_1}\cdots s_{i_l}\Lambda^k$ in the generators of $(\widehat{W}\times \widehat{W})^\sharp$ such that the curves of its glued Thurston diagram correspond to the sides of $\Delta$. Choose a coordinate system on the plane $H_1(\Sigma,\mathbb{R})$ such that none of the sides of the polygon is vertical. Then every side of the polygon gives a cycle on the torus $H_1(\Sigma,\mathbb{R})/H_1(\Sigma,\mathbb{Z})$. We can represent all these cycles by curves with nowhere vertical tangent and having minimal number of transversal intersections. Cutting the torus by a vertical line we get a collection of curves on a cylinder, the number of these curves $N$ will be twice the width of the polygon (half of the curves go to the left, and half to the right). Forgetting the left-going curves, we get a wiring diagram and the corresponding to it word in the generators $\{s_i|i\geq 0\}$ multiplied by $\Lambda^k$. The power $k$  is given by the difference in height of the leftmost and the rightmost vertex of the polygon. Forgetting curves going to the right we get in the same way a word in the generators $\{s_{\bar{i}}|\bar{i}\leq 0\}$ multiplied by $\Lambda^k$ with the same $k$. Shuffling these two words in any way and multiplying the result by $\Lambda^k$ we get a word in the generators of $(\widehat{W}\times \widehat{W})^\sharp$ such that the corresponding Thurston diagram has the desired Newton polygon $\Delta$.

For given reduced word $u\in (\widehat{W}\times \widehat{W})^\sharp$ there is a natural bijection between the set of variables parametrising the cell $G^u/\Ad H$ and the set of white faces of the Thurston diagram (and thus - of the faces of the corresponding bipartite graph $\Gamma$ on torus $\Sigma$). Indeed, given a reduced decomposition of $u$, the cell parametrising variables are in bijection with the generators $s_i$, which correspond to the triple intersections of the Thurston diagram. On the other hand every white face of the Thurston diagram on a torus with a chosen vertical direction and nowhere vertical tangents to the curves has its leftmost intersection point, and thus the white faces are also enumerated by the triple intersection points.

\begin{theorem}
Given a word $u=s_{i_1}\cdots s_{i_l}\Lambda^k$ in the coextended double Weyl group of the group $\widehat{PGL}(N)$, let $S_u(\lambda,\mu;\boldsymbol{x})=\det(\mathcal{A}_u(\lambda,\boldsymbol{x})-\mu)$ be the characteristic polynomial of the matrix, parametrising a double Bruhat cell $G^u/\Ad H$ in cluster coordinates $\boldsymbol{x}$, corresponding to given decomposition. Then this function coincides up to normalisation to the dimer face partition function of the bipartite graph on a torus, constructed out of the Thurston diagram, corresponding to the word $s_{i_1}\cdots s_{i_l}\Lambda^k$.
\end{theorem}

Most of the proof of this theorem is already done. Given a Thurston diagram corresponding to a word $u=s_{i_1}\cdots s_{i_l}\Lambda^k$, and at the same time to the graph $\Gamma$, one can cut it by vertical circles into $l$ elementary pieces having exactly one triple points inside. Observe that every slice corresponds to one of the graphs shown on fig.~\ref{fi:pglcyl} and fig.~\ref{fi:graphaffine}. One should just compare the monodromies around faces and zig-zags of the bipartite graph, used to compute the partition function $S_u(\lambda,\mu,\boldsymbol{x})$, and check that they coincide with the monodromies around faces and the zig-zags, prescribed by the rules to compute the dimer partition function.

\section{Mutations and discrete flows}\label{ss:automorph}


\subsection{Equivalence of bipartite graphs
\label{ss:equgraphs}}

As shown in \cite{GK}, under certain transformations of a graph, the face partition function remains the same, up to a cluster transformation (composition of mutations and projections) of the face variables $\boldsymbol{x}$. One such transformation (fig.~\ref{fi:spider}A) just contracts two edges incident to a two-valent vertex and  keeps intact the variables on the faces.
The second one is the so-called \textit{spider move}, shown on fig.~\ref{fi:spider}C. The variables on the faces are now changed by mutation in the center of the lozenge - with the exchange matrix $\varepsilon^{ij}$, corresponding to the Poisson bracket \rf{GKbraG}.
The third one eliminates one of the edges out of two, having the same vertices (fig.~\ref{fi:spider}B). This transformation reduces the number of variables by one, and is a composition of mutation and  projection. Every graph can be reduced to a graph, called minimal, where the number of variables cannot be already reduced.
\begin{figure}[H]
\begin{center}
\begin{tikzpicture}

\draw[thick] (0,2)--(-0.2,2.2);
\draw[thick] (0,2)--(-0.2,1.8);
\node at (-0.3,2.1){$\vdots$};
\draw[thick] (2,2)--(2.2,2.2);
\draw[thick] (2,2)--(2.2,1.8);
\node at (2.3,2.1){$\vdots$};
\draw[thick] (1.2,0.2)--(0.8,-0.2);
\draw[thick] (1.2,-0.2)--(0.8,0.2);
\node at (0.7,0.1){$\vdots$};
\node at (1.3,0.1){$\vdots$};
\draw[thick] (0,2)--(2,2);
\fill (2,2) circle (0.08); \fill[color=gray!0] (2,2) circle (0.04);
\fill (0,2) circle (0.08); \fill[color=gray!0] (0,2) circle (0.04);
\fill (1,2) circle (0.08);
\fill (1,2) circle (0.08);
\fill (1,0) circle (0.08);\fill[color=gray!0] (1,0) circle (0.04);
\draw[<->,thick] (1,1.5)--(1,0.5);
\node at (1,-1.3) {A};
\end{tikzpicture}\hspace{1.5cm}
\begin{tikzpicture}
\node at (1,2) {$x$};
\node at (1,2.4) {$y$};
\node at (1,1.6) {$z$};
\node at (1.2,0.4) {$y(1+x^{-1})^{-1}$};
\node at (1,-0.4) {$z(1+x)$};

\draw[thick] (0,2)--(-0.2,2.2);
\draw[thick] (0,2)--(-0.2,1.8);
\node at (-0.3,2.1){$\vdots$};
\draw[thick] (2,2)--(2.2,2.2);
\draw[thick] (2,2)--(2.2,1.8);
\node at (2.3,2.1){$\vdots$};
\draw[thick] (0,0)--(-0.2,0.2);
\draw[thick] (0,0)--(-0.2,-0.2);
\node at (-0.3,0.1){$\vdots$};
\draw[thick] (2,0)--(2.2,0.2);
\draw[thick] (2,0)--(2.2,-0.2);
\node at (2.3,0.1){$\vdots$};
\draw[thick] (0,0)--(2,0);
\draw[thick] (0,2) ..controls(0.2,1.7) and (1.8,1.7).. (2,2);
\draw[thick] (0,2) ..controls(0.2,2.3) and (1.8,2.3).. (2,2);
\fill (0,2) circle (0.08); \fill[color=gray!0] (0,2) circle (0.04);
\fill (2,2) circle (0.08);
\fill (0,0) circle (0.08); \fill[color=gray!0] (0,0) circle (0.04);
\fill (2,0) circle (0.08);
\draw[->,thick] (1,1.3)--(1,0.7);
\node at (0.5,-1.3) {B};
\node at (3.5,-1.3) {~};
\end{tikzpicture}
\end{center}
\end{figure}
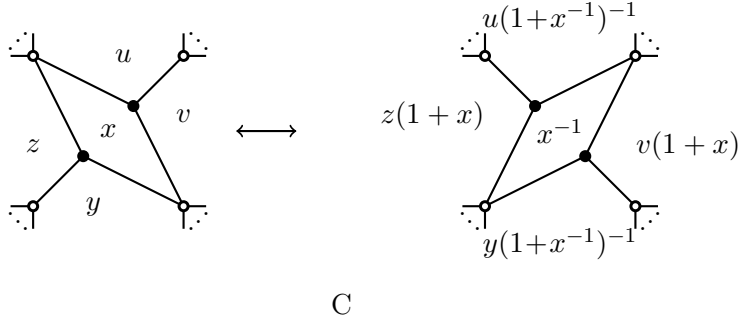
\begin{figure}[H]
\begin{center}
\begin{tikzpicture}
\node at (1,1) {$x$};
\node at (0.8,0) {$y$};
\node at (0,0.8) {$z$};
\node at (1.2,2) {$u$};
\node at (2,1.2) {$v$};
\draw[thick] (0,2.3)--(0,2)--(-0.3,2);
\node at (-0.1,2.3) {.};\node at (-0.2,2.2) {.};\node at (-0.3,2.1) {.};
\draw[thick] (0,-0.3)--(0,0)--(-0.3,0);
\node at (-0.1,-0.3) {.};\node at (-0.2,-0.2) {.};\node at (-0.3,-0.1) {.};
\draw[thick] (2.3,0)--(2,0)--(2,-0.3);
\node at (2.3,-0.1) {.};\node at (2.2,-0.2) {.};\node at (2.1,-0.3) {.};
\draw[thick] (2,2.3)--(2,2)--(2.3,2);
\node at (2.1,2.3) {.};\node at (2.2,2.2) {.};\node at (2.3,2.1) {.};
\draw[thick] (0,0)--(0.6666,0.6666);
\draw[thick] (1.333,1.333)--(2,2);
\draw[thick] (0,2)--(0.6666,0.6666)--(2,0)--(1.3333,1.3333)--(0,2);
\fill (1.333,1.333) circle (0.08);
\fill (0.666,0.666) circle (0.08);
\fill (0,2) circle (0.08); \fill[color=gray!0] (0,2) circle (0.04);
\fill (2,0) circle (0.08); \fill[color=gray!0] (2,0) circle (0.04);
\fill (2,2) circle (0.08); \fill[color=gray!0] (2,2) circle (0.04);
\fill (0,0) circle (0.08); \fill[color=gray!0] (0,0) circle (0.04);
\draw[<->,thick] (2.7,1)--(3.5,1);
\node at (4.1,-1.3) {C};
\end{tikzpicture}
\begin{tikzpicture}
\node at (1,1) {$x^{-1}$};
\node at (1,-0.5) {$y(1\!+\!x^{-1})^{-1}$};
\node at (-0.7,1.2) {$z(1+x)$};
\node at (1,2.5) {$u(1\!+\!x^{-1})^{-1}$};
\node at (2.7,0.8) {$v(1+x)$};
\draw[thick] (0,2.3)--(0,2)--(-0.3,2);
\node at (-0.1,2.3) {.};\node at (-0.2,2.2) {.};\node at (-0.3,2.1) {.};
\draw[thick] (0,-0.3)--(0,0)--(-0.3,0);
\node at (-0.1,-0.3) {.};\node at (-0.2,-0.2) {.};\node at (-0.3,-0.1) {.};
\draw[thick] (2.3,0)--(2,0)--(2,-0.3);
\node at (2.3,-0.1) {.};\node at (2.2,-0.2) {.};\node at (2.1,-0.3) {.};
\draw[thick] (2,2.3)--(2,2)--(2.3,2);
\node at (2.1,2.3) {.};\node at (2.2,2.2) {.};\node at (2.3,2.1) {.};
\draw[thick] (0,2)--(0.666,1.333);
\draw[thick] (2,0)--(1.333,0.666);
\draw[thick] (0,0)--(0.666,1.333)--(2,2)--(1.333,0.666)--(0,0);
\fill (0.666,1.333) circle (0.08);
\fill (1.333,0.666) circle (0.08);
\fill (0,2) circle (0.08); \fill[color=gray!0] (0,2) circle (0.04);
\fill (2,0) circle (0.08); \fill[color=gray!0] (2,0) circle (0.04);
\fill (2,2) circle (0.08); \fill[color=gray!0] (2,2) circle (0.04);
\fill (0,0) circle (0.08); \fill[color=gray!0] (0,0) circle (0.04);
\node[color=white] at (2,-1.3) {C};
\end{tikzpicture}
\end{center}
\caption{(A): Elimination of a two-valent vertex; (B): Elimination of a double edge; (C): Spider move.}
\label{fi:spider}
\end{figure}

In particular, if a sequence of such moves preserves the graph, but changes the variables (thus giving a nontrivial mapping class group element of the corresponding cluster variety) the corresponding map preserves the Poisson structure and the integrals of motion and thus commutes with the integrable flow.

\subsection{Discrete flow $\tau$\label{sstau}}

Recall that the transformation $\tau$ is defined on $G/\Ad H$ and on Poisson submanifolds $G^u/\mbox{Ad}H$ as by $\tau(g_-g_+)=g_+g_-$ where $g_-$ and $g_+$ belong to the lower and upper Borel subgroups, respectively. For group $\widehat{PGL}(N)$ one can generalise this transformation by $\tau(g_-g_+) = g_+\Lambda^k g_-\Lambda^{-k}$. This transformation obviously preserves the Hamiltonians of the integrable system given by the generating function $\det(g\Lambda^k-\mu)$.

Any $u\in (\widehat{W}\times \widehat{W})^\sharp$ admits a decomposition $u=s_{i_1}\cdots s_{i_l}\Lambda^k$ with $i_1,\ldots i_p$ positive and $i_{p+1}\ldots i_l$ negative\footnote{One can start from any decomposition and then using that negative and positive generators commute move the former to the right and the latter to the left.}. Such decomposition can be transformed in the following way:
$$ u=s_{i_1}\cdots s_{i_l}\Lambda^k = s_{i_1}\cdots s_{i_p}s_{i_{p+1}}\cdots s_{i_l}\Lambda^k
\sim s_{i_{p+1}}\cdots s_{i_l}\Lambda^k s_{i_1}\cdots s_{i_p}=$$ $$= s_{i_{p+1}}\cdots s_{i_l} s_{i_1+k}\cdots s_{i_p+k}\Lambda^k=s_{i_1+k}\cdots s_{i_p+k}s_{i_{p+1}}\cdots s_{i_l}\Lambda^k =u',
$$
where $\sim$ means a cyclic shift of the generators, the second and the fourth equality follow from the commutativity of negative and positive generators and the third follows from the relation $\Lambda s_i= s_{i+1}\Lambda$. This sequence of equalities and cyclic shifts induces a morphism from $G^u/\mbox{Ad}H$ to $G^{u'}/\mbox{Ad}H$. In particular it is clear that if $k=0$ than $u=u'$ and we get the automorphism of $G^u/\mbox{Ad}H$ preserving the integrals of motion. This automorphism is just the transformation $\tau$. If $k\neq 0$ the same transformation, repeated $N/\mbox{gcd}(k,N)$ times, gives an automorphism, which we will also denote by $\tau$.

\subsection{General discrete flows}

Take a Newton polygon $\Delta$ and consider a sequence of the corresponding Thurston diagrams related by Thurston moves and with the first and the last diagram coinciding. Such sequence defines a possibly trivial automorphism of the corresponding cluster variety $\mathcal{X}_\Delta$ (see Appendix~\ref{ap:cluster}), preserving the integrable system of sect.~\ref{ss:dimerstoloop}. Such automorphisms form a group $\mathcal{G}_\Delta$ of discrete automorphisms of our integrable system. Though we do not have a complete description of this group so far, we are going to present a combinatorial construction of a discrete Abelian group $\mathcal{G}^\prime_\Delta$, where $\Delta$ is a Newton polygon, and define a homomorphism of the group of discrete automorphisms of the corresponding integrable system into $\mathcal{G}^\prime_\Delta$. We conjecture that this homomorphism is in fact an isomorphism.

Let $\mathbb{Z}^{|V_\Delta|}$ be the group  of integral valued functions on the set of vertices $V_\Delta$ of the Newton polygon $\Delta$. Let $A$ be the subgroup of functions on $V_\Delta$ extendible to affine functions with integral coefficients, i.e. the functions  $aX+bY+c$ with $a,b,c\in \mathbb{Z}$ and $(X,Y)$ are coordinates on $H_1(\Sigma)=\mathbb{Z}^2$. Our group is defined as the quotient $\mathcal{G}'_\Delta=\mathbb{Z}^{|V_\Delta|}/A$.
It is obvious, that this group has rank equal to $|V_\Delta|-3$, but in general it has a nontrivial torsion.

Now we are going to describe the homomorphism $\mathcal{G}_\Delta\to\mathcal{G}^\prime_\Delta$. Take a sequence of Thurston diagrams related by Thurston moves representing an element of $\mathcal{G}_\Delta$.
Any such sequence can be interpolated by a one parameter family of collections of curves such, that for all but finitely many parameter values, all intersections are triple, and for each remaining value we allow one quadruple intersection corresponding to the moves. Since the first and the last diagrams coincide we can assume that the parameter belongs to a unit circle $S^1$. Consider a product $\Sigma\times S^1$. Every curve of the Thurston diagram traces a surface in this product. If $\Sigma$ is a torus $T^2$ then $\Sigma\times S^1$ is a three-dimensional torus $T^3$ and every cycle defines an element of the homology group $H_2(T^3,\mathbb{Z})=\mathbb{Z}^3$. Therefore we get a collection $\delta_1,\ldots,\delta_K$ of points in $H_2(T^3,\mathbb{Z})$, one for each side of the Newton polygon $\Delta$. To be specific assume, that the indices enumerate the cycles in the counter-clockwise order of the sides of $\Delta$, starting from the rightmost vertex. The sum $\sum_r\delta_r$ of the cycles of such collection vanishes in homology. Indeed, the connected components of the complement to the cycles inherit colouring in grey and white from the Thurston diagrams. Therefore the sum  of the cycles is the boundary of the union of the, say, grey components.

Choose a homology basis $e^1,e^2,e^3$ in $H_1(T^3,\mathbb{Z})$ with $e^1$ and $e^2$ being the basis $H_1(T^2,\mathbb{Z})$. This basis induces a basis $e^{kl}=e^k\wedge e^l$, where $k<l$ of $H_2(T^3,\mathbb{Z})$. Denote by $(a^i_{12},a^i_{23},a^i_{13})$ the coordinates of the classes $\delta_i$ in this basis.  Since $\sum a^i_{12}=0$, there is a collection of numbers $\{f^i\}$ on vertices of $\Delta$, defined up to a constant and such that $a^i_{12}=f^{i+1}-f^{i}$ for $i\in \mathbb{Z}/K\mathbb{Z}$. This collection defines the element of $\mathcal{G}^\prime_\Delta$.

Observe that if a curve of a Thurston diagram represents the cycle $Xe^1+Ye^2$ and it moves by $ae^1+be^2$ then the corresponding cycle on $T^3$ is given by $(Xe^1+Ye^2)\wedge (ae^1+be^2+e^3)=(bX-aY)e^{12}+Xe^{13}+Ye^{23}$. Consider the case when the shift $ae^1+be^2$ is the same for all curves. If $a$ and $b$ are integers, the corresponding automorphism is trivial and the corresponding element of the group $\{f_i\}\in \mathcal{G}^\prime_\Delta$ is also trivial as given by an affine function $bX-aY$. This argument shows that the construction of the homomorphism sounds. On the other hand if $X$ and $Y$ are rational but all coordinates $a_i^{kl}$ are integers then the corresponding element of the group $\mathcal{G}^\prime_\Delta$ is a torsion.

Now consider an automorphism $\tau$ defined in sect.~\ref{sstau}. In the language of Thurston diagrams it corresponds to moving the collection of cycles oriented from right to left as a whole and keeping the remaining cycles on their places. In order to come back to the original configuration of cycles we need to move the collection by the vector $(Ne_1+ke_2)/\mbox{gcd}(k,N)$. The vector $Ne_1+ke_2$ is just the vector connecting the leftmost and the rightmost vertices of the Newton polygon $\Delta$. Such vector cuts the Newton polygon in two parts. The corresponding function $f_i$ can be taken to be zero on the upper part and a minimal nonzero linear function vanishing on the vector and taking integral values on integral points on the lower part.

As it was indicated to us by A.~B.~Goncharov, the group $\mathcal{G}_\Delta$ has another interpretation. Namely, a Newton polygon defines a toric surface $M_\Delta$ and a homology class on it represented by the spectral curves (see \cite{Fulton}, sect.~3.4).  The sides of the polygon correspond to the components of the divisor at infinity. The Picard group on $M_\Delta$ restricts to the Picard group of a generic spectral curve. The group $\mathcal{G}_\Delta$ is the inverse image of lines bundles of degree zero. In other terms the group $\mathcal{G}_\Delta$ is the group of line bundles on a spectral curve extendible to line bundles on $M_\Delta$. We conjecture that the discrete flows act on the Jacobian of the spectral curve by multiplication by the corresponding line bundles.

\section{Examples
\label{ss:examples}}

\subsection{Triangle}

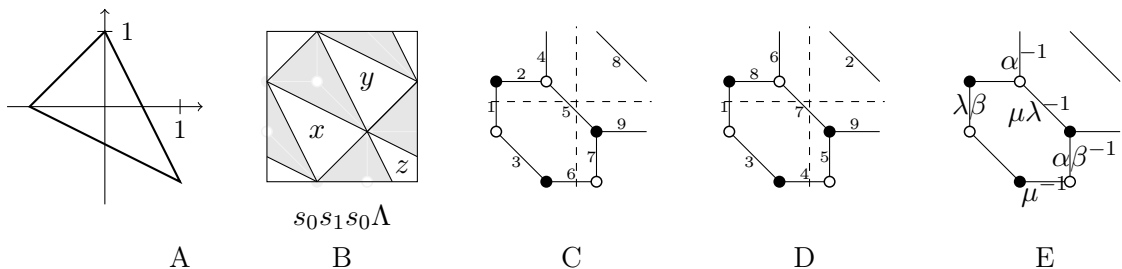
\begin{figure}[H]
\begin{center}
\begin{tikzpicture}
\draw[->](-1.3,0)--(1.3,0);
\draw[->](0,-1.3)--(0,1.3);
\draw (1,0.1)--(1,-0.1);\node at (1,-0.3) {\small 1};
\draw (0.1,1)--(-0.1,1);\node at (0.3,1) {\small 1};
\draw[thick] (1,-1)--(0,1)--(-1,0)--cycle;
\node at (1,-2) {A};
\end{tikzpicture}\qquad
\begin{tikzpicture}[scale=0.6666]
\fill[color=gray!20] (0,2)--(1,3)--(2,1)--cycle;
\fill[color=gray!20] (0,2)--(1,0)--(0,0.5)--cycle;
\fill[color=gray!20] (3,2)--(1,3)--(2.5,3)--cycle;
\fill[color=gray!20] (1,0)--(2.5,0)--(2,1)--cycle;
\fill[color=gray!20] (3,2)--(2,1)--(3,0.5)--cycle;
\draw[color=gray!10] (0,1)--(1,0)--(2,0)--(2,1)--(1,2)--(0,2)--cycle;
\draw[color=gray!10] (2,1)--(3,1);
\draw[color=gray!10] (1,2)--(1,3);
\draw[color=gray!10] (2,3)--(3,2);
\fill[color=gray!10] (2,1) circle (0.12);
\fill[color=gray!10] (0,2) circle (0.12);
\fill[color=gray!10] (1,0) circle (0.12);
\fill[color=gray!10] (1,2) circle (0.12);\fill[color=white] (1,2) circle (0.09);
\fill[color=gray!10] (2,0) circle (0.12);\fill[color=white] (2,0) circle (0.09);
\fill[color=gray!10] (0,1) circle (0.12);\fill[color=white] (0,1) circle (0.09);
\draw[thin] (3,0)--(0,0)--(0,3)--(3,3)--cycle;
\draw (0,2)--(1,0);
\draw (0,2)--(3,0.5);
\draw (1,3)--(2.5,0);
\draw (1,0)--(3,2);
\draw (2,1)--(3,2);
\draw (0,0.5)--(1,0);
\draw (0,2)--(1,3);
\draw (3,2)--(2.5,3);
\draw (1,3)--(3,2);
\node at (1,1) {$x$};
\node at (2,2) {$y$};
\node at (2.7,0.3) {$z$};
\node at (1.5,-0.7) {$s_0s_1s_0\Lambda$};
\node at (1.5,-1.5) {B};
\end{tikzpicture}\qquad
\begin{tikzpicture}[scale=0.6666]
\draw (0,1)--(1,0)--(2,0)--(2,1)--(1,2)--(0,2)--cycle;
\draw (2,1)--(3,1);
\draw (1,2)--(1,3);
\draw (2,3)--(3,2);
\fill (2,1) circle (0.12);
\fill (0,2) circle (0.12);
\fill (1,0) circle (0.12);
\fill (1,2) circle (0.12);\fill[color=white] (1,2) circle (0.09);
\fill (2,0) circle (0.12);\fill[color=white] (2,0) circle (0.09);
\fill (0,1) circle (0.12);\fill[color=white] (0,1) circle (0.09);
\node at (1.5,-1.5) {C};
\node at (-0.1,1.5) {\tiny 1};
\node at (0.5,2.15) {\tiny 2};
\node at (0.4,0.4) {\tiny 3};
\node at (1.4,1.4) {\tiny 5};
\node at (0.9,2.5) {\tiny 4};
\node at (1.5,0.15) {\tiny 6};
\node at (1.9,0.5) {\tiny 7};
\node at (2.4,2.4) {\tiny 8};
\node at (2.5,1.15) {\tiny 9};
\draw[dashed] (1.6,-0.1)--(1.6,3.1);
\draw[dashed] (-0.1,1.6)--(3.1,1.6);
\node at (1.5,-0.7) {~};
\end{tikzpicture}\qquad
\begin{tikzpicture}[scale=0.6666]
\draw (0,1)--(1,0)--(2,0)--(2,1)--(1,2)--(0,2)--cycle;
\draw (2,1)--(3,1);
\draw (1,2)--(1,3);
\draw (2,3)--(3,2);
\fill (2,1) circle (0.12);
\fill (0,2) circle (0.12);
\fill (1,0) circle (0.12);
\fill (1,2) circle (0.12);\fill[color=white] (1,2) circle (0.09);
\fill (2,0) circle (0.12);\fill[color=white] (2,0) circle (0.09);
\fill (0,1) circle (0.12);\fill[color=white] (0,1) circle (0.09);
\node at (1.5,-1.5) {D};
\node at (-0.1,1.5) {\tiny 1};
\node at (0.5,2.15) {\tiny 8};
\node at (0.4,0.4) {\tiny 3};
\node at (1.4,1.4) {\tiny 7};
\node at (0.9,2.5) {\tiny 6};
\node at (1.5,0.15) {\tiny 4};
\node at (1.9,0.5) {\tiny 5};
\node at (2.4,2.4) {\tiny 2};
\node at (2.5,1.15) {\tiny 9};
\draw[dashed] (1.6,-0.1)--(1.6,3.1);
\draw[dashed] (-0.1,1.6)--(3.1,1.6);
\node at (1.5,-0.7) {~};
\end{tikzpicture}\qquad
\begin{tikzpicture}[scale=0.6666]
\draw (0,1)--(1,0)--(2,0)--(2,1)--(1,2)--(0,2)--cycle;
\draw (2,1)--(3,1);
\draw (1,2)--(1,3);
\draw (2,3)--(3,2);
\fill (2,1) circle (0.12);
\fill (0,2) circle (0.12);
\fill (1,0) circle (0.12);
\fill (1,2) circle (0.12);\fill[color=white] (1,2) circle (0.09);
\fill (2,0) circle (0.12);\fill[color=white] (2,0) circle (0.09);
\fill (0,1) circle (0.12);\fill[color=white] (0,1) circle (0.09);
\node at (1.5,-1.5) {E};
\node at (0,1.5) {$\lambda\beta$};
\node at (1.4,1.4) {$\mu\lambda^{-1}$};
\node at (1.1,2.5) {$\alpha^{-1}$};
\node at (1.5,-0.2) {$\mu^{-1}$};
\node at (2.3,0.5) {$\alpha\beta^{-1}$};
\node at (1.5,-0.7) {~};
\end{tikzpicture}
\end{center}
\caption{Triangle. (A): Newton polygon $\Delta$; (B): Thurston diagram with face variables; (C): The bipartite graph $\Gamma$; (D): Bipartite graph $\Gamma$ on dual torus; (E): Edge variables and cycles defining the quadratic form.}
\label{fi:hexagon}
\end{figure}

Consider the Newton polygon (just triangle in this case) with vertices $(-1,0),(0,1)$ and $(1,-1)$ of the area $\mathsf{S} = 3/2$, shown on figure \ref{fi:hexagon}A. The corresponding Thurston diagram is drawn on fig.~\ref{fi:hexagon}B, its lines are parallel to the sides of the triangle. Three independent triple
intersections in the fundamental domain, if sides of the triangle are oriented counter-clockwise, can be put into correspondence with the generators $s_0$, $s_1$ and, again, $s_0$ respectively, so that the
corresponding Poisson submanifold is described by the word $u=s_0s_1s_0\Lambda$, an element of the
coextended Weyl group of $\hat{A}_1$.

The corresponding trivalent graph $\Gamma$ is shown on fig. \ref{fi:hexagon}C. The dual surface $\widetilde{\Sigma}$ is also a torus in this example, with isomorphic embeddings of the graph.
The Poisson bracket between the face variables $\{x,y,z: xyz=1\}$ follow from the fig.~\ref{fi:hexagon}B (by the rule of
construction of the exchange matrix from fig.~\ref{fi:thurston}B of Appendix~\ref{ap:thurston}):
$$
\{x,y\}=3xy;\quad \{y,z\}=3yz;\quad \{z,x\}=3yz
$$
A section $H^1(\Sigma)\times H^1(\tilde{\Sigma})\to C^1(\Gamma)$ is shown on fig.~\ref{fi:hexagon}E. Here $\lambda,\mu$ are coordinates on $H^1(\Sigma)$ and $\alpha,\beta$ are coordinates on $H^1(\tilde{\Sigma})$, the space $\mathcal{C}^1$ is trivial in this example.
On fig.~\ref{fi:hexagon}C and fig.~\ref{fi:hexagon}D the cycles in homology bases of $\Sigma$ and $\tilde{\Sigma}$ are shown respectively, and the variables on the edges just indicate the intersections of these cycles and the edges of the graph. We associate to every variable a dashed cycle and put it on the edges, intersected by this cycle.

\begin{figure}[H]
\begin{center}
\begin{tikzpicture}[scale=0.6666]
\draw (0,1)--(1,0)--(2,0)--(2,1)--(1,2)--(0,2)--cycle;
\draw (2,1)--(3,1);
\draw (1,2)--(1,3);
\draw (2,3)--(3,2);
\draw[ultra thick] (2,1)--(3,1);
\draw[ultra thick] (1,0)--(2,0);
\draw[ultra thick] (0,2)--(1,2);
\fill (2,1) circle (0.12);
\fill (0,2) circle (0.12);
\fill (1,0) circle (0.12);
\fill (1,2) circle (0.12);\fill[color=white] (1,2) circle (0.09);
\fill (2,0) circle (0.12);\fill[color=white] (2,0) circle (0.09);
\fill (0,1) circle (0.12);\fill[color=white] (0,1) circle (0.09);
\node at (1.5,-0.7) {$\mu^{-1}$};
\end{tikzpicture}\quad
\begin{tikzpicture}[scale=0.6666]
\draw (0,1)--(1,0)--(2,0)--(2,1)--(1,2)--(0,2)--cycle;
\draw (2,1)--(3,1);
\draw (1,2)--(1,3);
\draw (2,3)--(3,2);
\draw[ultra thick] (1,2)--(1,3);
\draw[ultra thick] (2,0)--(2,1);
\draw[ultra thick] (0,1)--(0,2);
\fill (2,1) circle (0.12);
\fill (0,2) circle (0.12);
\fill (1,0) circle (0.12);
\fill (1,2) circle (0.12);\fill[color=white] (1,2) circle (0.09);
\fill (2,0) circle (0.12);\fill[color=white] (2,0) circle (0.09);
\fill (0,1) circle (0.12);\fill[color=white] (0,1) circle (0.09);
\node at (1.5,-0.7) {$\lambda$};
\end{tikzpicture}\quad
\begin{tikzpicture}[scale=0.6666]
\draw (0,1)--(1,0)--(2,0)--(2,1)--(1,2)--(0,2)--cycle;
\draw (2,1)--(3,1);
\draw (1,2)--(1,3);
\draw (2,3)--(3,2);
\draw[very thick] (0,1)--(1,0);
\draw[very thick] (2,1)--(1,2);
\draw[very thick] (3,2)--(2,3);
\fill (2,1) circle (0.12);
\fill (0,2) circle (0.12);
\fill (1,0) circle (0.12);
\fill (1,2) circle (0.12);\fill[color=white] (1,2) circle (0.09);
\fill (2,0) circle (0.12);\fill[color=white] (2,0) circle (0.09);
\fill (0,1) circle (0.12);\fill[color=white] (0,1) circle (0.09);
\node at (1.5,-0.7) {$-\mu\lambda^{-1}$};
\end{tikzpicture}\quad
\begin{tikzpicture}[scale=0.6666]
\draw (0,1)--(1,0)--(2,0)--(2,1)--(1,2)--(0,2)--cycle;
\draw (2,1)--(3,1);
\draw (1,2)--(1,3);
\draw (2,3)--(3,2);
\draw[ultra thick] (0,1)--(1,0);
\draw[ultra thick] (2,0)--(2,1);
\draw[ultra thick] (1,2)--(0,2);
\fill (2,1) circle (0.12);
\fill (0,2) circle (0.12);
\fill (1,0) circle (0.12);
\fill (1,2) circle (0.12);\fill[color=white] (1,2) circle (0.09);
\fill (2,0) circle (0.12);\fill[color=white] (2,0) circle (0.09);
\fill (0,1) circle (0.12);\fill[color=white] (0,1) circle (0.09);
\node at (1.5,-0.7) {$\alpha\beta^{-1}$};
\end{tikzpicture}\quad
\begin{tikzpicture}[scale=0.6666]
\draw (0,1)--(1,0)--(2,0)--(2,1)--(1,2)--(0,2)--cycle;
\draw (2,1)--(3,1);
\draw (1,2)--(1,3);
\draw (2,3)--(3,2);
\draw[ultra thick] (1,2)--(1,3);
\draw[ultra thick] (2,1)--(3,1);
\draw[ultra thick] (2,3)--(3,2);
\fill (2,1) circle (0.12);
\fill (0,2) circle (0.12);
\fill (1,0) circle (0.12);
\fill (1,2) circle (0.12);\fill[color=white] (1,2) circle (0.09);
\fill (2,0) circle (0.12);\fill[color=white] (2,0) circle (0.09);
\fill (0,1) circle (0.12);\fill[color=white] (0,1) circle (0.09);
\node at (1.5,-0.7) {$\alpha^{-1}$};
\end{tikzpicture}\quad
\begin{tikzpicture}[scale=0.6666]
\draw (0,1)--(1,0)--(2,0)--(2,1)--(1,2)--(0,2)--cycle;
\draw (2,1)--(3,1);
\draw (1,2)--(1,3);
\draw (2,3)--(3,2);
\draw[ultra thick] (1,2)--(2,1);
\draw[ultra thick] (0,1)--(0,2);
\draw[ultra thick] (1,0)--(2,0);
\fill (2,1) circle (0.12);
\fill (0,2) circle (0.12);
\fill (1,0) circle (0.12);
\fill (1,2) circle (0.12);\fill[color=white] (1,2) circle (0.09);
\fill (2,0) circle (0.12);\fill[color=white] (2,0) circle (0.09);
\fill (0,1) circle (0.12);\fill[color=white] (0,1) circle (0.09);
\node at (1.5,-0.7) {$\beta$};
\end{tikzpicture}\quad
\end{center}
\caption{Dimer configuration and corresponding pairings $(-1)^{Q(D)}\langle \boldsymbol{A},D\rangle$.}
\label{fi:hexagon-dimers}
\end{figure}
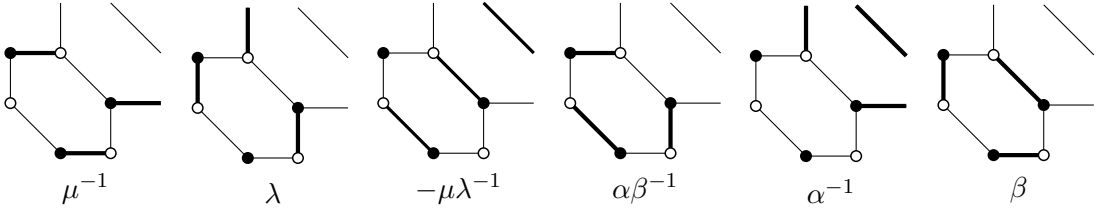

There are six dimer configurations on this graph shown on fig.~\ref{fi:hexagon-dimers},
\noindent
which sum up in the following the double partition function
\begin{equation}\label{Sftriang}
S^d=\mu^{-1}+\lambda-\mu\lambda^{-1}+\alpha^{-1}+\beta+\alpha\beta^{-1}
\end{equation}
The same expression can be obtained as the determinant $S^d=\det\mathfrak{D}$ of the Dirac operator
$$
\mathfrak{D}=\begin{pmatrix}1&1&-\lambda\beta\\
\mu\lambda^{-1}&\alpha^{-1}&1 \\
\alpha\beta^{-1}&-\mu^{-1}&1
\end{pmatrix}
$$
which can be read off the fig.~\ref{fi:hexagon}E (We have chosen the Kasteleyn marking in a way to reproduce the signs exactly and not up to equivalence.)

Solving the equations relating $H^1(\tilde{\Sigma})$ with $B^2(\Sigma)$ (one should take the products
around faces on fig.~\ref{fi:hexagon}E taking into account canonical orientation of the edges)
$$
x=\alpha\beta^{-2},\quad y=\alpha\beta,\quad z=\alpha^{-2}\beta
$$
and substituting the result for $\alpha$ and $\beta$ into the expression (\ref{Sftriang}), one gets the normalised face partition function:
$$
S^f = \lambda^{-1}+\mu-\lambda\mu^{-1}+ x^{-1/3}y^{-2/3} + x^{-1/3}y^{1/3} + x^{2/3}y^{1/3}
$$
On the other hand the characteristic polynomial of the corresponding matrix product is
\be
\mathcal{A}(\lambda,\mu;x,y,z) = \det\left(H_0(x)E_0H_1(y)E_1H_0(z)E_0\Lambda-\mu\right)=
\\
= \det\left[
\begin{pmatrix}1&0\\\lambda x&1\end{pmatrix}
\begin{pmatrix}y&0\\0&1\end{pmatrix}
\begin{pmatrix}1&1\\0&1\end{pmatrix}
T_xT_yT_z
\begin{pmatrix}1&0\\\lambda&1\end{pmatrix}
\begin{pmatrix}0&1\\\lambda&0\end{pmatrix}
-\mu \right]=
\\
=\mu^2 - \lambda y - \mu\lambda^2xy - \lambda\mu(1+y+xy)
\ee
Normalising by the substitution $\lambda\mapsto \lambda x^{-2/3}y^{-1/3}$, $\mu \mapsto \mu x^{-1/3}y^{1/3}$ and multiplying the result by $-x^{2/3}y^{-2/3}\lambda^{-1}\mu^{-1}$ one obtains the same expression for the partition function.

The group $\mathcal{G}^\prime_\Delta=\mathcal{G}_\Delta$ in this case is just $\mathbb{Z}/3\mathbb{Z}$ generated by the obvious cyclic permutation of variables $(x,y,z)\to(y,z,x)$ which corresponds to the element $(0,1,0)$ of the group $\mathcal{G}^\prime_\Delta$.

\subsection{The simplest relativistic Toda chain}

In our next example the rank of the Poisson structure is still two, but in contrast to the case of triangle there is a nontrivial Casimir function.

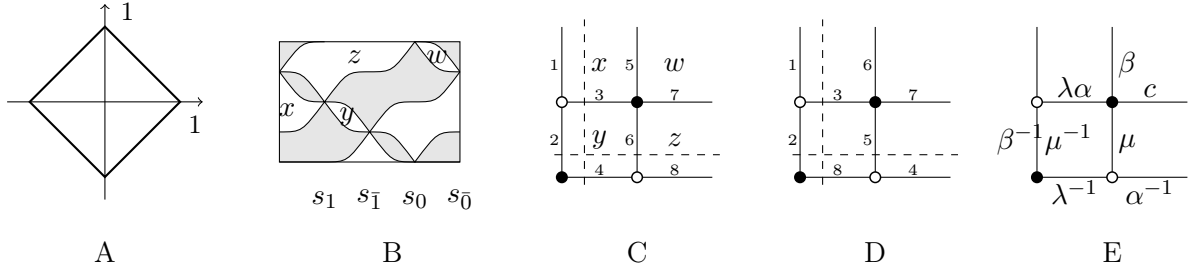
\begin{figure}[H]
\begin{center}
\begin{tikzpicture}
\draw[->](-1.3,0)--(1.3,0);
\draw[->](0,-1.3)--(0,1.3);
\node at (1.2,-0.3) {\small 1};
\node at (0.3,1.2) {\small 1};
\draw[thick] (1,0)--(0,1)--(-1,0)--(0,-1)--cycle;
\node at (0,-2) {A};
\end{tikzpicture}\qquad
\begin{tikzpicture}
\fill[color=gray!20] (2.4,1.6)--(2.4,1.2)..controls(2.1,1.6)..(1.8,1.6)--cycle;
\draw                (2.4,1.2)..controls(2.1,1.6)..(1.8,1.6);
\fill[color=gray!20] (0,1.2)..controls(0.3,1.6)..(0.6,1.6)--(0,1.6)--cycle;
\draw                (0,1.2)..controls(0.3,1.6)..(0.6,1.6);
\fill[color=gray!20] (0,1.2)..controls(0.3,1.2)..(0.6,0.8)..controls(0.3,0.8)..(0,1.2)--cycle;
\draw                (0,1.2)..controls(0.3,1.2)..(0.6,0.8)..controls(0.3,0.8)..(0,1.2);
\fill[color=gray!20] (1.2,0.4)..controls(1.5,0.4)..(1.8,0)..controls(1.5,0)..(1.2,0.4)--cycle;
\draw                (1.2,0.4)..controls(1.5,0.4)..(1.8,0)..controls(1.5,0)..(1.2,0.4);
\fill[color=gray!20] (2.4,0.4)..controls(2.1,0.4)..(1.8,0)--(2.4,0)--cycle;
\draw                (2.4,0.4)..controls(2.1,0.4)..(1.8,0)--(2.4,0);
\fill[color=gray!20] (2.4,1.2)..controls(2.1,1.2)..(1.8,1.6)..controls(1.5,1.2)..(1.2,1.2)..controls(0.9,1.2)..(0.6,0.8)..controls(0.9,0.8)..(1.2,0.4)..controls(1.5,0.8)..(1.8,0.8)..controls(2.1,0.8)..(2.4,1.2)--cycle;
\draw (2.4,1.2)..controls(2.1,1.2)..(1.8,1.6)..controls(1.5,1.2)..(1.2,1.2)..controls(0.9,1.2)..(0.6,0.8)..controls(0.9,0.8)..(1.2,0.4)..controls(1.5,0.8)..(1.8,0.8)..controls(2.1,0.8)..(2.4,1.2);
\fill[color=gray!20] (0,0)--(0.6,0)..controls (0.9,0)..(1.2,0.4)..controls(0.9,0.4)..(0.6,0.8)..controls(0.3,0.4)..(0,0.4)--cycle;
\draw                (0,0)--(0.6,0)..controls (0.9,0)..(1.2,0.4)..controls(0.9,0.4)..(0.6,0.8)..controls(0.3,0.4)..(0,0.4);
\draw[thin] (0,0)--(2.4,0)--(2.4,1.6)--(0,1.6)--cycle;
\node at (0.1,0.7) {$x$};
\node at (1,1.4) {$z$};
\node at (0.9,0.6) {$y$};
\node at (2.1,1.4) {$w$};
\node at (0.6,-0.5) {$s_1$};
\node at (1.2,-0.5) {$s_{\bar{1}}$};
\node at (1.8,-0.5) {$s_0$};
\node at (2.4,-0.5) {$s_{\bar{0}}$};
\node at (1.5,-1.2) {B};
\end{tikzpicture}\qquad
\begin{tikzpicture}
\draw (0,0) -- +(0,2);
\draw (0,0) -- +(2,0);
\draw (0,0)+(0,1) -- +(2,1);
\draw (0,0)+(1,0) -- +(1,2);
\fill (0,0)+(0,0) circle (0.08);
\fill (0,0)+(0,1) circle (0.08);\fill[color=white] (0,1) circle (0.06);
\fill (0,0)+(1,0) circle (0.08);\fill[color=white] (1,0) circle (0.06);
\fill (0,0)+(1,1) circle (0.08);
\node at (1,-1) {C};
\draw[dashed] (-0.1,0.3)--(2.1,0.3);
\draw[dashed] (0.3,-0.1)--(0.3,2.1);
\node at (1.5,1.5) {$w$};
\node at (0.5,1.5) {$x$};
\node at (1.5,0.5) {$z$};
\node at (0.5,0.5) {$y$};
\node at (-0.1,1.5) {\tiny 1};
\node at (-0.1,0.5) {\tiny 2};
\node at (0.5,1.1) {\tiny 3};
\node at (0.5,0.1) {\tiny 4};
\node at (0.9,1.5) {\tiny 5};
\node at (0.9,0.5) {\tiny 6};
\node at (1.5,1.1) {\tiny 7};
\node at (1.5,0.1) {\tiny 8};
\end{tikzpicture}\qquad
\begin{tikzpicture}
\draw (0,0) -- +(0,2);
\draw (0,0) -- +(2,0);
\draw (0,0)+(0,1) -- +(2,1);
\draw (0,0)+(1,0) -- +(1,2);
\fill (0,0)+(0,0) circle (0.08);
\fill (0,0)+(0,1) circle (0.08);\fill[color=white] (0,1) circle (0.06);
\fill (0,0)+(1,0) circle (0.08);\fill[color=white] (1,0) circle (0.06);
\fill (0,0)+(1,1) circle (0.08);
\draw[dashed] (-0.1,0.3)--(2.1,0.3);
\draw[dashed] (0.3,-0.1)--(0.3,2.1);
\node at (1,-1) {D};
\node at (-0.1,1.5) {\tiny 1};%
\node at (-0.1,0.5) {\tiny 2};%
\node at (0.5,1.1) {\tiny 3};%
\node at (0.5,0.1) {\tiny 8};
\node at (0.9,1.5) {\tiny 6};
\node at (0.9,0.5) {\tiny 5};
\node at (1.5,1.1) {\tiny 7};%
\node at (1.5,0.1) {\tiny 4};
\end{tikzpicture}\quad
\begin{tikzpicture}
\draw (0,0) -- +(0,2);
\draw (0,0) -- +(2,0);
\draw (0,0)+(0,1) -- +(2,1);
\draw (0,0)+(1,0) -- +(1,2);
\fill (0,0)+(0,0) circle (0.08);
\fill (0,0)+(0,1) circle (0.08);\fill[color=white] (0,1) circle (0.06);
\fill (0,0)+(1,0) circle (0.08);\fill[color=white] (1,0) circle (0.06);
\fill (0,0)+(1,1) circle (0.08);
\node at (1,-1) {E};
\node at (0.1,0.5) {$\beta^{-1}\mu^{-1}$};
\node at (0.5,1.2) {$\lambda\alpha$};
\node at (0.5,-0.2) {$\lambda^{-1}$};
\node at (1.2,1.5) {$\beta$};
\node at (1.2,0.45) {$\mu$};
\node at (1.5,-0.2) {$\alpha^{-1}$};
\node at (1.5,1.15) {$c$};
\end{tikzpicture}
\end{center}
\caption{Two-particle relativistic Toda chain. (A): Newton polygon; (B): Thurston diagram with face variables; (C):  Bipartite graph with contracted two-valent vertices, (D): Bipartite graph on dual torus $\widetilde{\Sigma}$; (E): Edge variables and cycles defining the quadratic form.}
\label{fi:toda2poly}
\end{figure}

Consider the Newton polygon with the vertices $(1,0)$, $(0,1)$, $(-1,0)$, $(0,-1)$, shown on fig.~\ref{fi:toda2poly}A, the corresponding Thurston diagram is depicted at
fig.~\ref{fi:toda2poly}B.
The Poisson brackets between the face variables follow from fig.~\ref{fi:toda2poly}B (again, by the general rule, shown on fig.~\ref{fi:thurston}B in Appendix~\ref{ap:thurston}):
$$
\{x,y\}=2xy;\quad \{y,z\}=2yz;\quad \{z,w\}=2zw;\quad \{w,x\}=2wx;\quad \{x,z\}=\{y,w\}=0
$$
The embedding of $H^1(\Sigma)=\{(\lambda,\mu)\}$ into $H^1(\Gamma)$ is given, for example, by $A_3=\lambda$, $A_4=\lambda^{-1}$, $A_6=\mu$, $A_2=\mu^{-1}$ and all others $A_e=1$ (cf. figs.~\ref{fi:toda2poly}C and \ref{fi:toda2poly}E). Similarly the  embedding of $H^1(\tilde{\Sigma})=\{(\alpha,\beta)\}$ is given by $A_3=\alpha$, $A_8=\alpha^{-1}$, $A_5=\beta$, $A_2=\beta$ and all others $A_e=1$, where the embedding can be read off the fig.~\ref{fi:toda2poly}D. The Casimir space $\mathcal{C}^1$ is one dimensional, we denote the Casimir function (coupling constant in Toda chain) by $c$ and put on the edge 7, $A_7=c$.

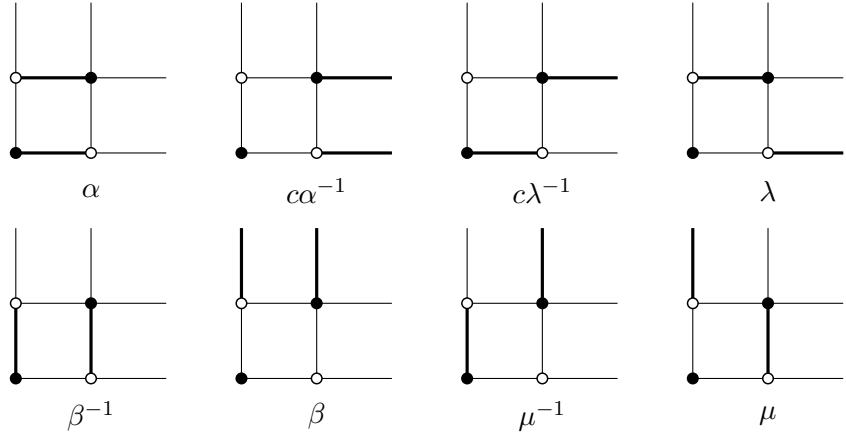
\begin{figure}[H]
\begin{center}
\begin{tikzpicture}
\draw[very thick](0,3)+(0,0)-- +(1,0);
\draw[very thick](0,3)+(0,1)-- +(1,1);
\node at (1,-0.5) {$\beta^{-1}$};
\draw[very thick](3,3)+(1,0)-- +(2,0);
\draw[very thick](3,3)+(1,1)-- +(2,1);
\node at (4,-0.5) {$\beta$};
\draw[very thick](6,3)+(0,0)-- +(1,0);
\draw[very thick](6,3)+(1,1)-- +(2,1);
\node at (7,-0.5) {$\mu^{-1}$};
\draw[very thick](9,3)+(1,0)-- +(2,0);
\draw[very thick](9,3)+(0,1)-- +(1,1);
\node at (10,-0.5) {$\mu$};
\draw[very thick](0,0)+(0,0)-- +(0,1);
\draw[very thick](0,0)+(1,0)-- +(1,1);
\node at (1,2.5) {$\alpha$};
\draw[very thick](3,0)+(0,1)-- +(0,2);
\draw[very thick](3,0)+(1,1)-- +(1,2);
\node at (4,2.5) {$c\alpha^{-1}$};
\draw[very thick](6,0)+(0,0)-- +(0,1);
\draw[very thick](6,0)+(1,1)-- +(1,2);
\node at (7,2.5) {$c\lambda^{-1}$};
\draw[very thick](9,0)+(0,1)-- +(0,2);
\draw[very thick](9,0)+(1,0)-- +(1,1);
\node at (10,2.5) {$\lambda$};

\foreach \x in {0,3,...,9}
{
   \foreach \y in {0,3}
        {
             \draw (\x,\y) -- +(0,2);
             \draw (\x,\y) -- +(2,0);
             \draw (\x,\y)+(0,1) -- +(2,1);
             \draw (\x,\y)+(1,0) -- +(1,2);
             \fill (\x,\y)+(0,0) circle (0.08);
             \fill (\x,\y)+(0,1) circle (0.08);\fill[color=white] (\x,\y)+(0,1) circle (0.06);
             \fill (\x,\y)+(1,0) circle (0.08);\fill[color=white] (\x,\y)+(1,0) circle (0.06);
             \fill (\x,\y)+(1,1) circle (0.08);
	}
}
\end{tikzpicture}
\end{center}
\caption{Dimer configurations for the graph shown on figure \ref{fi:toda2poly}C.}
\label{fi:dimers-sl2-2}
\end{figure}

There are eight dimer configurations on this graph depicted at fig.~\ref{fi:dimers-sl2-2},
\noindent
which give rise to the following double partition function~\footnote{It is interesting to point out, that already this equation in its quantum version gives rise to a nontrivial spectral problem of the
Hofstadter model \cite{Hofst}. We are planning to returm to this problem elsewhere.}
\begin{equation}
\label{Sftoda2}
 S^d=\alpha + c\alpha^{-1} + \beta + \beta^{-1} + \lambda + c\lambda^{-1}+\mu^{-1} + \mu
\end{equation}
Solving equations, relating $H^1(\tilde{\Sigma})$ and $B^2(\Sigma)$, following from fig.~\ref{fi:toda2poly}C and fig.~\ref{fi:toda2poly}E
$$
 x=\alpha\beta^{-1},\quad y=\alpha^{-1}\beta^{-1},\quad z=c\alpha^{-1}\beta,\quad w=c^{-1}\alpha\beta
$$
we change the variables to $\alpha=c^{1/2}x^{1/2}w^{1/2} = x^{1/2}y^{-1/2}$, $\beta=z^{1/2}w^{1/2}=x^{-1/2}y^{-1/2}$ since $c=xz = (yw)^{-1}$. Expressing the partition function (\ref{Sftoda2}) in terms of these variables, one gets the normalised face partition function:
\be
\label{Toda1}
S^f  
= x^{1/2}y^{-1/2} + cx^{-1/2}y^{1/2}+x^{1/2}y^{1/2}+x^{-1/2}y^{-1/2}+ \lambda + c\lambda^{-1} + \mu + \mu^{-1}
\ee
On the other hand the characteristic polynomial of the corresponding matrix product, can be
computed e.g. as
\be
\mathcal{A}(\lambda,\mu;x,y,z,w) =\det\left[H_1(x){E}_{1}H_1(y)E_{\bar{1}}H_0(z)E_{0}H_0(w)E_{\bar{0}}-\mu\right]=
\\
=\det\left[\begin{pmatrix}x&0\\ 0&1\end{pmatrix}
\begin{pmatrix}1&1\\0&1\end{pmatrix}
\begin{pmatrix}y&0\\0&1\end{pmatrix}
\begin{pmatrix}1&0\\1&1\end{pmatrix}T_xT_yT_zT_w
\begin{pmatrix}1&0\\\lambda/w&1\end{pmatrix}
\begin{pmatrix}1&\lambda^{-1}\\0&1\end{pmatrix}-\mu\right]=
\\
=\det\left[ \begin{pmatrix}xy+x+\lambda w^{-1}x&\lambda^{-1}xy+\lambda^{-1}x+xw^{-1}+x\\
1+\lambda w^{-1}&\lambda^{-1} + w^{-1}+ 1\end{pmatrix}-\mu\right]=
\\
=\mu^2 - \mu(w^{-1}+1+x+xy) - \lambda\mu w^{-1}x - \lambda^{-1}\mu  + xy
\ee
Substituting $w=c^{-1}y^{-1}$, normalising by $\lambda \to x^{-1/2}y^{-1/2}c^{-1}\lambda,$ $\mu \to -x^{1/2}y^{1/2}\mu$, and multiplying this expression by $(xy\mu)^{-1}$, one reproduces exactly the normalised dimer face partition function \rf{Toda1}.

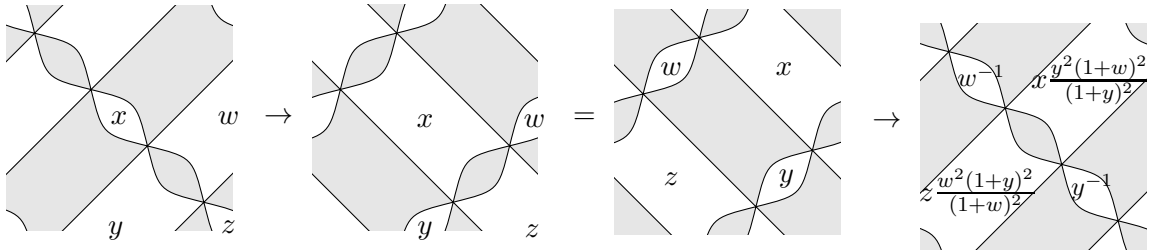
\begin{figure}[h]
 \begin{center}
\begin{tikzpicture}[scale=1.5]
\node at (2.5,1.25) {$~~~\rightarrow$};
\node at (1.25,1.25) {$x$};
\node at (2.22,0.28) {$z$};
\node at (2.22,1.25) {$w$};
\node at (1.22,0.25) {$y$};
\path[clip] (0.25,0.25) rectangle (2.25,2.25);
\foreach \x/\y in {0cm/0cm,1cm/1cm,2cm/2cm,-1cm/1cm,1cm/-1cm,0cm/2cm,2cm/0cm,4cm/2cm}
\draw[fill=gray!20,xshift=\x,yshift=\y](0.5,0)--(1.5,1)..controls(1.1,1.1)..(1,1.5)--(0,0.5)..controls(0.4,0.4)..(0.5,0);
\foreach \x/\y in {-0.5cm/0.5cm,0.5cm/-0.5cm,1.5cm/0.5cm,2.5cm/1.5cm,0.5cm/1.5cm,1.5cm/2.5cm}
\draw
[fill=gray!20,xshift=\x,yshift=\y](0.5,0)..controls(0.1,0.1)..(0,0.5)..controls(0.4,0.4)..(0.5,0);

\end{tikzpicture}
\begin{tikzpicture}[rotate=90,scale=1.5]
\node at (2.25,0) {$~~~=$};
\node at (1.25,1.25) {$y$};
\node at (2.22,0.28) {$w$};
\node at (2.22,1.25) {$x$};
\node at (1.25,0.3) {$z$};
\path[clip] (1.25,0.25) rectangle (3.25,2.25);
\foreach \x/\y in {0cm/0cm,1cm/1cm,2cm/2cm,-1cm/1cm,1cm/-1cm,0cm/2cm,2cm/0cm,3cm/1cm}
\draw[fill=gray!20,xshift=\x,yshift=\y](0.5,0)--(1.5,1)..controls(1.1,1.1)..(1,1.5)--(0,0.5)..controls(0.4,0.4)..(0.5,0);
\foreach \x/\y in {-0.5cm/0.5cm,0.5cm/-0.5cm,1.5cm/0.5cm,2.5cm/1.5cm,0.5cm/1.5cm,1.5cm/2.5cm,2.5cm/-0.5cm}
\draw
[fill=gray!20,xshift=\x,yshift=\y](0.5,0)..controls(0.1,0.1)..(0,0.5)..controls(0.4,0.4)..(0.5,0);
\end{tikzpicture}
\begin{tikzpicture}[rotate=90,scale=1.5]
\node at (1.25,1.23) {$y$};
\node at (2.22,2.25) {$w$};
\node at (2.22,1.25) {$x$};
\node at (1.25,2.25) {$z$};
\node at (0.71,2.25) {$~$};
\node at (1.75,0.5) {$~~~\rightarrow$};
\path[clip] (0.75,0.75) rectangle (2.75,2.75);
\foreach \x/\y in {0cm/0cm,1cm/1cm,2cm/2cm,-1cm/1cm,1cm/-1cm,0cm/2cm,2cm/0cm,4cm/2cm}
\draw[fill=gray!20,xshift=\x,yshift=\y](0.5,0)--(1.5,1)..controls(1.1,1.1)..(1,1.5)--(0,0.5)..controls(0.4,0.4)..(0.5,0);
\foreach \x/\y in {-0.5cm/0.5cm,0.5cm/-0.5cm,1.5cm/0.5cm,2.5cm/1.5cm,0.5cm/1.5cm,1.5cm/2.5cm,2.5cm/-0.5cm}
\draw
[fill=gray!20,xshift=\x,yshift=\y](0.5,0)..controls(0.1,0.1)..(0,0.5)..controls(0.4,0.4)..(0.5,0);
\end{tikzpicture}
\begin{tikzpicture}[scale=1.5]
\path[clip] (0.75,-0.25) rectangle (2.75,1.75);
\foreach \x/\y in {0cm/0cm,1cm/1cm,2cm/2cm,-1cm/1cm,1cm/-1cm,0cm/2cm,2cm/0cm,4cm/2cm}
\draw[fill=gray!20,xshift=\x,yshift=\y](0.5,0)--(1.5,1)..controls(1.1,1.1)..(1,1.5)--(0,0.5)..controls(0.4,0.4)..(0.5,0);
\foreach \x/\y in {-0.5cm/0.5cm,0.5cm/-0.5cm,1.5cm/0.5cm,2.5cm/1.5cm,0.5cm/1.5cm,1.5cm/2.5cm,2.5cm/-0.5cm}
\draw
[fill=gray!20,xshift=\x,yshift=\y](0.5,0)..controls(0.1,0.1)..(0,0.5)..controls(0.4,0.4)..(0.5,0);
\node at (1.29,1.29) {\small $w^{-1}$};
\node at (2.27,0.3) {\small $y^{-1}$};
\node at (2.26,1.25) {$x\frac{y^2(1+w)^2}{(1+y)^2}$};
\node at (1.26,0.25) {$z\frac{w^2(1+y)^2}{(1+w)^2}$};
\node at (1.22,-0.3) {$~$};
\end{tikzpicture}
\end{center}
\caption{Discrete flow. The first move is a composition of two grey mutations, the second is just the fundamental domain change and the third is a composition of two white mutations. }\label{discrete-flow}
\end{figure}

The group $\mathcal{G}'_\Delta$ is a group of rank one with torsion $\mathbb{Z}/2\mathbb{Z}$. The torsion generator $(1,1,0,0)$ is given by the permutation of coordinates $(x,y,z,w)\mapsto (z,w,x,y)$. The generator of infinite order $(0,1,0,0)$ is given by
 $$(x,y,z,w)\mapsto\left( w^{-1}, zw^2\frac{(1+y)^2}{(1+w)^2}, y^{-1}, xy^2\frac{(1+w)^2}{(1+y)^2}\right).$$
and corresponds to two grey and two white mutations as shown on fig.~\ref{discrete-flow}.

\subsection{Degeneration to non-affine Toda system}

This example illustrates that integrable systems, corresponding to a smaller Newton polygon, can be obtained from a larger one by a certain degeneration procedure. Another reason to consider this example is that it formally does not fit into the GK scheme \cite{GK}, since the corresponding bipartite graph is not minimal in their sense.

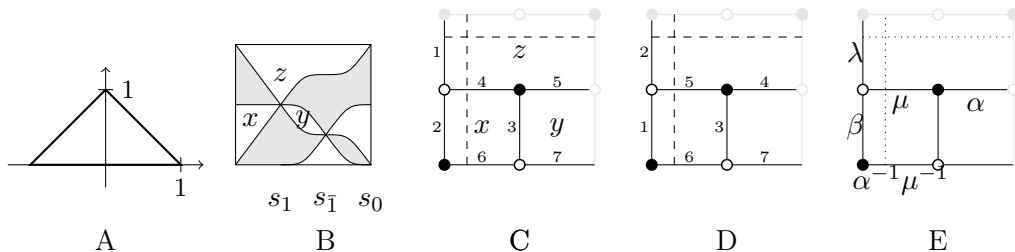
\begin{figure}[H]
\begin{center}
\begin{tikzpicture}
\draw[->](-1.3,0)--(1.3,0);
\draw[->](0,-0.3)--(0,1.3);
\draw (1,0.1)--(1,-0.1);\node at (1,-0.3) {\small 1};
\draw (0.1,1)--(-0.1,1);\node at (0.3,1) {\small 1};
\draw[thick] (-1,0)--(0,1)--(1,0)--cycle;
\node at (0,-1) {A};
\end{tikzpicture}\quad
\begin{tikzpicture}
\fill[color=gray!20] (0,1.6)--(0.6,0.8)--(0,0.8)--cycle;
\draw                (0,1.6)--(0.6,0.8)--(0,0.8);
\fill[color=gray!20] (1.2,0.4)..controls(1.5,0.4)..(1.8,0)..controls(1.5,0)..(1.2,0.4)--cycle;
\draw                (1.2,0.4)..controls(1.5,0.4)..(1.8,0)..controls(1.5,0)..(1.2,0.4)--cycle;
\fill[color=gray!20] (1.8,1.6)..controls(1.5,1.2)..(1.2,1.2)..controls(0.9,1.2)..(0.6,0.8)..controls(0.9,0.8)..(1.2,0.4)..controls(1.5,0.8)..(1.8,0.8)--cycle;
\draw (1.8,1.6)..controls(1.5,1.2)..(1.2,1.2)..controls(0.9,1.2)..(0.6,0.8)..controls(0.9,0.8)..(1.2,0.4)..controls(1.5,0.8)..(1.8,0.8);
\fill[color=gray!20] (0,0)--(0.6,0)..controls (0.9,0)..(1.2,0.4)..controls(0.9,0.4)..(0.6,0.8)--cycle;
\draw                (0.6,0)..controls (0.9,0)..(1.2,0.4)..controls(0.9,0.4)..(0.6,0.8)--(0,0);
\draw[thin] (0,0)--(1.8,0.0)--(1.8,1.6)--(0,1.6)--cycle;
\node at (0.2,0.6) {$x$};
\node at (0.9,0.6) {$y$};
\node at (0.6,1.2) {$z$};
\node at (0.6,-0.5) {$s_1$};
\node at (1.2,-0.5) {$s_{\bar{1}}$};
\node at (1.8,-0.5) {$s_0$};
\node at (1.2,-1) {B};
\end{tikzpicture}\quad
\begin{tikzpicture}
\draw (0,0) -- +(0,2);
\draw (0,0) -- +(2,0);
\draw (0,0)+(0,1) -- +(2,1);
\draw (0,0)+(1,0) -- +(1,1);
\draw[color=gray!20] (0,2)--(2,2)--(2,0);
\fill (0,0)+(0,0) circle (0.08);
\fill[color=gray!20] (0,2) circle (0.08); \fill[color=gray!20] (1,2) circle (0.08); \fill[color=gray!20] (2,2) circle (0.08);
\fill[color=gray!20] (2,2) circle (0.08); \fill[color=gray!20] (2,1) circle (0.08);
\fill[color=white] (1,2) circle (0.06);\fill[color=white] (2,1) circle (0.06);
\fill (0,0)+(0,1) circle (0.08);\fill[color=white] (0,1) circle (0.06);
\fill (0,0)+(1,0) circle (0.08);\fill[color=white] (1,0) circle (0.06);
\fill (0,0)+(1,1) circle (0.08);
\node at (1,-1) {C};
\node at (-0.1,1.5) {\tiny 1};
\node at (-0.1,0.5) {\tiny 2};
\node at (0.9,0.5) {\tiny 3};
\node at (0.5,1.1) {\tiny 4};
\node at (1.5,1.1) {\tiny 5};
\node at (0.5,0.1) {\tiny 6};
\node at (1.5,0.1) {\tiny 7};
\draw[dashed] (0.3,0)--(0.3,2);\draw[dashed] (0,1.7)--(2,1.7);
\node at (1,1.5) {$z$};
\node at (1.5,0.5) {$y$};
\node at (0.5,0.5) {$x$};
\node at (1,-1) {C};
\end{tikzpicture}\quad
\begin{tikzpicture}
\draw (0,0) -- +(0,2);
\draw (0,0) -- +(2,0);
\draw (0,0)+(0,1) -- +(2,1);
\draw (0,0)+(1,0) -- +(1,1);
\draw[color=gray!20] (0,2)--(2,2)--(2,0);
\fill (0,0)+(0,0) circle (0.08);
\fill[color=gray!20] (0,2) circle (0.08); \fill[color=gray!20] (1,2) circle (0.08); \fill[color=gray!20] (2,2) circle (0.08);
\fill[color=gray!20] (2,2) circle (0.08); \fill[color=gray!20] (2,1) circle (0.08);
\fill[color=white] (1,2) circle (0.06);\fill[color=white] (2,1) circle (0.06);
\fill (0,0)+(0,1) circle (0.08);\fill[color=white] (0,1) circle (0.06);
\fill (0,0)+(1,0) circle (0.08);\fill[color=white] (1,0) circle (0.06);
\fill (0,0)+(1,1) circle (0.08);
\node at (-0.1,1.5) {\tiny 2};
\node at (-0.1,0.5) {\tiny 1};
\node at (0.9,0.5) {\tiny 3};
\node at (0.5,1.1) {\tiny 5};
\node at (1.5,1.1) {\tiny 4};
\node at (0.5,0.1) {\tiny 6};
\node at (1.5,0.1) {\tiny 7};
\draw[dashed] (0.3,0)--(0.3,2);\draw[dashed] (0,1.7)--(2,1.7);
\node at (1,-1) {D};
\end{tikzpicture}\quad
\begin{tikzpicture}
\draw (0,0) -- +(0,2);
\draw (0,0) -- +(2,0);
\draw (0,0)+(0,1) -- +(2,1);
\draw (0,0)+(1,0) -- +(1,1);
\draw[color=gray!20] (0,2)--(2,2)--(2,0);
\fill (0,0)+(0,0) circle (0.08);
\fill[color=gray!20] (0,2) circle (0.08); \fill[color=gray!20] (1,2) circle (0.08); \fill[color=gray!20] (2,2) circle (0.08);
\fill[color=gray!20] (2,2) circle (0.08); \fill[color=gray!20] (2,1) circle (0.08);
\fill[color=white] (1,2) circle (0.06);\fill[color=white] (2,1) circle (0.06);
\fill (0,0)+(0,1) circle (0.08);\fill[color=white] (0,1) circle (0.06);
\fill (0,0)+(1,0) circle (0.08);\fill[color=white] (1,0) circle (0.06);
\fill (0,0)+(1,1) circle (0.08);
\node at (-0.1,1.5) {$\lambda$};
\node at (-0.1,0.5) {$\beta$};
\node at (0.5,0.8) {$\mu$};
\node at (1.5,0.8) {$\alpha$};
\node at (0.5,-0.2) {$\alpha^{-1}\mu^{-1}$};
\draw[dotted] (0.3,0)--(0.3,2);\draw[dotted] (0,1.7)--(2,1.7);
\node at (1,-1) {E};
\end{tikzpicture}
\end{center}
\caption{Non-affine relativistic two-particle Toda chain. (A): Newton polygon; (B): Thurston diagram; (C): Bipartite graph on a torus $\Sigma$ and homology basis; (D): Bipartite graph on dual torus $\tilde{\Sigma}$ and a homology basis; (E): Edge variables.}
\label{fi:todasimple}
\end{figure}

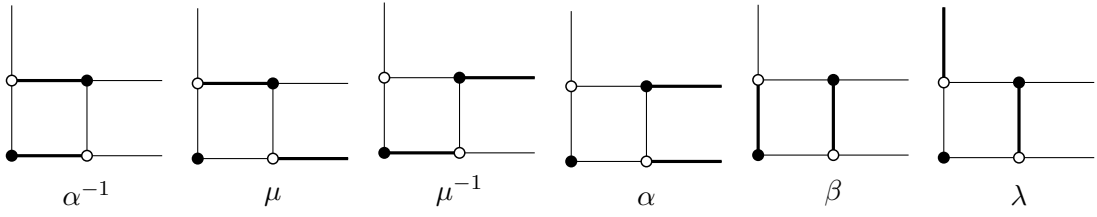
\begin{figure}[H]
\begin{center}
\begin{tikzpicture}
\draw[very thick] (0,0) -- (1,0);
\draw[very thick] (0,1) -- (1,1);
\draw (0,0) -- (0,2);
\draw (0,0) -- +(2,0);
\draw (0,0)+(0,1) -- +(2,1);
\draw (0,0)+(1,0) -- +(1,1);
\fill (0,0)+(0,0) circle (0.08);
\fill (0,0)+(0,1) circle (0.08);\fill[color=white] (0,1) circle (0.06);
\fill (0,0)+(1,0) circle (0.08);\fill[color=white] (1,0) circle (0.06);
\fill (0,0)+(1,1) circle (0.08);
\node at (1,-0.5) {$\alpha^{-1}$};
\end{tikzpicture}\quad
\begin{tikzpicture}
\draw[very thick] (1,0) -- (2,0);
\draw[very thick] (0,1) -- (1,1);
\draw (0,0) -- +(0,2);
\draw (0,0) -- +(2,0);
\draw (0,0)+(0,1) -- +(2,1);
\draw (0,0)+(1,0) -- +(1,1);
\fill (0,0)+(0,0) circle (0.08);
\fill (0,0)+(0,1) circle (0.08);\fill[color=white] (0,1) circle (0.06);
\fill (0,0)+(1,0) circle (0.08);\fill[color=white] (1,0) circle (0.06);
\fill (0,0)+(1,1) circle (0.08);
\node at (1,-0.5) {$\mu$};
\end{tikzpicture}\quad
\begin{tikzpicture}
\draw[very thick] (0,0) -- (1,0);
\draw[very thick] (1,1) -- (2,1);
\draw (0,0) -- +(0,2);
\draw (0,0) -- +(2,0);
\draw (0,0)+(0,1) -- +(2,1);
\draw (0,0)+(1,0) -- +(1,1);
\fill (0,0)+(0,0) circle (0.08);
\fill (0,0)+(0,1) circle (0.08);\fill[color=white] (0,1) circle (0.06);
\fill (0,0)+(1,0) circle (0.08);\fill[color=white] (1,0) circle (0.06);
\fill (0,0)+(1,1) circle (0.08);
\node at (1,-0.5) {$\mu^{-1}$};
\end{tikzpicture}\quad
\begin{tikzpicture}
\draw[very thick] (1,0) -- (2,0);
\draw[very thick] (1,1) -- (2,1);
\draw (0,0) -- +(0,2);
\draw (0,0) -- +(2,0);
\draw (0,0)+(0,1) -- +(2,1);
\draw (0,0)+(1,0) -- +(1,1);
\fill (0,0)+(0,0) circle (0.08);
\fill (0,0)+(0,1) circle (0.08);\fill[color=white] (0,1) circle (0.06);
\fill (0,0)+(1,0) circle (0.08);\fill[color=white] (1,0) circle (0.06);
\fill (0,0)+(1,1) circle (0.08);
\node at (1,-0.5) {$\alpha$};
\end{tikzpicture}\quad
\begin{tikzpicture}
\draw[very thick] (0,0) -- (0,1);
\draw[very thick] (1,0) -- (1,1);
\draw (0,0) -- +(0,2);
\draw (0,0) -- +(2,0);
\draw (0,0)+(0,1) -- +(2,1);
\draw (0,0)+(1,0) -- +(1,1);
\fill (0,0)+(0,0) circle (0.08);
\fill (0,0)+(0,1) circle (0.08);\fill[color=white] (0,1) circle (0.06);
\fill (0,0)+(1,0) circle (0.08);\fill[color=white] (1,0) circle (0.06);
\fill (0,0)+(1,1) circle (0.08);
\node at (1,-0.5) {$\beta$};
\end{tikzpicture}\quad
\begin{tikzpicture}
\draw[very thick] (0,1) -- (0,2);
\draw[very thick] (1,0) -- (1,1);
\draw (0,0) -- +(0,2);
\draw (0,0) -- +(2,0);
\draw (0,0)+(0,1) -- +(2,1);
\draw (0,0)+(1,0) -- +(1,1);
\fill (0,0)+(0,0) circle (0.08);
\fill (0,0)+(0,1) circle (0.08);\fill[color=white] (0,1) circle (0.06);
\fill (0,0)+(1,0) circle (0.08);\fill[color=white] (1,0) circle (0.06);
\fill (0,0)+(1,1) circle (0.08);
\node at (1,-0.5) {$\lambda$};
\end{tikzpicture}
\end{center}
\caption{Dimer configuration for the simplest non-affine Toda system.\label{fi:todasimpledimers}}
\end{figure}

Consider the Newton polygon (triangle) with vertices at $(-1,0),(1,0)$ and $(0,1)$ of area $\mathsf{S}=1$, shown on fig.~\ref{fi:todasimple}A.
There are six dimer configurations shown on fig.~\ref{fi:todasimpledimers},
\noindent
and the corresponding double partition function is
$$
S^d=\alpha+\alpha^{-1}+\beta +\mu+\mu^{-1}+\lambda
$$
The relation between cluster and face variables is
$$
x=\beta\alpha,\quad y=\alpha\beta^{-1},\quad z=\alpha^{-2}
$$
Thus, the normalised face partition function is
\begin{equation}\label{Sftodasimple}
S^f=y^{1/2}x^{1/2}+y^{-1/2}x^{-1/2}+y^{-1/2}x^{1/2} +\mu+\mu^{-1}+\lambda
\end{equation}
and it should coincide with given by the characteristic polynomial
\be
\det(H_1(x)E_1H_1(y)E_{\bar{1}}H_0(z)E_0(\lambda)-\mu)=\mu^2-\mu(1+x+yx)-\lambda\mu x + yx
\ee
Indeed, normalising the r.h.s. by the substitution $\mu\to -y^{1/2}x^{1/2}\mu$, $\lambda\to y^{1/2}x^{-1/2}\lambda$ and dividing it by $\mu yx$ one recovers the expression \rf{Sftodasimple}.

One finds, that in this case the spectral curve is rational, and it can be obtained as degeneration
of the elliptic curve \rf{Toda1} of the simplest affine relativistic Toda in spirit of \cite{KriVa,BraMa},
indeed we just get \rf{Sftodasimple} in the limit $c\to 0$ of \rf{Toda1}.

\subsection{Relativistic Toda chain of rank two
\label{ss:toda2}}

In both previous examples the integrable system of sect.~\ref{ss:dimerstoloop} was trivial, since the phase space was two-dimensional, and it has single integral of motion - the Hamiltonian. Let us now turn to
the simplest example with more than one commuting integrable flows - relativistic Toda chain with four-dimensional phase space and two independent integrals of motion.

\begin{figure}[H]
\begin{center}
\begin{tikzpicture}
\draw[->](-1.3,0)--(1.3,0);
\draw[->](0,-2.3)--(0,1.3);
\node at (1.2,-0.3) {\small 1};
\node at (0.3,1.2) {\small 1};
\draw[thick] (-1,0)--(0,1)--(1,-1)--(0,-2)--cycle;
\node at (0,-2.5) {A};
\end{tikzpicture}\quad
\begin{tikzpicture}
\fill[color=gray!20] (3,1.6)..controls(2.7,1.6)..(2.4,1.2)..controls(2.1,1.6)..(1.8,1.6)--cycle;
\draw                (3,1.6)..controls(2.7,1.6)..(2.4,1.2)..controls(2.1,1.6)..(1.8,1.6);
\fill[color=gray!20] (2.4,1.2)..controls(2.7,1.2)..(3,0.8)..controls(2.7,0.8)..(2.4,1.2)--cycle;
\draw                (2.4,1.2)..controls(2.7,1.2)..(3,0.8)..controls(2.7,0.8)..(2.4,1.2)--cycle;
\fill[color=gray!20] (0,1.2)..controls(0.3,1.6)..(0.6,1.6)--(0,1.6)--cycle;
\draw                (0,1.2)..controls(0.3,1.6)..(0.6,1.6);
\fill[color=gray!20] (0,1.2)..controls(0.3,1.2)..(0.6,0.8)..controls(0.3,0.8)..(0,1.2)--cycle;
\draw                (0,1.2)..controls(0.3,1.2)..(0.6,0.8)..controls(0.3,0.8)..(0,1.2);
\fill[color=gray!20] (1.2,0.4)..controls(1.5,0.4)..(1.8,0)..controls(1.5,0)..(1.2,0.4)--cycle;
\draw                (1.2,0.4)..controls(1.5,0.4)..(1.8,0)..controls(1.5,0)..(1.2,0.4);
\fill[color=gray!20] (1.8,0)..controls(2.1,0.4)..(2.4,0.4)..controls(2.7,0.4)..(3,0.8)..controls(3.3,0.4)..(3.6,0.4)..controls(3.3,0)..(3,0)--cycle;
\draw                (1.8,0)..controls(2.1,0.4)..(2.4,0.4)..controls(2.7,0.4)..(3,0.8)..controls(3.3,0.4)..(3.6,0.4)..controls(3.3,0)..(3,0);
\fill[color=gray!20] (3.6,1.2)..controls(3.3,1.2)..(3,0.8)..controls(3.3,0.8)..(3.6,0.4)--cycle;
\draw                (3.6,1.2)..controls(3.3,1.2)..(3,0.8)..controls(3.3,0.8)..(3.6,0.4);
\fill[color=gray!20] (2.4,1.2)..controls(2.1,1.2)..(1.8,1.6)..controls(1.5,1.2)..(1.2,1.2)..controls(0.9,1.2)..(0.6,0.8)..controls(0.9,0.8)..(1.2,0.4)..controls(1.5,0.8)..(1.8,0.8)..controls(2.1,0.8)..(2.4,1.2)--cycle;
\draw (2.4,1.2)..controls(2.1,1.2)..(1.8,1.6)..controls(1.5,1.2)..(1.2,1.2)..controls(0.9,1.2)..(0.6,0.8)..controls(0.9,0.8)..(1.2,0.4)..controls(1.5,0.8)..(1.8,0.8)..controls(2.1,0.8)..(2.4,1.2);
\fill[color=gray!20] (0,0)--(0.6,0)..controls (0.9,0)..(1.2,0.4)..controls(0.9,0.4)..(0.6,0.8)..controls(0.3,0.4)..(0,0.4)--cycle;
\draw                (0,0)--(0.6,0)..controls (0.9,0)..(1.2,0.4)..controls(0.9,0.4)..(0.6,0.8)..controls(0.3,0.4)..(0,0.4);
\draw[thin] (0,0)--(3.6,0)--(3.6,1.6)--(0,1.6)--cycle;
\node at (0.1,0.7) {$x$};
\node at (1,1.4) {$z$};
\node at (0.9,0.6) {$y$};
\node at (2.1,1.4) {$w$};
\node at (2.3,0.6) {$u$};
\node at (3.3,0.6) {$v$};
\node at (0.6,-0.5) {$s_1$};
\node at (1.2,-0.5) {$s_{\bar{1}}$};
\node at (1.8,-0.5) {$s_0$};
\node at (2.4,-0.5) {$s_{\bar{0}}$};
\node at (3,-0.5) {$s_1$};
\node at (3.6,-0.5) {$s_{\bar{1}}$};
\node at (1.5,-1.2) {B};
\end{tikzpicture}\quad
\begin{tikzpicture}
\draw[->](-1.3,0)--(2.3,0);
\draw[->](0,-1.3)--(0,1.3);
\node at (1.2,-0.3) {\small 1};
\node at (0.3,1.2) {\small 1};
\draw[thick] (-1,0)--(0,1)--(2,0)--(1,-1)--cycle;
\node at (0,-2) {A'};
\end{tikzpicture}\quad
\begin{tikzpicture}
\fill[color=gray!20] (0,0)--(0.6,0)..controls(0.3,0)..(0,0.4)--cycle;
\draw                (0.6,0)..controls(0.3,0)..(0,0.4);
\fill[color=gray!20] (0,0.4)..controls(0.3,0.8)..(0.6,0.8)--(1.8,0.8)..controls(2.1,0.8)..(2.4,1.2)..controls(2.7,0.8)..(3,0.8)..controls(2.7,0.4)..(2.4,0.4)--(1.2,0.4)..controls(0.9,0.4)..(0.6,0)..controls(0.3,0.4)..(0,0.4)--cycle;
\draw (0,0.4)..controls(0.3,0.8)..(0.6,0.8)--(1.8,0.8)..controls(2.1,0.8)..(2.4,1.2)..controls(2.7,0.8)..(3,0.8)..controls(2.7,0.4)..(2.4,0.4)--(1.2,0.4)..controls(0.9,0.4)..(0.6,0)..controls(0.3,0.4)..(0,0.4);
\fill[color=gray!20] (0,1.2)--(1.2,1.2)..controls(1.5,1.2)..(1.8,1.6)..controls(1.5,1.6)..(1.2,2)..controls(0.9,1.6)..(0.6,1.6)--(0,1.6)--cycle;
\draw (0,1.2)--(1.2,1.2)..controls(1.5,1.2)..(1.8,1.6)..controls(1.5,1.6)..(1.2,2)..controls(0.9,1.6)..(0.6,1.6)--(0,1.6);
\fill[color=gray!20] (0,2)..controls(0.3,2)..(0.6,2.4)--(0,2.4)--cycle;
\draw (0,2)..controls(0.3,2)..(0.6,2.4);
\fill[color=gray!20] (3.6,0.4)..controls(3.3,0.8)..(3,0.8)..controls(3.3,0.4)..(3.6,0.4)--cycle;
\draw                (3.6,0.4)..controls(3.3,0.8)..(3,0.8)..controls(3.3,0.4)..(3.6,0.4);
\fill[color=gray!20] (2.4,1.2)..controls(2.1,1.6)..(1.8,1.6)..controls(2.1,1.2)..(2.4,1.2)--cycle;
\draw                (2.4,1.2)..controls(2.1,1.6)..(1.8,1.6)..controls(2.1,1.2)..(2.4,1.2);
\fill[color=gray!20] (1.2,2)..controls(0.9,2.4)..(0.6,2.4)..controls(0.9,2)..(1.2,2)--cycle;
\draw                (1.2,2)..controls(0.9,2.4)..(0.6,2.4)..controls(0.9,2)..(1.2,2);
\fill[color=gray!20] (3,0)..controls(3.3,0)..(3.6,0.4)--(3.6,0)--cycle;
\draw                (3,0)..controls(3.3,0)..(3.6,0.4)--(3.6,0);
\fill[color=gray!20] (3.6,1.2)..controls(3.3,1.2)..(3,0.8)..controls(2.7,1.2)..(2.4,1.2)..controls(2.7,1.6)..(3,1.6)--(3.6,1.6)--cycle;
\draw (3.6,1.2)..controls(3.3,1.2)..(3,0.8)..controls(2.7,1.2)..(2.4,1.2)..controls(2.7,1.6)..(3,1.6)--(3.6,1.6);
\fill[color=gray!20] (3.6,2)--(2.4,2)..controls(2.1,2)..(1.8,1.6)..controls(1.5,2)..(1.2,2)..controls(1.5,2.4)..(1.8,2.4)--(3.6,2.4)--cycle;
\draw (3.6,2)--(2.4,2)..controls(2.1,2)..(1.8,1.6)..controls(1.5,2)..(1.2,2)..controls(1.5,2.4)..(1.8,2.4)--(3.6,2.4);
\node at (0.6,-0.5) {$s_0$};
\node at (1.2,-0.5) {$s_{\bar{0}}$};
\node at (1.8,-0.5) {$s_1$};
\node at (2.4,-0.5) {$s_{\bar{1}}$};
\node at (3,-0.5) {$s_2$};
\node at (3.6,-0.5) {$s_{\bar{2}}$};
\node at (0.3,0.2) {$w$};
\node at (1.5,1.8) {$y$};
\node at (2.7,1) {$v$};
\node at (1.8,1.1) {$z$};
\node at (0.6,1.9) {$u$};
\node at (3,0.3) {$x$};
\node at (3.3,-1.2) {B'};
\end{tikzpicture}

\begin{tikzpicture}
\draw (0,0) -- +(0,2);
\draw (1,0) -- +(0,2);
\draw (2,0) -- +(0,2);
\draw (0,0) -- +(3,0);
\draw (0,1) -- +(3,0);
\draw[color=gray!20] (0,2) -- (3,2)--(3,0);
\fill (0,0)+(0,0) circle (0.08);
\fill (0,0)+(1,1) circle (0.08);
\fill (0,0)+(2,0) circle (0.08);
\fill[color=gray!20] (0,0)+(0,2) circle (0.08);
\fill[color=gray!20] (0,0)+(2,2) circle (0.08);
\fill[color=gray!20] (0,0)+(3,1) circle (0.08);
\fill (0,0)+(0,1) circle (0.08);\fill[color=white] (0,1) circle (0.06);
\fill (0,0)+(1,0) circle (0.08);\fill[color=white] (1,0) circle (0.06);
\fill (0,0)+(2,1) circle (0.08);\fill[color=white] (2,1) circle (0.06);
\fill[color=gray!20] (0,0)+(1,2) circle (0.08);\fill[color=white] (1,2) circle (0.06);
\fill[color=gray!20] (0,0)+(3,2) circle (0.08);\fill[color=white] (3,2) circle (0.06);
\fill[color=gray!20] (0,0)+(3,0) circle (0.08);\fill[color=white] (3,0) circle (0.06);
\node at (-0.1,1.5) {\tiny 1};\node at (0.9,1.5) {\tiny 2};
\node at (1.9,1.5) {\tiny 3};\node at (2.9,1.5) {\tiny 4};
\node at (-0.1,0.5) {\tiny 4};\node at (0.9,0.5) {\tiny 5};
\node at (1.9,0.5) {\tiny 6};\node at (2.9,0.5) {\tiny 1};
\node at (0.5,0.1) {\tiny 7};\node at (1.5,0.1) {\tiny 8};\node at (2.5,0.1) {\tiny 9};
\node at (0.5,1.1) {\tiny 10};\node at (1.5,1.1) {\tiny 11};\node at (2.5,1.1) {\tiny 12};
\node at (0.5,2.1) {\tiny 7};\node at (1.5,2.1) {\tiny 8};\node at (2.5,2.1) {\tiny 9};
\node at (0.5,0.5) {$z$};
\node at (1.5,0.5) {$v$};
\node at (2.5,0.5) {$u$};
\node at (0.5,1.5) {$y$};
\node at (1.5,1.5) {$x$};
\node at (2.5,1.5) {$w$};
\draw[dashed] (0.3,0)--(0.3,2);
\draw[dashed] (0,1.3)--(3,0.7);
\node at (1,-1) {C};
\end{tikzpicture}\quad
\begin{tikzpicture}
\draw (0.5,1)..controls (1.2,0.8)..(1.5,0.3);\node at (1.25,0.85) {\tiny 4};
\draw[color=gray!40] (1.5,0.3)..controls (0.8,0.6)..(0.5,1);\node at (0.75,0.55) {\tiny 1};
\draw (-0.5,1)..controls (-1.2,0.8)..(-1.5,0.3);\node at (-1.25,0.85) {\tiny 6};
\draw[color=gray!40] (-1.5,0.3)..controls (-0.8,0.6)..(-0.5,1);\node at (-0.75,0.55) {\tiny 3};
\draw (0.5,-1)..controls (1.2,-0.8)..(1.5,-0.3);\node at (1.25,-0.85) {\tiny 6};
\draw[color=gray!40] (1.5,-0.3)..controls (0.8,-0.6)..(0.5,-1);\node at (0.75,-0.55) {\tiny 3};
\draw (-0.5,-1)..controls (-1.2,-0.8)..(-1.5,-0.3);\node at (-1.25,-0.85) {\tiny 4};
\draw[color=gray!40] (-1.5,-0.3)..controls (-0.8,-0.6)..(-0.5,-1);\node at (-0.75,-0.55) {\tiny 1};
\draw (1.5,0.3)..controls (1.3,0)..(1.5,-0.3);\node at (1.2,0) {\tiny 12};
\draw (-1.5,-0.3)..controls (-1.3,0)..(-1.5,0.3);\node at (-1.2,0) {\tiny 9};
\draw (0.5,1)..controls (0,0.4)..(-0.5,1);\node at (0.4,0.7) {\tiny 10};\node at (-0.4,0.7) {\tiny 11};
\draw (-0.5,-1)..controls (0,-0.4)..(0.5,-1);\node at (0.4,-0.7) {\tiny 8};\node at (-0.4,-0.7) {\tiny 7};
\draw[color=gray!40] (0,-0.55)..controls (0.2,0)..(0,0.55);\node at (0.3,0) {\tiny 5};
\draw (0,-0.55)..controls (-0.2,0)..(0,0.55);\node at (-0.3,0) {\tiny 2};
\fill (0.5,1) circle (0.08);\fill[color=white] (0.5,1) circle (0.06);
\fill (-0.5,1) circle (0.08);\fill[color=white] (-0.5,1) circle (0.06);
\fill (0,-0.55) circle (0.08);\fill[color=white] (0,-0.55) circle (0.06);
\fill (-1.5,-0.3) circle (0.08);\fill[color=white] (-1.5,-0.3) circle (0.06);
\fill (1.5,-0.3) circle (0.08);\fill[color=white] (1.5,-0.3) circle (0.06);
\fill (0.5,-1) circle (0.08);
\fill (-0.5,-1) circle (0.08);
\fill (0,0.55) circle (0.08);
\fill (-1.5,0.3) circle (0.08);
\fill (1.5,0.3) circle (0.08);
\draw[color=gray!40, dashed] (-1,0.55)..controls(-0.5,0) and (0.5,0)..(1,-0.55);
\draw[dashed] (1.13,0.75)..controls(0.5,0) and (-0.5,0)..(-1.13,-0.75);
\draw[dashed] (0.22,0.7)..controls(0.9,0.4)..(1.33,0.0);
\draw[color=gray!40, dashed] (0.22,0.7)..controls(0.7,0.2)..(1.33,0.0);
\draw[dashed] (-0.22,0.7)..controls(-0.9,0.4)..(-1.33,0.0);
\draw[color=gray!40, dashed] (-0.22,0.7)..controls(-0.7,0.2)..(-1.33,0.0);
\node at (0,-2) {D};
\end{tikzpicture}\quad
\begin{tikzpicture}
\draw (0,0) -- +(0,2);
\draw (1,0) -- +(0,2);
\draw (2,0) -- +(0,2);
\draw (0,0) -- +(3,0);
\draw (0,1) -- +(3,0);
\draw[color=gray!20] (0,2) -- (3,2)--(3,0);
\fill (0,0)+(0,0) circle (0.08);
\fill (0,0)+(1,1) circle (0.08);
\fill (0,0)+(2,0) circle (0.08);
\fill[color=gray!20] (0,0)+(0,2) circle (0.08);
\fill[color=gray!20] (0,0)+(2,2) circle (0.08);
\fill[color=gray!20] (0,0)+(3,1) circle (0.08);
\fill (0,0)+(0,1) circle (0.08);\fill[color=white] (0,1) circle (0.06);
\fill (0,0)+(1,0) circle (0.08);\fill[color=white] (1,0) circle (0.06);
\fill (0,0)+(2,1) circle (0.08);\fill[color=white] (2,1) circle (0.06);
\fill[color=gray!20] (0,0)+(1,2) circle (0.08);\fill[color=white] (1,2) circle (0.06);
\fill[color=gray!20] (0,0)+(3,2) circle (0.08);\fill[color=white] (3,2) circle (0.06);
\fill[color=gray!20] (0,0)+(3,0) circle (0.08);\fill[color=white] (3,0) circle (0.06);
\draw[dotted] (0.3,0)--(0.3,2);
\draw[dotted] (0,1.3)--(3,0.7);
\node at (1,1.5) {$\alpha^{-1}\mu^{-1}$};
\node at (2.1,1.5) {$\gamma$};
\node at (-0.1,0.5) {$\alpha$};\node at (1.1,0.5) {$\gamma^{-1}$};
\node at (-0.1,1.5) {$\mu$};
\node at (2.1,0.5) {$\mu^{-1}$};
\node at (0.5,-0.2) {$\lambda$};\node at (2.5,-0.2) {$\delta^{-1}$};\node at (1.5,-0.2) {$c$};
\node at (0.5,1) {$\beta\lambda^{-1}$};\node at (1.5,1) {$\delta\mu^{-1}$};\node at (2.5,1) {$\beta^{-1}$};
\node at (1,-1) {E};
\end{tikzpicture}
\end{center}
\caption{Affine $\widehat{PGL(3)}$ relativistic Toda chain. (A): Newton polygon, (B): Thurston diagram (right and left sides are glued after a shift by a half-turn), and the corresponding word in the double Weyl group of $\widehat{PGL(2)}$; (A'): Rotated Newton polygon; (B'): Rotated Thurston diagram and the corresponding word in the double Weyl group of $\widehat{PGL(3)}$; (C): Bipartite graph on a torus $\Sigma$, two homology cycles and the face variables; (D): Bipartite graph on dual surface $\widetilde{\Sigma}$ of genus two and the homology basis; (E): Edge variables $A_e=A_e(\lambda,\mu;\alpha,\beta,\gamma,\delta;c)$.\label{fi:toda23}}
\end{figure}
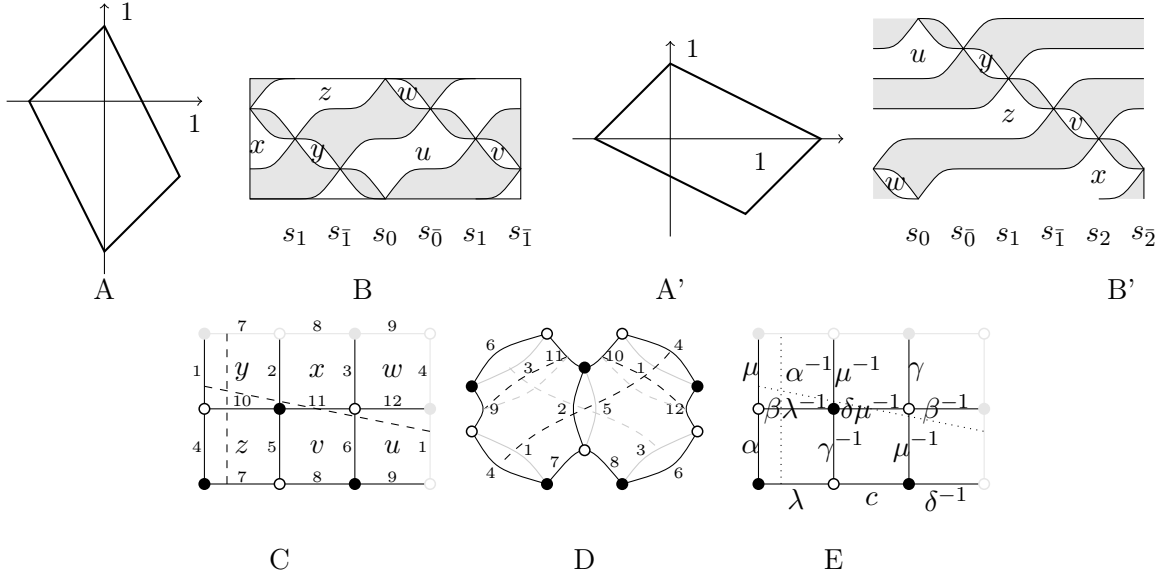

Consider Newton polygon - the parallelogram with the vertices $(-1,0)$, $(0,1)$, $(2,0)$ and $(1,-1)$, shown on fig.~\ref{fi:toda23}A. The corresponding Thurston diagram on a torus $\Sigma$ is shown on fig.~\ref{fi:toda23}B. After applying an automorphism of the lattice, the same Newton polygon can be presented as parallelogram from fig.~\ref{fi:toda23}A'. The corresponding to fig.~\ref{fi:toda23}A' Thurston diagram, depicted at fig.~\ref{fi:toda23}B', differs from that of fig.~\ref{fi:toda23}B just by the torus automorphism, however, after cutting tori on fig.~\ref{fi:toda23}B and fig.~\ref{fi:toda23}B' in vertical direction, the corresponding words in the double Weyl group generators do not coincide, and even correspond to the different groups. The first one belongs to the double Weyl group of $\widehat{PSL}(2)$, while in the second case one gets $\widehat{PSL}(3)$. This is of course the well-known fact about Toda chains, which can be described in terms of $2\times 2$ or $N\times N$ Lax matrices, $N=3$ for this example, (see \cite{AMJGP} for more details and references), and here it comes from possibility to cut torus $\Sigma$ in two different ways.

The coordinates on cohomologies of the torus $H^1(\Sigma)$ are $\lambda$ and $\mu$, they correspond to the cycles, shown on fig.~\ref{fi:toda23}C. The cohomology generators $\alpha,\beta,\gamma,\delta$ of the dual surface $H^1(\widetilde{\Sigma})$ correspond to the cycles, shown on fig.~\ref{fi:toda23}D.

The Poisson bracket between the face variables, following from the Thurston diagrams on fig.~\ref{fi:toda23}B or fig.~\ref{fi:toda23}B' are:
\be
\{y,z\}=2yz,\quad \{v,x\}=2vx,\quad \{w,u\}=2yz,
\\
\{z,v\}=zv,\quad\{z,w\}=zw,\quad
\{x,y\}=xy,
\\
 \{x,w\}=xw,\quad
\{u,y\}=uy,\quad \{u,w\}=uw,
\ee
while all others pairs of variables Poisson commute. This Poisson structure can be encoded by the exchange graph dual to the bipartite graph on fig.~\ref{fi:toda23}C, or by the Cartan matrix of the affine $\widehat{PSL}(3)$.

\begin{figure}[H]
\begin{center}
\begin{tikzpicture}[scale=0.8]
\draw (0,0) -- +(0,2);
\draw (1,0) -- +(0,2);
\draw (2,0) -- +(0,2);
\draw (0,0) -- +(3,0);
\draw (0,1) -- +(3,0);
\draw[very thick] (0,0)--(0,1);\draw[very thick] (1,1)--(1,2);\draw[very thick] (2,0)--(2,1);
\fill (0,0)+(0,0) circle (0.08);
\fill (0,0)+(1,1) circle (0.08);
\fill (0,0)+(2,0) circle (0.08);
\fill (0,0)+(0,1) circle (0.08);\fill[color=white] (0,1) circle (0.06);
\fill (0,0)+(1,0) circle (0.08);\fill[color=white] (1,0) circle (0.06);
\fill (0,0)+(2,1) circle (0.08);\fill[color=white] (2,1) circle (0.06);
\node at (1.5,-0.5) {$\mu^{-2}$};
\end{tikzpicture}\quad
\begin{tikzpicture}[scale=0.8]
\draw (0,0) -- +(0,2);
\draw (1,0) -- +(0,2);
\draw (2,0) -- +(0,2);
\draw (0,0) -- +(3,0);
\draw (0,1) -- +(3,0);
\draw[very thick] (0,0)--(0,1);\draw[very thick] (1,0)--(1,1);\draw[very thick] (2,0)--(2,1);
\fill (0,0)+(0,0) circle (0.08);
\fill (0,0)+(1,1) circle (0.08);
\fill (0,0)+(2,0) circle (0.08);
\fill (0,0)+(0,1) circle (0.08);\fill[color=white] (0,1) circle (0.06);
\fill (0,0)+(1,0) circle (0.08);\fill[color=white] (1,0) circle (0.06);
\fill (0,0)+(2,1) circle (0.08);\fill[color=white] (2,1) circle (0.06);
\node at (1.5,-0.5) {$\alpha\gamma^{-1}\mu^{-1}$};
\end{tikzpicture}\quad
\begin{tikzpicture}[scale=0.8]
\draw (0,0) -- +(0,2);
\draw (1,0) -- +(0,2);
\draw (2,0) -- +(0,2);
\draw (0,0) -- +(3,0);
\draw (0,1) -- +(3,0);
\draw[very thick] (0,1)--(0,2);\draw[very thick] (1,0)--(1,1);\draw[very thick] (2,1)--(2,2);
\fill (0,0)+(0,0) circle (0.08);
\fill (0,0)+(1,1) circle (0.08);
\fill (0,0)+(2,0) circle (0.08);
\fill (0,0)+(0,1) circle (0.08);\fill[color=white] (0,1) circle (0.06);
\fill (0,0)+(1,0) circle (0.08);\fill[color=white] (1,0) circle (0.06);
\fill (0,0)+(2,1) circle (0.08);\fill[color=white] (2,1) circle (0.06);
\node at (1.5,-0.5) {$\mu$};
\end{tikzpicture}\quad
\begin{tikzpicture}[scale=0.8]
\draw (0,0) -- +(0,2);
\draw (1,0) -- +(0,2);
\draw (2,0) -- +(0,2);
\draw (0,0) -- +(3,0);
\draw (0,1) -- +(3,0);
\draw[very thick] (0,1)--(0,2);\draw[very thick] (1,1)--(1,2);\draw[very thick] (2,1)--(2,2);
\fill (0,0)+(0,0) circle (0.08);
\fill (0,0)+(1,1) circle (0.08);
\fill (0,0)+(2,0) circle (0.08);
\fill (0,0)+(0,1) circle (0.08);\fill[color=white] (0,1) circle (0.06);
\fill (0,0)+(1,0) circle (0.08);\fill[color=white] (1,0) circle (0.06);
\fill (0,0)+(2,1) circle (0.08);\fill[color=white] (2,1) circle (0.06);
\node at (1.5,-0.5) {$\gamma\alpha^{-1}$};
\end{tikzpicture}
\end{center}
\begin{center}
\begin{tikzpicture}[scale=0.8]
\draw (0,0) -- +(0,2);
\draw (1,0) -- +(0,2);
\draw (2,0) -- +(0,2);
\draw (0,0) -- +(3,0);
\draw (0,1) -- +(3,0);
\draw[very thick] (0,0)--(0,1);\draw[very thick] (1,0)--(1,1);\draw[very thick] (2,1)--(2,2);
\fill (0,0)+(0,0) circle (0.08);
\fill (0,0)+(1,1) circle (0.08);
\fill (0,0)+(2,0) circle (0.08);
\fill (0,0)+(0,1) circle (0.08);\fill[color=white] (0,1) circle (0.06);
\fill (0,0)+(1,0) circle (0.08);\fill[color=white] (1,0) circle (0.06);
\fill (0,0)+(2,1) circle (0.08);\fill[color=white] (2,1) circle (0.06);
\node at (1.5,-0.5) {$\alpha$};
\end{tikzpicture}\quad
\begin{tikzpicture}[scale=0.8]
\draw (0,0) -- +(0,2);
\draw (1,0) -- +(0,2);
\draw (2,0) -- +(0,2);
\draw (0,0) -- +(3,0);
\draw (0,1) -- +(3,0);
\draw[very thick] (0,0)--(0,1);\draw[very thick] (1,1)--(1,2);\draw[very thick] (2,1)--(2,2);
\fill (0,0)+(0,0) circle (0.08);
\fill (0,0)+(1,1) circle (0.08);
\fill (0,0)+(2,0) circle (0.08);
\fill (0,0)+(0,1) circle (0.08);\fill[color=white] (0,1) circle (0.06);
\fill (0,0)+(1,0) circle (0.08);\fill[color=white] (1,0) circle (0.06);
\fill (0,0)+(2,1) circle (0.08);\fill[color=white] (2,1) circle (0.06);
\node at (1.5,-0.5) {$\mu^{-1}\gamma$};
\end{tikzpicture}\quad
\begin{tikzpicture}[scale=0.8]
\draw (0,0) -- +(0,2);
\draw (1,0) -- +(0,2);
\draw (2,0) -- +(0,2);
\draw (0,0) -- +(3,0);
\draw (0,1) -- +(3,0);
\draw[very thick] (0,1)--(0,2);\draw[very thick] (1,0)--(1,1);\draw[very thick] (2,0)--(2,1);
\fill (0,0)+(0,0) circle (0.08);
\fill (0,0)+(1,1) circle (0.08);
\fill (0,0)+(2,0) circle (0.08);
\fill (0,0)+(0,1) circle (0.08);\fill[color=white] (0,1) circle (0.06);
\fill (0,0)+(1,0) circle (0.08);\fill[color=white] (1,0) circle (0.06);
\fill (0,0)+(2,1) circle (0.08);\fill[color=white] (2,1) circle (0.06);
\node at (1.5,-0.5) {$\gamma^{-1}$};
\end{tikzpicture}\quad
\begin{tikzpicture}[scale=0.8]
\draw (0,0) -- +(0,2);
\draw (1,0) -- +(0,2);
\draw (2,0) -- +(0,2);
\draw (0,0) -- +(3,0);
\draw (0,1) -- +(3,0);
\draw[very thick] (0,1)--(0,2);\draw[very thick] (1,1)--(1,2);\draw[very thick] (2,0)--(2,1);
\fill (0,0)+(0,0) circle (0.08);
\fill (0,0)+(1,1) circle (0.08);
\fill (0,0)+(2,0) circle (0.08);
\fill (0,0)+(0,1) circle (0.08);\fill[color=white] (0,1) circle (0.06);
\fill (0,0)+(1,0) circle (0.08);\fill[color=white] (1,0) circle (0.06);
\fill (0,0)+(2,1) circle (0.08);\fill[color=white] (2,1) circle (0.06);
\node at (1.5,-0.5) {$\alpha^{-1}\mu^{-1}$};
\end{tikzpicture}
\end{center}
\begin{center}
\begin{tikzpicture}[scale=0.8]
\draw (0,0) -- +(0,2);
\draw (1,0) -- +(0,2);
\draw (2,0) -- +(0,2);
\draw (0,0) -- +(3,0);
\draw (0,1) -- +(3,0);
\draw[very thick] (0,0)--(1,0);\draw[very thick] (1,1)--(2,1);\draw[very thick] (2,0)--(3,0);
\fill (0,0)+(0,0) circle (0.08);
\fill (0,0)+(1,1) circle (0.08);
\fill (0,0)+(2,0) circle (0.08);
\fill (0,0)+(0,1) circle (0.08);\fill[color=white] (0,1) circle (0.06);
\fill (0,0)+(1,0) circle (0.08);\fill[color=white] (1,0) circle (0.06);
\fill (0,0)+(2,1) circle (0.08);\fill[color=white] (2,1) circle (0.06);
\node at (1.5,-0.5) {$-\lambda\mu^{-1}$};
\end{tikzpicture}\quad
\begin{tikzpicture}[scale=0.8]
\draw (0,0) -- +(0,2);
\draw (1,0) -- +(0,2);
\draw (2,0) -- +(0,2);
\draw (0,0) -- +(3,0);
\draw (0,1) -- +(3,0);
\draw[very thick] (0,1)--(1,1);\draw[very thick] (1,0)--(2,0);\draw[very thick] (2,1)--(3,1);
\fill (0,0)+(0,0) circle (0.08);
\fill (0,0)+(1,1) circle (0.08);
\fill (0,0)+(2,0) circle (0.08);
\fill (0,0)+(0,1) circle (0.08);\fill[color=white] (0,1) circle (0.06);
\fill (0,0)+(1,0) circle (0.08);\fill[color=white] (1,0) circle (0.06);
\fill (0,0)+(2,1) circle (0.08);\fill[color=white] (2,1) circle (0.06);
\node at (1.5,-0.5) {$c\lambda^{-1}$};
\end{tikzpicture}\quad
\begin{tikzpicture}[scale=0.8]
\draw (0,0) -- +(0,2);
\draw (1,0) -- +(0,2);
\draw (2,0) -- +(0,2);
\draw (0,0) -- +(3,0);
\draw (0,1) -- +(3,0);
\draw[very thick] (0,0)--(1,0);\draw[very thick] (0,1)--(1,1);\draw[very thick] (2,0)--(2,1);
\fill (0,0)+(0,0) circle (0.08);
\fill (0,0)+(1,1) circle (0.08);
\fill (0,0)+(2,0) circle (0.08);
\fill (0,0)+(0,1) circle (0.08);\fill[color=white] (0,1) circle (0.06);
\fill (0,0)+(1,0) circle (0.08);\fill[color=white] (1,0) circle (0.06);
\fill (0,0)+(2,1) circle (0.08);\fill[color=white] (2,1) circle (0.06);
\node at (1.5,-0.5) {$\beta\mu^{-1}$};
\end{tikzpicture}\quad
\begin{tikzpicture}[scale=0.8]
\draw (0,0) -- +(0,2);
\draw (1,0) -- +(0,2);
\draw (2,0) -- +(0,2);
\draw (0,0) -- +(3,0);
\draw (0,1) -- +(3,0);
\draw[very thick] (0,0)--(1,0);\draw[very thick] (0,1)--(1,1);\draw[very thick] (2,1)--(2,2);
\fill (0,0)+(0,0) circle (0.08);
\fill (0,0)+(1,1) circle (0.08);
\fill (0,0)+(2,0) circle (0.08);
\fill (0,0)+(0,1) circle (0.08);\fill[color=white] (0,1) circle (0.06);
\fill (0,0)+(1,0) circle (0.08);\fill[color=white] (1,0) circle (0.06);
\fill (0,0)+(2,1) circle (0.08);\fill[color=white] (2,1) circle (0.06);
\node at (1.5,-0.5) {$\beta\gamma$};
\end{tikzpicture}
\end{center}
\begin{center}
\begin{tikzpicture}[scale=0.8]
\draw (0,0) -- +(0,2);
\draw (1,0) -- +(0,2);
\draw (2,0) -- +(0,2);
\draw (0,0) -- +(3,0);
\draw (0,1) -- +(3,0);
\draw[very thick] (0,0)--(0,1);\draw[very thick] (1,0)--(2,0);\draw[very thick] (1,1)--(2,1);
\fill (0,0)+(0,0) circle (0.08);
\fill (0,0)+(1,1) circle (0.08);
\fill (0,0)+(2,0) circle (0.08);
\fill (0,0)+(0,1) circle (0.08);\fill[color=white] (0,1) circle (0.06);
\fill (0,0)+(1,0) circle (0.08);\fill[color=white] (1,0) circle (0.06);
\fill (0,0)+(2,1) circle (0.08);\fill[color=white] (2,1) circle (0.06);
\node at (1.5,-0.5) {$c\alpha\delta\mu^{-1}$};
\end{tikzpicture}\quad
\begin{tikzpicture}[scale=0.8]
\draw (0,0) -- +(0,2);
\draw (1,0) -- +(0,2);
\draw (2,0) -- +(0,2);
\draw (0,0) -- +(3,0);
\draw (0,1) -- +(3,0);
\draw[very thick] (0,1)--(0,2);\draw[very thick] (1,0)--(2,0);\draw[very thick] (1,1)--(2,1);
\fill (0,0)+(0,0) circle (0.08);
\fill (0,0)+(1,1) circle (0.08);
\fill (0,0)+(2,0) circle (0.08);
\fill (0,0)+(0,1) circle (0.08);\fill[color=white] (0,1) circle (0.06);
\fill (0,0)+(1,0) circle (0.08);\fill[color=white] (1,0) circle (0.06);
\fill (0,0)+(2,1) circle (0.08);\fill[color=white] (2,1) circle (0.06);
\node at (1.5,-0.5) {$c\delta$};
\end{tikzpicture}\quad
\begin{tikzpicture}[scale=0.8]
\draw (0,0) -- +(0,2);
\draw (1,0) -- +(0,2);
\draw (2,0) -- +(0,2);
\draw (0,0) -- +(3,0);
\draw (0,1) -- +(3,0);
\draw[very thick] (1,0)--(1,1);\draw[very thick] (2,0)--(3,0);\draw[very thick] (2,1)--(3,1);
\fill (0,0)+(0,0) circle (0.08);
\fill (0,0)+(1,1) circle (0.08);
\fill (0,0)+(2,0) circle (0.08);
\fill (0,0)+(0,1) circle (0.08);\fill[color=white] (0,1) circle (0.06);
\fill (0,0)+(1,0) circle (0.08);\fill[color=white] (1,0) circle (0.06);
\fill (0,0)+(2,1) circle (0.08);\fill[color=white] (2,1) circle (0.06);
\node at (1.5,-0.5) {$\beta^{-1}\gamma^{-1}\delta^{-1}$};
\end{tikzpicture}\quad
\begin{tikzpicture}[scale=0.8]
\draw (0,0) -- +(0,2);
\draw (1,0) -- +(0,2);
\draw (2,0) -- +(0,2);
\draw (0,0) -- +(3,0);
\draw (0,1) -- +(3,0);
\draw[very thick] (1,1)--(1,2);\draw[very thick] (2,0)--(3,0);\draw[very thick] (2,1)--(3,1);
\fill (0,0)+(0,0) circle (0.08);
\fill (0,0)+(1,1) circle (0.08);
\fill (0,0)+(2,0) circle (0.08);
\fill (0,0)+(0,1) circle (0.08);\fill[color=white] (0,1) circle (0.06);
\fill (0,0)+(1,0) circle (0.08);\fill[color=white] (1,0) circle (0.06);
\fill (0,0)+(2,1) circle (0.08);\fill[color=white] (2,1) circle (0.06);
\node at (1.5,-0.5) {$\alpha^{-1}\beta^{-1}\delta^{-1}\mu^{-1}$};
\end{tikzpicture}
\end{center}
\caption{Dimer configurations and pairings $(-1)^{Q(D)}\langle \boldsymbol{A},D\rangle$.\label{fi:toda23dimers}}
\end{figure}
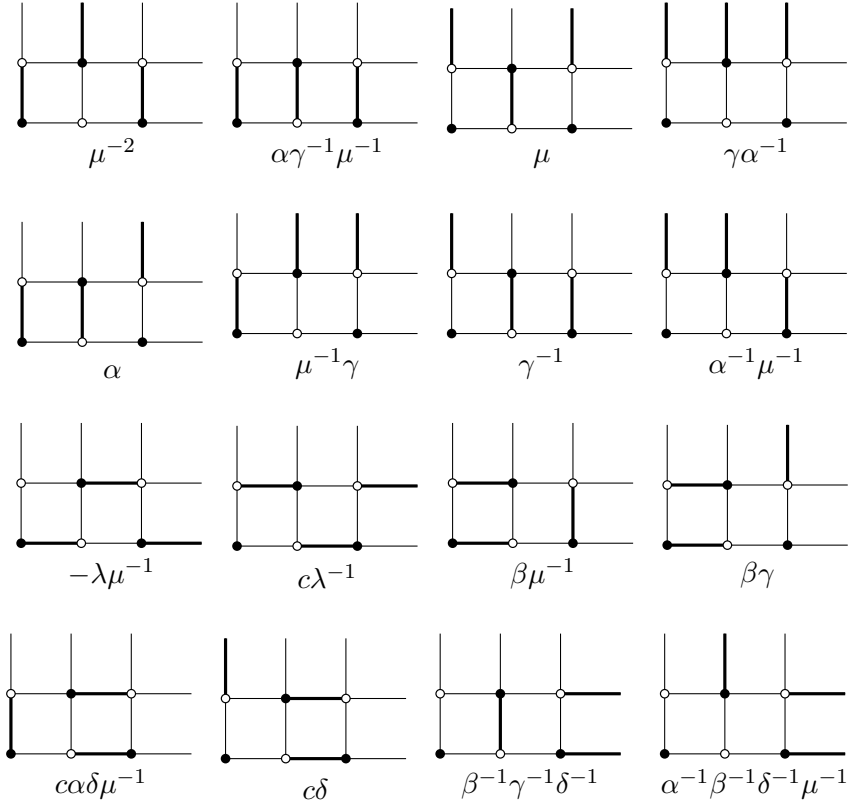

The are 16 dimer configurations, shown on fig.~\ref{fi:toda23dimers}, summed up to the
double partition function:
\be
\label{SdToda2}
S^d= \mu^{-2} + \mu + \mu^{-1}(\alpha\gamma^{-1}+\gamma +\alpha^{-1}+\beta+c\alpha\delta+\alpha^{-1}\beta^{-1}\delta^{-1})+
\\
+(\gamma\alpha^{-1}+\alpha+\gamma^{-1}+\beta\gamma+c\delta+
\gamma^{-1}\beta^{-1}\delta^{-1})+c\lambda^{-1}-\lambda\mu^{-1}
\ee
Comparing fig.~\ref{fi:toda23}C and fig.~\ref{fi:toda23}E one finds, that the face variables are related to the cohomology variables by
$$
x=c^{-1}\alpha^{-1}\gamma\delta^{-1},\quad y=\alpha\beta,\quad z=\alpha\beta^{-1}\gamma^{-1},\quad u=\beta\delta,\quad v=c\gamma\delta,\quad w=(\alpha\beta\gamma\delta)^{-1}.
$$
Solving these relations for the latter, and substituting the result into \rf{SdToda2}, we get the face partition function
\be
\label{SfToda2}
S^f=\mu^{-2} + \mu + \mu^{-1}c^{1/3}(u^2y^2zv)^{-1/3}(1+u+uy+uyz+uyzv+c^{-1}yv)+
\\
+c^{1/3}(xyu^2v^2)^{-1/3}(1+u+uv+uvx+uvxy+c^{-1}uv)+ c\lambda^{-1}-\lambda\mu^{-1}
\ee
On the other hand, the same partition function is given up to normalisation by the characteristic polynomial, which can be read off the word from fig.~\ref{fi:toda23}B. Computing the $2\times 2$ determinant after some convenient location of the shift operators, one gets the expression
\begin{equation}\label{sf4}
\begin{array}{c}
\det(H_1(x)E_1H_1(y)E_{\bar{1}}H_0(z)E_0(\lambda)H_0(w)E_{\bar{0}}(\lambda)H_1(u)E_1H_1(v)E_{\bar{1}}\Lambda -\mu)=\\[7pt]
=\mu^2-\lambda\mu x(1+u+uy+uyz+uyzv+uyzvx)-\\[7pt]-\mu(1+u+uv+uvx+uvxy+uvxyz)-\lambda^2\mu x^2yzu-\lambda^{-1}\mu uv-\lambda xuyv,
\end{array}
\end{equation}
which reduces to (\ref{SfToda2}) under the substitutions
$$\lambda\mapsto \mu^{-1}x^{-2/3}y^{-1/3}z^{-1/3}v^{1/3},\quad \mu\mapsto \lambda\mu^{-1}x^{2/3}y^{1/3}z^{1/3}uv^{2/3}$$
and multiplication the whole expression in \rf{sf4} by $-\lambda^{-1}\mu x^{-4/3}y^{-2/3}z^{-2/3}u^{-2}v^{-4/3}$.

The same face partition function can be obtained from the characteristic polynomial of the $3\times 3$ matrix, read off the word from fig.~\ref{fi:toda23}B'. Computing the determinant, after pushing all
shift operators to the right, one gets:
\be
\label{Toda2pgl3}
-\det(H_0(w)E_0(\lambda)H_0(u)E_{\bar{0}}(\lambda)H_1(y)E_1H_1(z)E_{\bar{1}}H_2(v)E_2H_2(x)E_{\bar{2}} -\mu)=
\\
= \mu^3-vwy\mu^2\lambda-(1+u+uv+uvx+uvxy+uvxyz)u^{-1}\mu^2+
\\
+(1+u+uy+uyz+uyzv+uyzvx)vxu^{-1}\mu
-vxu^{-1}w^{-1}\mu\lambda^{-1}-v^2x^2yz
\ee
After the substitutions
$$ \mu\mapsto -\mu x^{2/3}y^{1/3}z^{1/3}v^{2/3},\quad \lambda\mapsto \lambda^{-1}x^{2/3}y^{1/3}z^{1/3}v^{2/3}$$
and multiplication expression in \rf{Toda2pgl3} by $\mu^{-2}x^{-2}y^{-1}z^{-1}v^{-2}$, one comes back to \rf{SfToda2}.

The symmetry group $\mathcal{G}_\Delta' = \mathcal{G}_\Delta$ in this case has rank one and torsion $\mathbb{Z}/3\mathbb{Z}$. The torsion generator $(0,1,1,0)$ is the permutation $(x,y,z,u,v,w)\mapsto(u,v,x,z,w,y)$. The infinite order generator $(0,1,0,0)$ is the transformation $\tau$ given by
$$(x,y,z,u,v,w)\mapsto \left(y\frac{(1+x)(1+u)}{(1+z^{-1})^2},x^{-1},v\frac{(1+u)(1+z)}{(1+x^{-1})^2},u^{-1},w\frac{(1+z)(1+x)}{(1+u^{-1})^2},z^{-1}\right)$$

\subsection{Parallelograms of arbitrary size and the pentagram map
\label{ss:pentagram}}

Finally we consider a subclass of integrable systems with Newton polygons being parallelograms of arbitrary size.
For $G=\widehat{PSL}(N)$ the word $s_0s_{\bar{0}}s_1s_{\bar{1}}\ldots s_{l-1}s_{\overline{l-1}}\Lambda^{l}$, where the indices are understood modulo $N$, is cyclically reduced and the corresponding Newton polygon is a parallelogram with vertices $(0,0),(1,0),(N,l),(N-1,l)$ of double area $2l$. In this way one can obtain any parallelogram with at leas one primitive side.

Our main example here will be parallelogram of width 3, i.e. just a little wider than the parallelograms, giving rise to relativistic Toda lattices. The main reason for considering such example is that the discrete evolution for this system, discovered in \cite{Schwartz} and studied in \cite{OT} was one of the motivations for our work.

\begin{figure}[H]
\begin{center}
\begin{tikzpicture}
\draw (-0.1,0.5)--(-0.5,2);
\draw (-0.5,0.7)--(0.5,2.7);
\draw (-0.5,1.7)--(1.8,3.1);
\draw (0,2.5)--(3,3);
\draw (1.3,3.1)--(4.7,2.2);
\draw (2.5,3)--(5.6,1.6);
\draw (4,2.7)--(5.6,1.2);
\draw (5.4,1.8)--(5.2,0.5);
\fill (-0.28,1.15) circle (0.06); \node [right] at (-0.28,1.15) {$p_0$};
\fill (0.22,2.15) circle (0.06); \node [right] at (0.2,2) {$p_1$};
\fill (1.15,2.7) circle (0.06); \node [below] at (1.3,2.7) {$p_2$};
\fill (2.2,2.86) circle (0.06); \node [below] at (2.2,2.86) {$p_3$};
\fill (3.65,2.48) circle (0.06); \node [below] at (3.65,2.48) {$p_4$};
\fill (4.78,1.95) circle (0.06); \node [below] at (4.78,1.95) {$p_5$};
\fill (5.35,1.45) circle (0.06); \node [left] at (5.45,1.35) {$p_6$};

\fill (-0.43,1.74) circle (0.06); \node [left] at (-0.43,1.74) {$\tau(p_1)$};
\fill (0.435,2.573) circle (0.06); \node [above left] at (0.5,2.45) {$\tau(p_2)$};
\fill (1.64,3.01) circle (0.06); \node [above] at (1.6,3) {$\tau(p_3)$};
\fill (2.64,2.93) circle (0.06); \node [above] at (2.7,2.9) {$\tau(p_4)$};
\fill (4.47,2.25) circle (0.06); \node [above right] at (4.4,2.15) {$\tau(p_5)$};
\fill (5.38,1.7) circle (0.06); \node [right] at (5.35,1.8) {$\tau(p_6)$};
\end{tikzpicture}
\end{center}
\caption{Pentagram map.}
\label{fi:pentagram}
\end{figure}
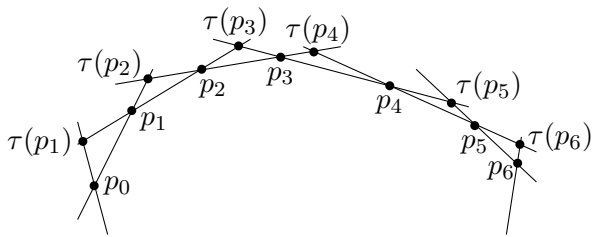

Recall briefly the construction of the discrete integrable system by R.~Schwartz.
Consider the space of sequences $\dotsc,p_{-1},p_0,p_1,\dotsc$ of points in projective plane. This space has a (birational) automorphism $\tau$ defined by a simple geometric construction, shown on fig.~\ref{fi:pentagram}. Namely, $\tau(p_i)$ is the intersection point of the lines $(p_{i-1},p_i)$ and $(p_{i+1},p_{i+2})$.

The group $PGL(3)$ acts on this space by projective transformations. The space of such sequences contains a subquotient of projective classes of $l$-quasiperiodic sequences~\footnote{By  an $l$-quasiperiodic sequence we mean a sequence $\dotsc, p_0,p_1,\dotsc$ such that the shifted sequence $\dotsc,p_l,p_{l+1},\dotsc$ is projectively isomorphic to it.},
that are finite dimensional and which we shall denote by $\mathcal{P}_l$. The space $\mathcal{P}_l$ obviously has dimension $2l$ for any $l\in \mathbb{N}$.

Our aim is to show that $\mathcal{P}_l$ is canonically isomorphic to the double Bruhat cell of $\widehat{PGL}(3)$ corresponding to the word $u=(\prod_{i=1}^ls_is_{\bar{i}})\Lambda^l\in (\widehat{W}\times\widehat{W})^\sharp$, where all indices are understood modulo 3. Moreover we will check, that the pentagram map $\tau$ is just a particular case of the discrete flow, described in sect.~\ref{ss:automorph}.

For a generic sequence of points $\dotsc,p_1,p_2\dotsc$ in $\mathbb{RP}^2$ one can define its lift $\dotsc,e_1,e_2,\dotsc$ to a collection of vectors of $\mathbb{R}^3$ by the requirement that $e_{i}+e_{i+1}$ belongs to the plane generated by $e_{i+2}$ and $e_{i+3}$, i.e. there exist coefficients $a_i$ and $b_i$, such that
\begin{equation}\label{normalization}
e_{i+3}=a_i(e_{i+1}+e_i)+b_ie_{i+2}.
\end{equation}
The transition matrix from the basis $e_{i},e_{i+1},e_{i+2}$ to the basis $e_{i+1},e_{i+2},e_{i+3}$ can be written as
\be
\label{Mmatr}
 M_i = \begin{pmatrix}
0&0&a_i\\
1&0&a_i\\
0&1&b_i
\end{pmatrix}
\ee
Recall the R.~Schwartz coordinates \cite{Schwartz} on the space of sequences of points on a projective plane\footnote{Our indexing convention is shifted from \cite{Schwartz} by two.}. For every $i$ let $x_i$ and $y_i$ be the cross-ratios of quadruples of collinear points:
$$x_i=[p_i,p_{i+1},(p_i,p_{i+1})\cap (p_{i+2},p_{i+3}),(p_i,p_{i+1})\cap (p_{i+3},p_{i+4})],$$
$$y_i=[(p_{i+4},p_{i+3})\cap (p_{i},p_{i+1}),(p_{i+4},p_{i+3})\cap (p_{i+1},p_{i+2}),p_{i+3},p_{i+4}].$$
Here $(p_i,p_j)$ denotes a line through $p_i$ and $p_j$; the cross-ratio is chosen to be
$$[a,b,c,d]=\frac{(a-b)(c-d)}{(a-c)(b-d)}$$
The Schwartz coordinates $x_i$ and $y_i$  can be calculated (see \rf{cross-ratios} in Appendix~\ref{ap:schwartz}) as
\begin{equation}\label{xyab}x_{i}= b_i/a_i,\quad\quad y_{i} = -a_{i+1}/b_ib_{i+1}
\end{equation}
Conversely
\begin{equation}
\label{abxy}a_i=-(x_ix_{i+1}y_i)^{-1},\quad\quad b_i = -(x_{i+1}y_i)^{-1}
\end{equation}
so that expression \rf{Mmatr} in the Schwartz coordinates, up to a scalar multiple is
\begin{equation}
\label{minschwartz}
M_i= \begin{pmatrix}
0&0&x_i^{-1}\\
-y_ix_{i+1}&0&x_i^{-1}\\
0&-y_ix_{i+1}&1
\end{pmatrix}
\end{equation}
The transition matrix from $e_{k},e_{k+1},e_{k+2}$ to $e_{m},e_{m+1},e_{m+2}$ is therefore given by the product
\begin{equation}
\label{monodromy}
M_{k\rightarrow m}=M_kM_{k+1}\cdots M_{m-1}.
\end{equation}
of the matrices \rf{minschwartz}.

Our key observation is that the transition matrix $M_{k\to m}$ can be written as a simple expression of generators of the group $\widehat{PSL}^\sharp(3,\mathbb R)$. Namely, consider the product
\be
\label{hatM}
\hat{M}(u,v,\lambda)=H_0(u)E_0 \Lambda E_{\bar{0}}  H_2(v)=
\\
= T_u
\begin{pmatrix}
1&0&0\\
0&1&0\\
\lambda&0&1
  \end{pmatrix}
\begin{pmatrix}
0&0&\lambda^{-1}\\
1&0&0\\
0&1&0
  \end{pmatrix}
  \begin{pmatrix}
1&0&\lambda^{-1}\\
0&1&0\\
0&0&1
  \end{pmatrix}
  \begin{pmatrix}
v&0&0\\
0&v&0\\
0&0&1
  \end{pmatrix}T_v=
\\
=\begin{pmatrix}
0&0&\lambda^{-1}u^{-1}\\
v&0&\lambda^{-1}u^{-1}\\
0&v&1
  \end{pmatrix}T_{uv}=
M(u,v,\lambda)T_{uv}
\ee
The matrix part of this expression
is directly related to the matrix \rf{minschwartz}
$$
M_i = M(a_i/b_i,1/b_i,1)=M(x_i,y_ix_{i+1},1)
$$
and has the obvious property:
$$
T_w M(u,v,\lambda)= M(uw,v,\lambda)T_w
$$
Hence, the product $\hat{M}(u_k,v_k,\lambda)\cdots \hat{M}(u_{m-1},v_{m-1},\lambda)$ can be transformed by taking all shift operators to the right:
$$
\hat{M}(u_k,v_k,\lambda)\cdots \hat{M}(u_{m-1},v_{m-1},\lambda)=
$$
$$
=M(u_k,v_k,\lambda)M(u_{k+1}v_ku_k,v_{k+1},\lambda)\cdots M(u_{m-1}v_{m-2}u_{m-2}\cdots v_ku_k,v_{m-1},\lambda)T_{u_kv_k\cdots u_{m-1}v_{m-1}}.
$$
Comparing this expression to the product (\ref{monodromy}), with substituted (\ref{minschwartz}),
$$ M_{k\to m} = M(x_k,-y_kx_{k+1},1)\cdots M(x_m,-y_mx_{m+1},1)$$
one finds, that they coincide under identification
\begin{equation}\label{schwortztoglick}
u_{i}=-(x_{i-1}y_{i-1})^{-1},\quad v_i = -y_ix_{i+1}.
\end{equation}
The variables $u_i$ and $v_i$ are coordinates on the space of sequences of points, introduced by M.~Glick in \cite{Glick}, where it is shown than the pentagram map in this coordinates is a sequence of mutations.

On the other hand the product $\hat{M}(u_k,v_k)\cdots \hat{M}(u_{m-1},v_{m-1})$ can be transformed by taking all the automorphisms $\Lambda$ to the right:
\be
\label{mon3}
\hat{M}(u_k,v_k)\cdots \hat{M}(u_{m-1},v_{m-1}) = H_0(u_k)E_0 \Lambda E_{\bar{0}}  H_{-1}(v_k)\cdots H_0(u_{m-1}) E_0 \Lambda E_{\bar{0}}  H_{-1}(v_{m-1})=
\\
=
H_0(u_k) E_0 E_{\bar{1}}  H_{0}(v_k)\cdots H_{m-k-1}(u_j)E_{m-k-1} E_{\overline{m-k}} H_{m-k-1}(v_{m-1})\Lambda^{m-k}
\ee
Therefore this expression is the standard expression for the map \rf{laxmap}, corresponding to the word $u=s_0s_{\bar{1}}s_1s_{\bar{2}}s_{2}\cdots s_{m-k-1} s_{\overline{m-k}}\Lambda^{m-k}$, which is obviously cyclically reduced. According to the standard rules from the appendix \ref{ap:exchange} the Poisson brackets coincides with that of \cite{Glick}:
\be
\label{pbglick}
\{u_j,v_i\}=u_jv_i \mbox{ if } j=i+2 \mbox{ or } j=i-1,
\\
\{u_j,v_i\}=-u_jv_i \mbox{ if } j=i \mbox{ or } j=i+1
\ee
and all the other brackets vanish. The characteristic polynomial of \rf{mon3} gives rise to a spectral curve with the Newton polygon being a parallelogram of width 3 (matrices $3\times 3$) and arbitrary length, see \cite{AMJGP}.

\paragraph{Pentagram map.} Consider now the pentagram map $\tau$. According to the definition of the map and the normalisation (\ref{normalization}) of our sequence of vectors $\{e_i\}$, the evolved basis $\{\tilde{e}_i=\tau(e_i)\}$ should be given by the simple formula $\tilde{e}_{i}=\alpha_i(e_{i+1}+e_i)$, with some nonzero constants $\alpha_i$. On the other hand the new ordinates $\tilde{a}_i$ and $\tilde{b}_i$ of the sequence are defined by the relation $$\tilde{e}_{i+3}=\tilde{a}_i(\tilde{e}_{i+1}+\tilde{e}_i)+\tilde{b}_i\tilde{e}_{i+2}.$$ Putting these relations together one gets an explicit expression for the pentagram map (see Appendix \ref{ap:schwartz}):
\begin{equation}\label{evolution}
\begin{array}{rcl}
\alpha_i &=&\displaystyle\frac{a_i}{b_i+1}\\[10pt]
a_i &\mapsto& \tilde{a}_i = \displaystyle\frac{a_{i+3}(b_{i+1}+1)}{(b_{i+3}+1)},\\[10pt]
b_i &\mapsto&\tilde{b}_i = \displaystyle\frac{b_i(b_{i+1}+1)(b_{i+2}+1)a_{i+3}}{(b_{i}+1)(b_{i+3}+1)a_{i+2}}.
\end{array}
\end{equation}
Using (\ref{xyab}) and (\ref{abxy}) one gets in the Schwartz coordinates
\begin{equation}\label{evolutionSchwartz}
\begin{array}{rcl}
x_i &\mapsto& x_{i+2}\displaystyle\frac{1-x_{i+3}y_{i+2}}{1-x_{i+1}y_{i}}\\[10pt]
y_i &\mapsto& y_{i+1}\displaystyle\frac{1-x_{i+1}y_{i}}{1-x_{i+3}y_{i+2}}.
\end{array}
\end{equation}
what coincides, up to an index shift, with \cite{Schwartz} (see  page 523) or formula (2.5) from \cite{OT}.

Expressing this evolution in terms of cluster coordinates $u_i$ and $v_i$ we get
\begin{equation}\label{evolutionuv}
\begin{array}{rcl}
u_{i} &\mapsto& v_{i-1}^{-1},\\[5pt]
v_{i} &\mapsto& u_{i+2}(1+v_{i})(1+v_{i+3})(1+v_{i+1}^{-1})^{-1}(1+v_{i+2}^{-1})^{-1}.
\end{array}
\end{equation}
This transformation is a composition of a cluster mutation
$$
\begin{array}{rcl}
v_{i} &\mapsto& v_{i}^{-1},\\[5pt]
u_{i} &\mapsto& u_{i}(1+v_{i-2})(1+v_{i+1})(1+v_{i-1}^{-1})^{-1}(1+v_{i}^{-1})^{-1}.
\end{array}
$$
in all variables $v_i$ composed with the renumbering
$$
\begin{array}{rcl}
v_{i} &\mapsto& u_{i+2}\\[8pt]
u_{i} &\mapsto& v_{i-1}.
\end{array}
$$

\paragraph{Other dimensions.} Observe, that the construction works if we replace $\widehat{PGL(3)}$ by  $\widehat{PGL(N)}$ with any $N\geqslant 2$.
For generic $N>2$ the operator \rf{hatM} becomes
$$\hat{M}=T_u
\begin{pmatrix}
0&0&0&\cdots&0&\lambda^{-1}u\\
v&0&0&\cdots&0&\lambda^{-1}u\\
0&v&0&\cdots&0&0\\
\vdots&&&\ddots&\vdots&\vdots\\
0&0&0&\cdots&0&0\\
0&0&0&\cdots&v&1
\end{pmatrix}
$$
The Poisson brackets are just analogous to \rf{pbglick}:
$$\{u_j,v_i\}=u_jv_i \mbox{ if } j=i+N-1 \mbox{ or } j=i-1,$$
$$\{u_j,v_i\}=-u_jv_i \mbox{ if } j=i \mbox{ or } j=i+N-2.$$
and zero otherwise.
For $N=2$ the operator $\hat{M}$ becomes:
$$\hat{M}=T_u
\begin{pmatrix}1&0\\\lambda&1\end{pmatrix}
\begin{pmatrix}0&\lambda^{-1}\\1&0\end{pmatrix}
\begin{pmatrix}1&\lambda^{-1}\\0&1\end{pmatrix}
\begin{pmatrix}v&0\\0&1\end{pmatrix}T_v=
\begin{pmatrix}0&\lambda^{-1}u^{-1}\\v&\lambda^{-1}u^{-1}+1\end{pmatrix}T_{uv}
$$
and the Poisson brackets become:
$$\{u_j,v_i\}=u_jv_i \mbox{ if } j=i+1 \mbox{ or } j=i-1,$$
$$\{u_j,v_i\}=-2u_jv_i \mbox{ if } j=i$$
which has been already used in sect.~\ref{ss:toda2} (cf. e.g. with \rf{sf4}).

\section*{Acknowledgements}

We would like to thank the QGM, University of Aarhus, the Max Planck Institute, and the Haussdorf Institute for Mathematics in Bonn where essential parts of this work have been done.
The work of V.F. has been partially supported by ANR GTAA and ANR ETTT grants. The work of A.M. has been partially supported by the research grant 13-05-0006 of NRU HSE, by joint RFBR project 12-02-92108, by the Program of Support of Scientific Schools (NSh-3349.2012.2), and by the Russian Ministry of Education under the contract 8207.

\appendix

\section*{Appendix}
\setcounter{equation}0
\renewcommand{\thesubsection}{\Alph{subsection}}
\renewcommand{\theequation}{\Alph{subsection}\arabic{equation}}
\subsection{Cluster varieties of type $\mathcal{X}$\label{ap:cluster}}

Here we collect necessary definitions and properties of cluster varieties of type $\mathcal{X}$. For simplicity of the presentation we restrict ourself to the so-called simply-laced case, where the exchange matrices are skew-symmetric.

A \textit{seed torus} is just a pair $(\mathcal{X},\varepsilon)$ consisting of a split algebraic torus $\mathcal{X}=(\mathbb{C}^\times)^N$ and a skew-symmetric matrix $\varepsilon^{ij}$, $1 \leqslant i,j \leqslant N$ called \textit{exchange matrix}. The standard coordinates on $\mathcal{X}$ are denoted by $\{ x_i\}$.
An isomorphism between two seed tori $(\mathcal{X},\varepsilon)$ and $(\mathcal{X}',\varepsilon')$ is a map preserving the splitting and sending one exchange matrix to another. It is given by a bijection $\sigma:[1,\ldots,N]\to [1,\ldots,N]$ of the sets of coordinates such that $\varepsilon^{ij}=\varepsilon'_{\sigma(i)\sigma(j)}$.

A map between two seeds $(\mathcal{X},\varepsilon)$ and $(\mathcal{X}',\varepsilon')$ is called a \textit{mutation in the $k$-th coordinate} if  $\varepsilon^{ik}\in \mathbb Z$ for any $i$ and the exchange matrices are related by
$$
\varepsilon'_{ij}= \left\{
  \begin{array}{l}
  -\varepsilon^{ij} \mbox{ if } i=k \mbox{ or } j=k\\
  \varepsilon^{ij} \mbox{ if } \varepsilon^{ik}\varepsilon^{kj}<0\\
  \varepsilon^{ij}+ \varepsilon^{ik}|\varepsilon^{kj}| \mbox{ otherwise }
  \end{array}
 \right.
$$
and the coordinates are related by
$$
x'_i = \left\{
  \begin{array}{l}
 1/x_i \mbox{ if } i=k\\
 x_i(1+x_k)^{\varepsilon^{ij}} \mbox{ if } \varepsilon^{ij}\geqslant 0\\
 x_i(1+1/x_k)^{\varepsilon^{ij}} \mbox{ if } \varepsilon^{ij}< 0
\end{array}
\right.
$$
A cluster transformation is a composition of mutations. A cluster map is a composition of cluster transformations and projections along the standard coordinate axes.

A \textit{cluster variety} $X$ of type $\mathcal{X}$ is an algebraic variety covered up to codimension two by (possibly infinite) set of \textit{cluster charts}, which are the maps $\phi^\alpha(\mathcal{X}^\alpha)\to X$, where $\{(\mathcal{X}^\alpha,\varepsilon^\alpha)\}$ is a collection of seeds, and such that any transition function $(\phi^\alpha)^{-1}\phi^\beta$ is a cluster transformation.
A cluster variety possesses a canonical Poisson structure given by
\be
\label{clubra}
 \{x_i,x_j\}=\varepsilon^{ij}x_ix_j
\ee
in any cluster chart $(X^\alpha,\varepsilon^\alpha)$.

A regular function on $X$ being Laurent polynomial with positive integral coefficients in any chart is called a \textit{cluster function}. Cluster function which can not be presented as a sum of two cluster functions is called \textit{indecomposable}. Conjecturally the set of indecomposable cluster functions form a basis called \textit{canonical basis} in the space of regular functions on $X$.

In order to establish a (birational) isomorphism between two cluster varieties $X$ and $Y$ it is enough to establish an isomorphism of a seed of $X$ with a seed of $Y$. In particular, any isomorphism between a seed and another seed, obtained from the first by a cluster transformation, defines an automorphism of the whole cluster variety. The group of such transformations is called the \textit{mapping class group} of the cluster variety $X$.  This term comes from the analogy with Teichm\"uller theory, where this group is indeed the mapping class group of the corresponding surface: one should not understand this term here literally.

\subsection{Relations among the generators of a simply laced Lie group}\label{ap:relx}
\begin{enumerate}
\item $H_i(x)H_j(y)=H_j(y)H_i(x),$
\item $H_i(x)H_i(y)=H_i(xy),$
\item $E_iH_i(-1)E_i=H_i(-1)$
\item $H_i(x)E_j=E_jH_i(x)$ for $i\neq |j|,$
\item $E_i E_j = E_jE_i$ if $C_{ij}=0$,
\item $ E_i H_i(x)E_i
= H_i(1+x) E_i H_i(1+x^{-1})^{-1}, $

\item
$E_i H_i(x) E_j E_i =  H_i(1+x)H_j(1+x^{-1})^{-1} E_j H_j(x)^{-1}
E_i E_j H_i(1+x^{-1})^{-1} H_j(1+x)$\\ for $C_{ij}=-1$ and $i,j>0$,

\item
$E_{\bar{i}} H_i(x) E_{\bar{j}} E_{\bar{i}} =  H_i(1+x^{-1})^{-1}H_j(1+x) E_{\bar{j}} H_j(x)^{-1}
E_{\bar{i}} E_{\bar{i}} H_i(1+x) H_j(1+x^{-1})^{-1}$\\ for $C_{ij}=-1$ and $i,j>0$,

\item $E_{\bar{i}}H_i(x)
E_i = \prod\limits_{j \neq  i} H_j(1+x)^{-C_{ij}}
H_i(1+x^{-1})^{-1} E_i H_i(x^{-1}) E_{\bar{i}} H_i(1+x^{-1})^{-1}$ for $i>0$.

\end{enumerate}

\setcounter{equation}0
\subsection{Exchange graphs and decompositions of $u\in W\times W$ 
\label{ap:exchange}}

As it is shown in \cite{Drinfeld}, the Poisson brackets induced by the Drinfeld-Jimbo $r$-matrix \rf{rDD} on the parameters $x_i$ is log-constant and half-integral, i.e. given by the formula \rf{clubra},
 where $\varepsilon^{ij}$ is a skew-symmetric matrix taking integral or half-integral values. The same matrix plays a role of the exchange matrix (if the group $G$ is simply laced) for the corresponding cluster variety. Now we shall give the description of the matrix $\varepsilon^{ij}$ for the parametrisations of cells, we use in sect.~\ref{ss:clustcell} and sect.~\ref{s:loop}.

\begin{figure}
\begin{center}

\begin{tikzpicture}
\draw[semithick] (0,0) -- (6.5,0);
\draw[semithick] (0,1) -- (6.5,1);
\draw[semithick] (0,2) -- (6.5,2);
\draw[ultra thick,color=gray!55] (1,1) -- (1.5,2);
\draw[ultra thick,color=gray!55,shorten >=0.45cm,-stealth] (1,1) -- (1.5,2);
\draw[ultra thick,color=gray!55] (0.5,2) -- (1,1);
\draw[ultra thick,color=gray!55,shorten >=0.45cm,-stealth] (0.5,2) -- (1,1);
\draw[ultra thick,shorten >=0.4cm,-stealth] (1.5,2) -- (0.5,2);
\draw[ultra thick] (1.5,2) -- (0.5,2);
\fill (1,0) circle (0.08);
\fill (1,1) circle (0.08);
\fill (0.5,2) circle (0.08);
\fill (1.5,2) circle (0.08);
\draw[ultra thick,color=gray!55] (2.5,1) -- (3,2);
\draw[ultra thick,color=gray!55,shorten <=0.45cm,stealth-] (2.5,1) -- (3,2);
\draw[ultra thick,color=gray!55] (2,2) -- (2.5,1);
\draw[ultra thick,color=gray!55,shorten <=0.45cm,stealth-] (2,2) -- (2.5,1);
\draw[ultra thick,shorten <=0.4cm,stealth-] (3,2) -- (2,2);
\draw[ultra thick] (3,2) -- (2,2);
\fill (2.5,0) circle (0.08);
\fill (2.5,1) circle (0.08);
\fill (2,2) circle (0.08);
\fill (3,2) circle (0.08);
\draw[ultra thick,color=gray!55] (3.5,1) -- (4,2);
\draw[ultra thick,color=gray!55,shorten >=0.45cm,-stealth] (3.5,1) -- (4,2);
\draw[ultra thick,color=gray!55] (4,2) -- (4.5,1);
\draw[ultra thick,color=gray!55,shorten >=0.45cm,-stealth] (4,2) -- (4.5,1);
\draw[ultra thick,color=gray!55] (3.5,1) -- (4,0);
\draw[ultra thick,color=gray!55,shorten >=0.45cm,-stealth] (3.5,1) -- (4,0);
\draw[ultra thick,color=gray!55] (4,0) -- (4.5,1);
\draw[ultra thick,color=gray!55,shorten >=0.45cm,-stealth] (4,0) -- (4.5,1);
\draw[ultra thick,shorten <=0.4cm,stealth-] (3.5,1) -- (4.5,1);
\draw[ultra thick] (3.5,1) -- (4.5,1);
\fill (4,0) circle (0.08);
\fill (4,2) circle (0.08);
\fill (3.5,1) circle (0.08);
\fill (4.5,1) circle (0.08);
\draw[ultra thick,color=gray!55] (5.5,1) -- (6,0);
\draw[ultra thick,color=gray!55,shorten <=0.45cm,stealth-] (5.5,1) -- (6,0);
\draw[ultra thick,color=gray!55] (5,0) -- (5.5,1);
\draw[ultra thick,color=gray!55,shorten <=0.45cm,stealth-]  (5,0) -- (5.5,1);
\draw[ultra thick,shorten <=0.4cm,stealth-] (6,0) -- (5,0);
\draw[ultra thick] (6,0) -- (5,0);
\fill (5.5,1) circle (0.08);
\fill (5.5,2) circle (0.08);
\fill (5,0) circle (0.08);
\fill (6,0) circle (0.08);
\draw[->] (7,1) -- (9,1);
\draw[semithick] (9.5,0) -- (12.5,0);
\draw[semithick] (9.5,1) -- (12.5,1);
\draw[semithick] (9.5,2) -- (12.5,2);
\draw[ultra thick,shorten <=0.4cm,stealth-] (10,2) -- (11,2);
\draw[ultra thick] (10,2) -- (11,2);
\draw[ultra thick,shorten >=0.4cm,-stealth] (11,2) -- (12,2);
\draw[ultra thick] (11,2) -- (12,2);
\draw[ultra thick,shorten <=0.4cm,stealth-] (11,2) -- (10.5,1);
\draw[ultra thick] (11,2) -- (10.5,1);
\draw[ultra thick,shorten >=0.4cm,-stealth] (11,0) -- (12,0);
\draw[ultra thick] (11,0) -- (12,0);
\draw[ultra thick,shorten <=0.4cm,stealth-] (10.5,1) -- (11.5,1);
\draw[ultra thick] (10.5,1) -- (11.5,1);
\draw[ultra thick,color=gray!55,shorten >=0.45cm,-stealth] (10,2) -- (10.5,1);
\draw[ultra thick,color=gray!55] (10,2) -- (10.5,1);
\draw[ultra thick,color=gray!55,shorten >=0.45cm,-stealth] (10.5,1) -- (11,0);
\draw[ultra thick,color=gray!55] (10.5,1) -- (11,0);
\draw[ultra thick,color=gray!55,shorten >=0.45cm,-stealth] (12,2) -- (11.5,1);
\draw[ultra thick,color=gray!55] (12,2) -- (11.5,1);
\draw[ultra thick,color=gray!55,shorten >=0.45cm,-stealth] (12,0) -- (11.5,1);
\draw[ultra thick,color=gray!55] (12,0) -- (11.5,1);
\fill (10,2) circle (0.08);
\fill (11,2) circle (0.08);
\fill (12,2) circle (0.08);
\fill (10.5,1) circle (0.08);
\fill (11.5,1) circle (0.08);
\fill (11,0) circle (0.08);
\fill (12,0) circle (0.08);
\end{tikzpicture}
\end{center}
\caption{Example of a graph, made of an ordered sequence of chords (left), for the word $\bar{1}1\bar{2}3$ in the case of group $PGL(4)$. The result of gluing into a single graph is presented on the right picture.}
\label{fi:sl4-simple}
\end{figure}
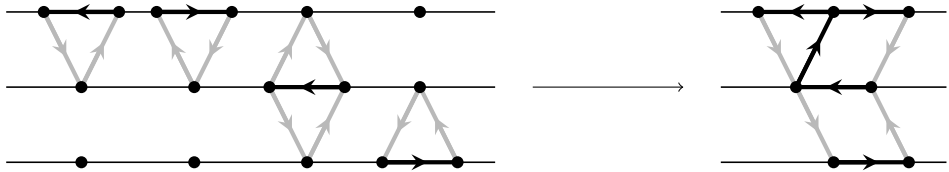

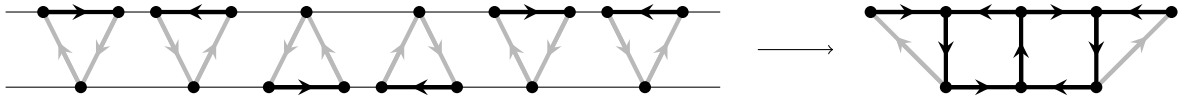
\begin{figure}
\begin{center}
\begin{tikzpicture}
\draw[ultra thick,shorten >=0.4cm,-stealth] (0.5,1) -- (1.5,1);
\draw[ultra thick] (0.5,1) -- (1.5,1);
\draw[ultra thick,color=gray!55,shorten >=0.45cm,-stealth] (1.5,1) -- (1,0);
\draw[ultra thick,color=gray!55] (1.5,1) -- (1,0);
\draw[ultra thick,color=gray!55,shorten >=0.45cm,-stealth] (1,0) -- (0.5,1);
\draw[ultra thick,color=gray!55] (1,0) -- (0.5,1);
\fill (0.5,1) circle (0.08);
\fill (1.5,1) circle (0.08);
\fill (1,0) circle (0.08);
\draw[ultra thick,shorten <=0.4cm,stealth-] (2,1) -- (3,1);
\draw[ultra thick] (2,1) -- (3,1);
\draw[ultra thick,color=gray!55,shorten <=0.45cm,stealth-] (3,1) -- (2.5,0);
\draw[ultra thick,color=gray!55] (3,1) -- (2.5,0);
\draw[ultra thick,color=gray!55,shorten <=0.45cm,stealth-] (2.5,0) -- (2,1);
\draw[ultra thick,color=gray!55] (2.5,0) -- (2,1);
\fill (2,1) circle (0.08);
\fill (3,1) circle (0.08);
\fill (2.5,0) circle (0.08);
\draw[ultra thick,shorten >=0.4cm,-stealth] (3.5,0) -- (4.5,0);
\draw[ultra thick] (3.5,0) -- (4.5,0);
\draw[ultra thick,color=gray!55,shorten >=0.45cm,-stealth] (4.5,0) -- (4,1);
\draw[ultra thick,color=gray!55] (4.5,0) -- (4,1);
\draw[ultra thick,color=gray!55,shorten >=0.45cm,-stealth] (4,1) -- (3.5,0);
\draw[ultra thick,color=gray!55] (4,1) -- (3.5,0);
\fill (3.5,0) circle (0.08);
\fill (4.5,0) circle (0.08);
\fill (4,1) circle (0.08);
\draw[ultra thick,shorten <=0.4cm,stealth-] (5,0) -- (6,0);
\draw[ultra thick] (5,0) -- (6,0);
\draw[ultra thick,color=gray!55,shorten <=0.45cm,stealth-] (6,0) -- (5.5,1);
\draw[ultra thick,color=gray!55] (6,0) -- (5.5,1);
\draw[ultra thick,color=gray!55,shorten <=0.45cm,stealth-] (5.5,1) -- (5,0);
\draw[ultra thick,color=gray!55] (5.5,1) -- (5,0);
\fill (5,0) circle (0.08);
\fill (6,0) circle (0.08);
\fill (5.5,1) circle (0.08);
\draw[ultra thick,shorten >=0.4cm,-stealth] (6.5,1) -- (7.5,1);
\draw[ultra thick] (6.5,1) -- (7.5,1);
\draw[ultra thick,color=gray!55,shorten >=0.45cm,-stealth] (7.5,1) -- (7,0);
\draw[ultra thick,color=gray!55] (7.5,1) -- (7,0);
\draw[ultra thick,color=gray!55,shorten >=0.45cm,-stealth] (7,0) -- (6.5,1);
\draw[ultra thick,color=gray!55] (7,0) -- (6.5,1);
\fill (6.5,1) circle (0.08);
\fill (7.5,1) circle (0.08);
\fill (7,0) circle (0.08);
\draw[ultra thick,shorten <=0.4cm,stealth-] (8,1) -- (9,1);
\draw[ultra thick] (8,1) -- (9,1);
\draw[ultra thick,color=gray!55,shorten <=0.45cm,stealth-] (9,1) -- (8.5,0);
\draw[ultra thick,color=gray!55] (9,1) -- (8.5,0);
\draw[ultra thick,color=gray!55,shorten <=0.45cm,stealth-] (8.5,0) -- (8,1);
\draw[ultra thick,color=gray!55] (8.5,0) -- (8,1);
\fill (8,1) circle (0.08);
\fill (9,1) circle (0.08);
\fill (8.5,0) circle (0.08);
\draw (0,0)--(9.5,0);
\draw (0,1)--(9.5,1);
\draw[->] (10,0.5)-- (11,0.5);
\draw[ultra thick,shorten >=0.4cm,-stealth] (11.5,1) -- (12.5,1);
\draw[ultra thick] (11.5,1) -- (12.5,1);
\draw[ultra thick,shorten <=0.4cm,stealth-] (12.5,1) -- (13.5,1);
\draw[ultra thick] (12.5,1) -- (13.5,1);
\draw[ultra thick,shorten >=0.4cm,-stealth] (13.5,1) -- (14.5,1);
\draw[ultra thick] (13.5,1) -- (14.5,1);
\draw[ultra thick,shorten <=0.4cm,stealth-] (14.5,1) -- (15.5,1);
\draw[ultra thick] (14.5,1) -- (15.5,1);
\draw[ultra thick,shorten >=0.4cm,-stealth] (12.5,0) -- (13.5,0);
\draw[ultra thick] (12.5,0) -- (13.5,0);
\draw[ultra thick,shorten <=0.4cm,stealth-] (13.5,0) -- (14.5,0);
\draw[ultra thick] (13.5,0) -- (14.5,0);
\draw[ultra thick,shorten <=0.4cm,stealth-] (12.5,0) -- (12.5,1);
\draw[ultra thick] (12.5,0) -- (12.5,1);
\draw[ultra thick,shorten >=0.4cm,-stealth] (13.5,0) -- (13.5,1);
\draw[ultra thick] (13.5,0) -- (13.5,1);
\draw[ultra thick,shorten <=0.4cm,stealth-] (14.5,0) -- (14.5,1);
\draw[ultra thick] (14.5,0) -- (14.5,1);
\draw[ultra thick,shorten >=0.46cm,-stealth,color=gray!55] (12.5,0) -- (11.5,1);
\draw[ultra thick,color=gray!55] (12.5,0) -- (11.5,1);
\draw[ultra thick,shorten <=0.46cm,stealth-,color=gray!55] (15.5,1) -- (14.5,0);
\draw[ultra thick,color=gray!55] (15.5,1) -- (14.5,0);
\fill (11.5,1) circle (0.08);
\fill (12.5,1) circle (0.08);
\fill (13.5,1) circle (0.08);
\fill (14.5,1) circle (0.08);
\fill (15.5,1) circle (0.08);
\fill (12.5,0) circle (0.08);
\fill (13.5,0) circle (0.08);
\fill (14.5,0) circle (0.08);
\end{tikzpicture}
\end{center}
\caption{Example of a graph, constructed as on fig.~\ref{fi:sl4-simple} from glued chords, for the word $1\bar{1}2\bar{2}1\bar{1}$ (big cell) in the group $PGL(3)$.}
\label{fi:sl3-simple}
\end{figure}

Instead of writing formulae, we shall give a construction of a graph with oriented edges called \emph{exchange graph}~\footnote{This graph is dual to the bipartite graph $\Gamma$ on torus $\Sigma$ used to construct the dimer partition functions.}  with vertices in bijection with the indices $i$ (or the respective coordinate functions $x_i$) and $\varepsilon^{ij}$ arrows from $i$ to $j$ if $\varepsilon^{ij}\geqslant 0$. If the value of $\varepsilon^{ij}$ is half-integral, we shall indicate  the fractional part by a grey arrow, see fig.~\ref{fi:sl4-simple} and fig.~\ref{fi:sl3-simple}.

Let $u=s_{k_1}\cdots s_{k_l}$ be a word made from the generators of $W\times W$.
To construct the graph start with a staff: a collection of disjoint horizontal lines on the plane, enumerated by positive simple roots $i\in\Pi$. The points on the staff are ordered, like in music, by their projections onto a horizontal axis. For each $k_j$ draw a chord -- a graph $\Gamma(k_j)$ with the vertices on the staff and oriented edges. Different chords are put on the staff respecting the order of the generators $s_{k_j}$ in the word, e.g. the chord $\Gamma(k_{j+1})$ is located to the right from the chord $\Gamma(k_j)$.
A chord $\Gamma(k_j)$ has one leftmost vertex $L_j(k_j)$ and one rightmost vertex $R_j(k_j)$ on the $j$-th line $j\in\Pi$, together with the vertices $\{ S_l(k_j)|l\neq j\}$, located at each other line between $L_j(k_j)$ and $R_j(k_j)$ with arbitrary mutual order. If $k_j>0$ we draw on $j$-th line of the staff between $L_j(k_j)$ and $R_j(k_j)$ a forward oriented edge, and connect $L_j(k_j)$ and $R_j(k_j)$ with $S_l(k_j)$ by $-C_{jl}/2$ backward oriented edges (drawn in gray for most of the cases, indicating that often $-C_{jl}/2$ is equal to $1/2$). For $k_j < 0$ we draw the same chords, but the orientation of all arrows is opposite, see the left pictures at fig.~\ref{fi:sl4-simple} and fig.~\ref{fi:sl3-simple}.

To get the desired graph just contract all staff lines and remove pairs of arrows connecting the same vertices with opposite orientation, as illustrated at the right pictures on figs.~\ref{fi:sl4-simple},\ref{fi:sl3-simple}. The graph from fig.~\ref{fi:sl3-simple} exactly corresponds to the example considered above (see formula \rf{bigcell} in sect.~\ref{ss:clustcell}),
of the big cell in $PGL(3)$.

\begin{figure}[ht]
\begin{center}
\begin{tikzpicture}
 \draw[ultra thick,shorten >=0.4cm,-stealth] (0.5,1) -- (1.5,1);
\draw[ultra thick] (0.5,1) -- (1.5,1);
\draw[ultra thick,color=gray!55,shorten >=0.45cm,-stealth] (1.5,1) -- (1,0);
\draw[ultra thick,color=gray!55] (1.5,1) -- (1,0);
\draw[ultra thick,color=gray!55,shorten >=0.45cm,-stealth] (1,0) -- (0.5,1);
\draw[ultra thick,color=gray!55] (1,0) -- (0.5,1);
\fill (0.5,1) circle (0.08);
\fill (1.5,1) circle (0.08);
\fill (1,0) circle (0.08);
\draw[ultra thick,shorten <=0.4cm,stealth-] (2,1) -- (3,1);
\draw[ultra thick] (2,1) -- (3,1);
\draw[ultra thick,color=gray!55,shorten <=0.45cm,stealth-] (3,1) -- (2.5,0);
\draw[ultra thick,color=gray!55] (3,1) -- (2.5,0);
\draw[ultra thick,color=gray!55,shorten <=0.45cm,stealth-] (2.5,0) -- (2,1);
\draw[ultra thick,color=gray!55] (2.5,0) -- (2,1);
\fill (2,1) circle (0.08);
\fill (3,1) circle (0.08);
\fill (2.5,0) circle (0.08);
\draw[ultra thick,shorten >=0.4cm,-stealth] (3.5,0) -- (4.5,0);
\draw[ultra thick] (3.5,0) -- (4.5,0);
\draw[ultra thick,color=gray!55,shorten >=0.45cm,-stealth] (4.5,0) -- (4,1);
\draw[ultra thick,color=gray!55] (4.5,0) -- (4,1);
\draw[ultra thick,color=gray!55,shorten >=0.45cm,-stealth] (4,1) -- (3.5,0);
\draw[ultra thick,color=gray!55] (4,1) -- (3.5,0);
\fill (3.5,0) circle (0.08);
\fill (4.5,0) circle (0.08);
\fill (4,1) circle (0.08);
\draw[ultra thick,shorten <=0.4cm,stealth-] (5,0) -- (6,0);
\draw[ultra thick] (5,0) -- (6,0);
\draw[ultra thick,color=gray!55,shorten <=0.45cm,stealth-] (6,0) -- (5.5,1);
\draw[ultra thick,color=gray!55] (6,0) -- (5.5,1);
\draw[ultra thick,color=gray!55,shorten <=0.45cm,stealth-] (5.5,1) -- (5,0);
\draw[ultra thick,color=gray!55] (5.5,1) -- (5,0);
\fill (5,0) circle (0.08);
\fill (6,0) circle (0.08);
\fill (5.5,1) circle (0.08);
\draw (0,0)--(6.5,0);
\draw (0,1)--(6.5,1);
\draw[very thick] (0,0)--(0,1);
\draw[very thick](6.5,0)--(6.5,1);
\draw[thin] (0.13,0)--(0.13,1);
\draw[thin](6.37,0)--(6.37,1);
\fill (0.2,0.57) circle (0.03);
\fill (0.2,0.43) circle (0.03);
\fill (6.3,0.57) circle (0.03);
\fill (6.3,0.43) circle (0.03);
\draw[->] (7,0.5)-- (8,0.5);
\draw[very thick, double distance = 1pt] (8.5,0) -- (9.5,0);
\draw[ultra thick,shorten <=0.4cm,stealth-] (9.5,0) -- (9.5,1);
\draw[ultra thick] (9.5,0) -- (9.5,1);
\draw[very thick, double distance = 1pt] (9.5,1) -- (8.5,1);
\draw[ultra thick,shorten <=0.4cm,stealth-] (8.5,1) -- (8.5,0);
\draw[ultra thick] (8.5,1) -- (8.5,0);
\draw[ultra thick,-stealth] (8.9,1) -- (9,1);
\draw[ultra thick,stealth-] (8.9,0) -- (9.0,0);
\fill (8.5,0) circle (0.08);
\fill (8.5,1) circle (0.08);
\fill (9.5,0) circle (0.08);
\fill (9.5,1) circle (0.08);
\end{tikzpicture}
\end{center}
\caption{Example of exchange graph  for the word $1\bar{1}2\bar{2}$ for $PGL(3)/\mbox{Ad}H$ }
\label{fi:todasl3}
\end{figure}
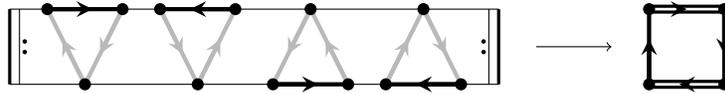
To construct the graphs for the cells of the quotient $G/\mbox{Ad}H$ just consider the staff on a cylinder instead  of a plane, so that the staff lines turn into circles (e.g. the graph for the symplectic leaf in $PGL(3)/\mbox{Ad}H$, corresponding to fig.~\ref{fi:sl3-simple} can be obtained just by identifying the leftmost point for every staff line  with the rightmost one). Similarly, the graph for the symplectic leave of the $PGL(3)$ Toda chain (cf. with \rf{pretoda}) is constructed exactly in the same way, as shown on fig.~\ref{fi:todasl3}.

\setcounter{equation}0
\subsection{Thurston diagrams
\label{ap:thurston}}
Here we briefly introduce and discuss a combinatorial object, introduced by D.~Thurston in \cite{Thurston}, and show its relations to triangulations, bipartite planar graphs, decompositions of permutations into product of generators, coverings of surfaces and link diagrams.

\begin{figure}[ht]
\begin{center}
\begin{tikzpicture}
\fill[color=gray!10] (0,0) circle (1.5);
\draw[thick]         (0,0) circle (1.5);
\fill[color=gray!0] (0,1)..controls(-0.4,0.5)..(0,0)..controls(0.4,0.5)..(0,1);
\draw                (0,1)..controls(-0.4,0.5)..(0,0)..controls(0.4,0.5)..(0,1);
\fill[color=gray!0] (0,1)--(107:1.5) arc (107:138:1.5)--(0,1);
\draw 		     (0,1)--(107:1.5) arc (107:138:1.5)--(0,1);
\fill[color=gray!0] (0,1)--(73:1.5) arc (73:30:1.5)..controls (0.5,1) .. (0,1);%
\draw                (0,1)--(73:1.5) arc (73:30:1.5)..controls (0.5,1) .. (0,1);%
\fill[color=gray!0] (-0.4,-0.5)--(0.4,-0.5)--(-80:1.5) arc (-80:-100:1.5)--cycle;
\draw                (-0.4,-0.5)--(0.4,-0.5)--(-80:1.5) arc (-80:-100:1.5)--cycle;
\fill[color=gray!0] (0,0).. controls (0.5,0)..(0.4,-0.5)--cycle;
\draw                (0,0).. controls (0.5,0)..(0.4,-0.5)--cycle;
\fill[color=gray!0] (0.4,-0.5)..controls (0.8,-0.5)..(-9:1.5) arc (-9:-51:1.5)--cycle;
\draw                (0.4,-0.5)..controls (0.8,-0.5)..(-9:1.5) arc (-9:-51:1.5)--cycle;
\fill[color=gray!0] (-0.4,-0.5)..controls (-1,-0.5).. (-153:1.5) arc (-153:-129:1.5)--cycle;
\draw                (-0.4,-0.5)..controls (-1,-0.5).. (-153:1.5) arc (-153:-129:1.5)--cycle;
\fill[color=gray!0]  (0,0) .. controls (-0.5,0)..(155:1.5) arc (155:180:1.5).. controls (-0.5,-0.3)..(-0.4,-0.5)--cycle;
\draw                 (0,0) .. controls (-0.5,0)..(155:1.5) arc (155:180:1.5).. controls (-0.5,-0.3)..(-0.4,-0.5)--cycle;
\node at (0,-2) {A};
\end{tikzpicture}\qquad\qquad
\begin{tikzpicture}
\fill[color=gray!10] (0,0) circle (1.5);
\draw[thick]         (0,0) circle (1.5);
\fill[color=gray!0] (0,1)..controls(-0.4,0.5)..(0,0)..controls(0.4,0.5)..(0,1);
\draw                (0,1)..controls(-0.4,0.5)..(0,0)..controls(0.4,0.5)..(0,1);
\fill[color=gray!0] (0,1)--(107:1.5) arc (107:138:1.5)--(0,1);
\draw 		     (0,1)--(107:1.5) arc (107:138:1.5)--(0,1);
\fill[color=gray!0] (0,1)--(73:1.5) arc (73:30:1.5)..controls (0.5,1) .. (0,1);%
\draw                (0,1)--(73:1.5) arc (73:30:1.5)..controls (0.5,1) .. (0,1);%
\fill[color=gray!0] (-0.4,-0.5)--(0.4,-0.5)--(-80:1.5) arc (-80:-100:1.5)--cycle;
\draw                (-0.4,-0.5)--(0.4,-0.5)--(-80:1.5) arc (-80:-100:1.5)--cycle;
\fill[color=gray!0] (0,0).. controls (0.5,0)..(0.4,-0.5)--cycle;
\draw                (0,0).. controls (0.5,0)..(0.4,-0.5)--cycle;
\fill[color=gray!0] (0.4,-0.5)..controls (0.8,-0.5)..(-9:1.5) arc (-9:-51:1.5)--cycle;
\draw                (0.4,-0.5)..controls (0.8,-0.5)..(-9:1.5) arc (-9:-51:1.5)--cycle;
\fill[color=gray!0] (-0.4,-0.5)..controls (-1,-0.5).. (-153:1.5) arc (-153:-129:1.5)--cycle;
\draw                (-0.4,-0.5)..controls (-1,-0.5).. (-153:1.5) arc (-153:-129:1.5)--cycle;
\fill[color=gray!0]  (0,0) .. controls (-0.5,0)..(155:1.5) arc (155:180:1.5).. controls (-0.5,-0.3)..(-0.4,-0.5)--cycle;
\draw                 (0,0) .. controls (-0.5,0)..(155:1.5) arc (155:180:1.5).. controls (-0.5,-0.3)..(-0.4,-0.5)--cycle;
\draw[-stealth] (-0.05,0.2)--(-0.3,-0.1);
\draw[-stealth] (-0.2,-0.15)--(0.2,-0.15);
\draw[-stealth] (0.3,-0.1) -- (0.05,0.2);
\draw[stealth-] (-0.05,0.8)--(-0.3,1.1);
\draw[stealth-] (-0.25,1.15)--(0.25,1.15);
\draw[stealth-] (0.3,1.1) -- (0.05,0.8);
\draw[-stealth] (-0.1,-0.6)--(-0.3,-0.3);
\draw[-stealth] (-0.65,-0.65)--(-0.2,-0.65);
\draw[-stealth] (-0.4,-0.3) -- (-0.7,-0.6);
\draw[-stealth] (0.65,-0.6)--(0.4,-0.3);
\draw[-stealth] (0.15,-0.65)--(0.6,-0.65);
\draw[-stealth] (0.3,-0.3) -- (0.1,-0.6);
\draw[stealth-,color=gray!50] (-85:1.4)--(-45:1.4);
\draw[stealth-,color=gray!50] (-17:1.4)--(45:1.4);
\draw[stealth-,color=gray!50] (70:1.4)--(110:1.4);
\draw[stealth-,color=gray!50] (129:1.4)--(162:1.4);
\draw[stealth-,color=gray!50] (177:1.4)--(212:1.4);
\draw[stealth-,color=gray!50] (225:1.4)--(267:1.4);
\node at (0,0.5) {\tiny{1}};
\node at (0.3,-0.2) {\tiny{2}};
\node at (-0.5,1.2) {\tiny{3}};
\node at (0.7,1.1) {\tiny{4}};
\node at (-0.8,0) {\tiny{5}};
\node at (-0.9,-0.7) {\tiny{6}};
\node at (0,-1) {\tiny{7}};
\node at (1,-0.7) {\tiny{8}};
;\node at (0,-2) {B};
\end{tikzpicture}\qquad\qquad
\begin{tikzpicture}
\fill[color=gray!5] (0,0) circle (1.5);
\fill[color=gray!0] (0,1)..controls(-0.4,0.5)..(0,0)..controls(0.4,0.5)..(0,1);
\draw[color=gray!10](0,1)..controls(-0.4,0.5)..(0,0)..controls(0.4,0.5)..(0,1);
\fill[color=gray!0] (0,1)--(107:1.5) arc (107:138:1.5)--(0,1);
\draw[color=gray!10](0,1)--(107:1.5) arc (107:138:1.5)--(0,1);
\fill[color=gray!0] (0,1)--(73:1.5) arc (73:30:1.5)..controls (0.5,1) .. (0,1);
\draw[color=gray!10](0,1)--(73:1.5) arc (73:30:1.5)..controls (0.5,1) .. (0,1);
\fill[color=gray!0] (-0.4,-0.5)--(0.4,-0.5)--(-80:1.5) arc (-80:-100:1.5)--cycle;
\draw[color=gray!10](-0.4,-0.5)--(0.4,-0.5)--(-80:1.5) arc (-80:-100:1.5)--cycle;
\fill[color=gray!0] (0,0).. controls (0.5,0)..(0.4,-0.5)--cycle;
\draw[color=gray!10](0,0).. controls (0.5,0)..(0.4,-0.5)--cycle;
\fill[color=gray!0]  (0.4,-0.5)..controls (0.8,-0.5)..(-9:1.5) arc (-9:-51:1.5)--cycle;
\draw [color=gray!10](0.4,-0.5)..controls (0.8,-0.5)..(-9:1.5) arc (-9:-51:1.5)--cycle;
\fill[color=gray!0] (-0.4,-0.5)..controls (-1,-0.5).. (-153:1.5) arc (-153:-129:1.5)--cycle;
\draw[color=gray!10](-0.4,-0.5)..controls (-1,-0.5).. (-153:1.5) arc (-153:-129:1.5)--cycle;
\fill[color=gray!0]  (0,0) .. controls (-0.5,0)..(155:1.5) arc (155:180:1.5).. controls (-0.5,-0.3)..(-0.4,-0.5)--cycle;
\draw[color=gray!10](0,0) .. controls (-0.5,0)..(155:1.5) arc (155:180:1.5).. controls (-0.5,-0.3)..(-0.4,-0.5)--cycle;
\draw[thick,color=gray!10]         (0,0) circle (1.5);
\draw (90:1)--(90:1.5);
\draw (0,1)--(10:1); \draw (10:1)--(10:1.5);
\draw (0,1)--(137:1); \draw (137:1)--(147:1.5);
\draw (0,0)--(137:1);
\draw (0,0)--(10:1);
\draw (0,0)--(0,-0.25);
\draw (0,-0.25)--(-0.4,-0.5);
\draw (0,-0.25)--(0.4,-0.5);
\draw (0,-0.25) circle (0.7mm);
\draw (-0.4,-0.5)--(199:1); \draw (199:1)--(195:1.5);
\draw (-0.4,-0.5)--(240:1.1); \draw (240:1.1)--(245:1.5);
\draw (0.4,-0.5)--(300:1.1); \draw (300:1.1)--(295:1.5);
\draw (0.4,-0.5)--(10:1);
\fill (90:1.25) circle (0.07cm);\fill[color=white] (90:1.25) circle (0.05cm);
\fill (10:1) circle (0.07cm);\fill[color=white] (10:1) circle (0.05cm);
\fill (137:1) circle (0.07cm);\fill[color=white] (137:1) circle (0.05cm);
\fill[color=white] (0,-0.25) circle (0.5mm);
\fill (199:1) circle (0.07cm);\fill[color=white] (199:1) circle (0.05cm);
\fill (240:1.1) circle (0.07cm);\fill[color=white] (240:1.1) circle (0.05cm);
\fill (300:1.1) circle (0.07cm);\fill[color=white] (300:1.1) circle (0.05cm);
\fill (0,0) circle (0.07cm);
\fill (-0.4,-0.5) circle (0.07cm);
\fill (0.4,-0.5) circle (0.07cm); );
\fill (0,1) circle (0.07cm); m);
\node at (0,-2) {C};
\end{tikzpicture}
\end{center}
\caption{(A): Thurston diagram on a disc; (B): Construction of the exchange matrix; (C): Construction of a bipartite graph.\label{fi:thurston} }
\end{figure}
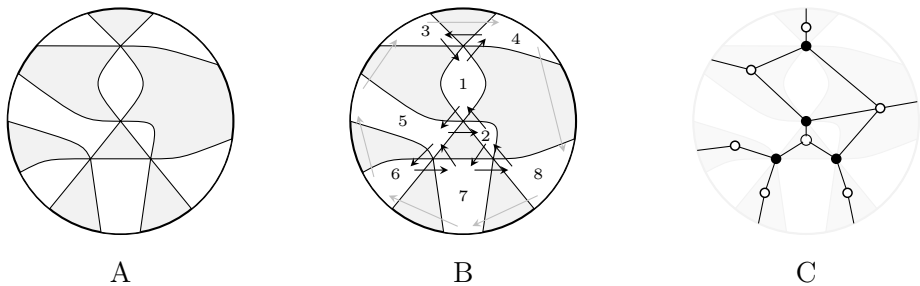

Thurston diagram on an oriented surface $\Sigma$ is an isotopy class of collections of oriented curves, such that all intersection points are triple and the orientation of the curves at every intersection point is alternating. An example of a Thurston diagram on a disc is shown on fig.~\ref{fi:thurston}A. The curves may be either closed or go from boundary to boundary. Connected components of the complement to a diagram are called \textit{faces}. The orientation condition at the intersection points implies that the segments of the curves binding a face are oriented either clockwise (in this case the face is called \textit{white}) or counter-clockwise (in this case the face is called \textit{grey}). All neighbouring faces of a white face are grey and visa versa. On the picture we do not indicate the orientations of the curves, since they can be restored from the colouring of the faces.

\begin{figure}[ht]
\begin{center}
\begin{tikzpicture}
\fill[color=gray!20] (0,0) circle (1);
\fill[color=gray!0] (22:1)..controls(22:0.7)..(0.4,0)..controls(-22:0.7)..(-22:1) arc (-22:22:1);
\draw               (22:1)..controls(22:0.7)..(0.4,0)..controls(-22:0.7)..(-22:1) arc (-22:22:1);
\fill[color=gray!0] (158:1)..controls(158:0.7)..(-0.4,0)..controls(-158:0.7)..(-158:1) arc (-158:-202:1);
\draw               (158:1)..controls(158:0.7)..(-0.4,0)..controls(-158:0.7)..(-158:1) arc (-158:-202:1);
\draw               (67:1)..controls(0.4,0.5)..(0.4,0)..controls(0,0.5)..(-0.4,0)..controls(-0.4,0.5)..(113:1) arc (113:67:1);
\fill[color=gray!0]              (67:1)..controls(0.4,0.5)..(0.4,0)..controls(0,0.5)..(-0.4,0)..controls(-0.4,0.5)..(113:1) arc (113:67:1);
\draw               (67:1)..controls(0.4,0.5)..(0.4,0)..controls(0,0.5)..(-0.4,0)..controls(-0.4,0.5)..(113:1) arc (113:67:1);
\fill[color=gray!0]       (-67:1)..controls(0.4,-0.5)..(0.4,0)..controls(0,-0.5)..(-0.4,0)..controls(-0.4,-0.5)..(-113:1) arc (-113:-67:1);
\draw               (-67:1)..controls(0.4,-0.5)..(0.4,0)..controls(0,-0.5)..(-0.4,0)..controls(-0.4,-0.5)..(-113:1) arc (-113:-67:1);
\draw[thick,color=gray!0] (0,0) circle (1);
\node at (1.8,0) {$\longleftrightarrow$};
\node at (1.8,-1.7) {A};
\end{tikzpicture}
\begin{tikzpicture}
\fill[color=gray!20] (0,0) circle (1);
\fill[color=gray!0] (112:1)..controls(112:0.7)..(0,0.4)..controls(68:0.7)..(68:1) arc (68:112:1);
\draw               (112:1)..controls(112:0.7)..(0,0.4)..controls(68:0.7)..(68:1) arc (68:112:1);
\fill[color=gray!0] (-112:1)..controls(-112:0.7)..(0,-0.4)..controls(-68:0.7)..(-68:1) arc (-68:-112:1);
\draw               (-112:1)..controls(-112:0.7)..(0,-0.4)..controls(-68:0.7)..(-68:1) arc (-68:-112:1);
\fill[color=gray!0] (157:1)..controls(-0.5,0.4)..(0,0.4)..controls(-0.5,0)..(0,-0.4)..controls(-0.5,-0.4)..(203:1) arc (203:157:1);
\draw         (157:1)..controls(-0.5,0.4)..(0,0.4)..controls(-0.5,0)..(0,-0.4)..controls(-0.5,-0.4)..(203:1) arc (203:157:1);
\fill[color=gray!0]              (-22:1)..controls(0.5,-0.4)..(0,-0.4)..controls(0.5,0)..(0,0.4)..controls(0.5,0.4)..(22:1) arc (22:-22:1);
\draw               (-22:1)..controls(0.5,-0.4)..(0,-0.4)..controls(0.5,0)..(0,0.4)..controls(0.5,0.4)..(22:1) arc (22:-22:1);
\draw[thick,color=gray!0] (0,0) circle (1);
\node at (0,-1.9) {~};
\end{tikzpicture}\hspace{2cm}
\begin{tikzpicture}
\fill[color=gray!20] (22:1)..controls(22:0.7)..(0.4,0)..controls(-22:0.7)..(-22:1) arc (-22:22:1);
\draw               (22:1)..controls(22:0.7)..(0.4,0)..controls(-22:0.7)..(-22:1) arc (-22:22:1);
\fill[color=gray!20] (158:1)..controls(158:0.7)..(-0.4,0)..controls(-158:0.7)..(-158:1) arc (-158:-202:1);
\draw               (158:1)..controls(158:0.7)..(-0.4,0)..controls(-158:0.7)..(-158:1) arc (-158:-202:1);
\fill[color=gray!20]              (67:1)..controls(0.4,0.5)..(0.4,0)..controls(0,0.5)..(-0.4,0)..controls(-0.4,0.5)..(113:1) arc (113:67:1);
\draw               (67:1)..controls(0.4,0.5)..(0.4,0)..controls(0,0.5)..(-0.4,0)..controls(-0.4,0.5)..(113:1) arc (113:67:1);
\fill[color=gray!20]       (-67:1)..controls(0.4,-0.5)..(0.4,0)..controls(0,-0.5)..(-0.4,0)..controls(-0.4,-0.5)..(-113:1) arc (-113:-67:1);
\draw               (-67:1)..controls(0.4,-0.5)..(0.4,0)..controls(0,-0.5)..(-0.4,0)..controls(-0.4,-0.5)..(-113:1) arc (-113:-67:1);
\draw[thick,color=gray!0] (0,0) circle (1);
\node at (1.8,0) {$\longleftrightarrow$};
\node at (1.8,-1.7) {B};
\end{tikzpicture}
\begin{tikzpicture}
\fill[color=gray!20] (112:1)..controls(112:0.7)..(0,0.4)..controls(68:0.7)..(68:1) arc (68:112:1);
\draw               (112:1)..controls(112:0.7)..(0,0.4)..controls(68:0.7)..(68:1) arc (68:112:1);
\fill[color=gray!20] (-112:1)..controls(-112:0.7)..(0,-0.4)..controls(-68:0.7)..(-68:1) arc (-68:-112:1);
\draw               (-112:1)..controls(-112:0.7)..(0,-0.4)..controls(-68:0.7)..(-68:1) arc (-68:-112:1);
\fill[color=gray!20] (157:1)..controls(-0.5,0.4)..(0,0.4)..controls(-0.5,0)..(0,-0.4)..controls(-0.5,-0.4)..(203:1) arc (203:157:1);
\draw         (157:1)..controls(-0.5,0.4)..(0,0.4)..controls(-0.5,0)..(0,-0.4)..controls(-0.5,-0.4)..(203:1) arc (203:157:1);
\fill[color=gray!20]              (-22:1)..controls(0.5,-0.4)..(0,-0.4)..controls(0.5,0)..(0,0.4)..controls(0.5,0.4)..(22:1) arc (22:-22:1);
\draw               (-22:1)..controls(0.5,-0.4)..(0,-0.4)..controls(0.5,0)..(0,0.4)..controls(0.5,0.4)..(22:1) arc (22:-22:1);
\draw[thick,color=gray!0] (0,0) circle (1);
\node at (0,-1.9) {~};
\end{tikzpicture}
\end{center}
\begin{center}
\begin{tikzpicture}
\fill[color=gray!20] (0,0)..controls(-1,1)and(1,1)..(0,0)--cycle;
\draw (0,0)..controls(-1,1)and(1,1)..(0,0);
\fill[color=gray!20] (210:1)..controls(225:0.5)..(0,0)--(180:1)  arc (180:210:1)--cycle;
\draw (210:1)..controls(225:0.5)..(0,0)--(180:1);
\fill[color=gray!20] (-30:1)..controls(-45:0.5)..(0,0)--(0:1)  arc (0:-30:1)--cycle;
\draw (-30:1)..controls(-45:0.5)..(0,0)--(0:1);
\draw[thick,color=gray!0] (0,0) circle (1);
\node at (1.8,0) {$\longleftrightarrow$};
\node at (1.8,-1.7) {C};
\end{tikzpicture}
\begin{tikzpicture}
\fill[color=gray!20] (-1,0) arc (180:210:1)--(-30:1) arc (-30:0:1)--cycle;
\draw (-30:1)--(210:1);
\draw (0:1)--(180:1);
draw[thick,color=gray!0] (0,0) circle (1);
\node at (0,-1.9) {~};
\end{tikzpicture}\hspace{2cm}
\begin{tikzpicture}
\fill[color=gray!20] (-1,0)--(0,0)..controls(-1,1)and(1,1)..(0,0)--(1,0) arc (0:180:1)--cycle;
\draw (-1,0)--(0,0)..controls(-1,1)and(1,1)..(0,0)--(1,0);
\fill[color=gray!20] (0,0)..controls(225:0.5)..(210:1) arc (210:330:1)..controls(315:0.5)..(0,0)--cycle;
\draw (0,0)..controls(225:0.5)..(210:1) arc (210:330:1)..controls(315:0.5)..(0,0);
\draw[thick,color=gray!0] (0,0) circle (1);
\node at (1.8,0) {$\longleftrightarrow$};
\node at (1.8,-1.7) {D};
\end{tikzpicture}
\begin{tikzpicture}
\fill[color=gray!20] (-1,0) arc (180:0:1)--cycle;
\draw (-1,0)--(1,0);
\fill[color=gray!20] (210:1) arc (210:330:1)--cycle;
\draw (210:1)--(330:1);
draw[thick,color=gray!0] (0,0) circle (1);
\node at (0,-1.9) {~};
\end{tikzpicture}
\end{center}
\caption{(A): Grey Thurston move; (B): White Thurston move; (C): Grey Thurston reduction; (D): White Thurston reduction.\label{fi:mutation}}
\end{figure}
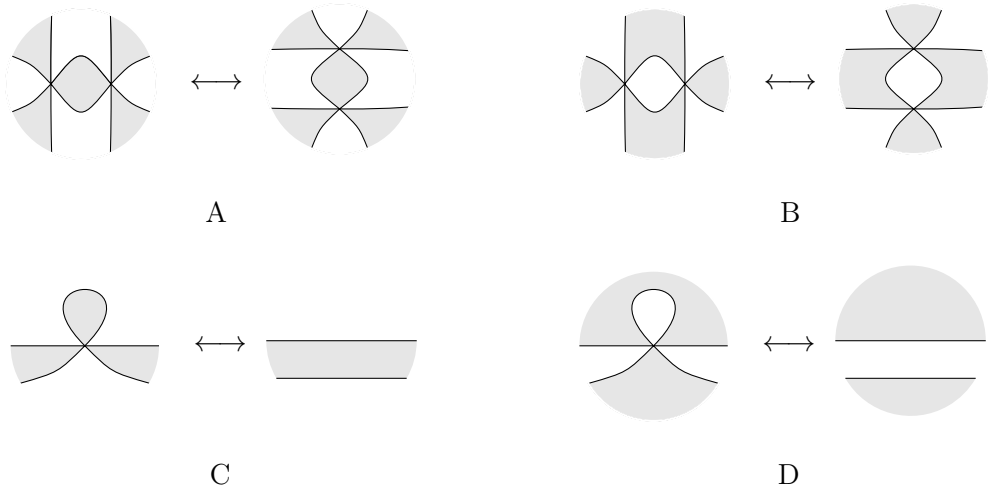

Thurston diagrams admit four kinds of standard modifications, shown on fig.~\ref{fi:mutation}, ABCD. The first two, called Thurston moves, do not change the number of faces. The second two reduce the number of faces and are called Thurston reductions.
A move can be performed each time, when there is a face with one or two corners, and diagrams related by Thurston moves are called equivalent. Diagrams not equivalent to those, where a Thurston reduction can be applied, are called minimal.

\paragraph{Thurston diagrams and cluster varieties.}
Following D.~Thurston and A.~Henriques \cite{Thurston} we define a cluster variety starting from an equivalence class of Thurston diagrams. The charts
\begin{wrapfigure}{l}{0.25\textwidth}
\begin{tikzpicture}
\fill[color=gray!20] (0:1) arc (0:60:1)--(0,0)--cycle;
\fill[color=gray!20] (120:1) arc (120:180:1)--(0,0)--cycle;
\fill[color=gray!20] (240:1) arc (240:300:1)--(0,0)--cycle;
\draw (180:1)--(0:1);
\draw (300:1)--(120:1);
\draw (60:1)--(240:1);
\draw[-stealth] (100:0.8)--(200:0.8);
\draw[-stealth] (220:0.8)--(320:0.8);
\draw[-stealth] (340:0.8)--(80:0.8);
\node at (0,-1.5) {A};
\end{tikzpicture}\qquad
\begin{tikzpicture}
\fill[color=gray!20] (159:3) arc (159:166:3)--(166:2) arc (166:159:2)--cycle;
\draw (166:2)--(166:3);
\fill[color=gray!20] (180:3) arc (180:194:3)--(194:2) arc (194:180:2)--cycle;
\draw[thick] (201:3) arc (201:159:3);
\draw (180:2)--(180:3);
\draw (194:2)--(194:3);
\draw[color=gray!80,-stealth] (198:2.8)--(176:2.8);
\draw[color=gray!80] (170:2.8)--(159:2.8);
\node at (-2.5,-1.5) {B};
\end{tikzpicture}
\caption{Exchange matrix from Thurston diagram.\label{fi:brackets}}
\end{wrapfigure}
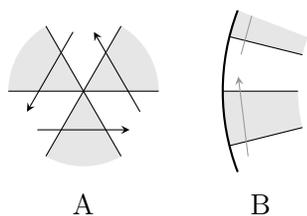
 of this manifold correspond to particular diagrams in a given class and Thurston moves correspond to transition functions called mutation in the cluster language. Cluster variables parametrising a chart are assigned to the white faces. For surfaces with boundary the faces neighbouring the boundary correspond to frozen variables.

The exchange matrix $\varepsilon^{ij}$ is defined as follows.
\begin{itemize}
  \item Draw three arrows connecting white faces around every triple intersection and directed counter-clockwise as shown on fig.~\ref{fi:brackets}A.
  \item For every two segments of the boundary belonging to white faces and separated by a segment, belonging to a grey face, connect them by a \emph{grey} arrow, pointing to the right if viewed from inside the surface (as shown on the fig.~\ref{fi:brackets}B).
\end{itemize}

Then the value of $\varepsilon^{ij}$ is equal to the number of arrows from the white face $i$ to the white face $j$ minus the number of arrows in the inverse direction. The grey arrows are counted with the coefficient one-half. One can easily check, that the exchange matrix does not change under the mutation A, and it changes according to the cluster rule under the mutation B. The example of such collection of arrows together with the indices, enumerating white faces are shown on fig.~\ref{fi:thurston}B. From this picture one can see, for example, that $\varepsilon^{12}=-1$ (one solid arrow going from the face 2 to the face 1) and $\varepsilon^{34}=-1/2$ (one solid arrow from 4 to three and one grey one in the backward direction).

\paragraph{Thurston diagrams and bipartite graphs.}
For every Thurston diagram one can associate a bipartite graph. Conversely, every bipartite graph with three-valent white vertices corresponds to a Thurston diagram.

To construct a bipartite graph out of a Thurston diagram just put a white vertex inside every grey face and a terminal vertex at every grey segment of the boundary. Then, put a black vertex at every triple intersection point and connect it to the three black vertices in the three faces touching the vertex. Connect also the terminal vertices with the corresponding white ones. An example of a graph corresponding to the Thurston diagram is shown on fig.~\ref{fi:thurston}C.

Observe, that this correspondence has the following properties:
\begin{itemize}
 \item A grey Thurston move corresponds to two GK moves of type A (fig.~\ref{fi:spider}).
 \item A white Thurston move corresponds to the GK spider move (type C, fig.~\ref{fi:spider}).
 \item Faces of the bipartite graph correspond to white faces of the Thurston diagram.
 \item The zig-zag paths of the bipartite graph correspond to curves of the Thurston diagram.
\end{itemize}

\paragraph{Thurston diagrams and triangulations.}
For any triangulation of a surface and for an integer $N\geqslant 2$ one can associate a Thurston diagram. We will illustrate this construction for $N=2$ and $N=3$.

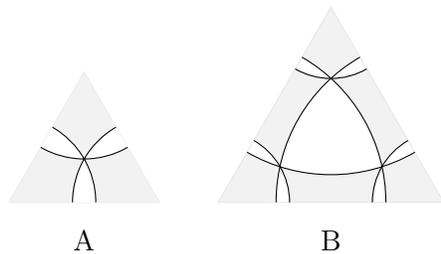
\begin{figure}[H]
\begin{center}
\begin{tikzpicture}
\draw[color=gray!20] (-1,0)--(0,1.732)--(1,0)--cycle;
\fill[color=gray!10] (-1,0)--(0,1.732)--(1,0)--cycle;
\fill[color=gray!0]  (-0.155,0) arc (180:150:1.154) arc (30:0:1.154)--cycle;
\draw (-0.155,0) arc (180:150:1.154) arc (30:0:1.154);
\fill[color=gray!0] (-0.577,0.732) arc (240:270:1.154) arc (30:60:1.154)--cycle;
\draw (-0.577,0.732) arc (240:270:1.154) arc (30:60:1.154);
\fill[color=gray!0] (0.577,0.732) arc (-60:-90:1.154) arc (150:120:1.154)--cycle;
\draw (0.577,0.732) arc (-60:-90:1.154) arc (150:120:1.154);
\node at (0,-0.5) {A};
\end{tikzpicture}\qquad
\begin{tikzpicture}
\fill[color=gray!10] (-1.5,0)+(0:2.23) arc (0:12:2.23) arc (-72:-60:2.23)--(1.5,0)--cycle;
\fill[color=gray!10] (0,2.598)+(240:2.23) arc (240:252:2.23) arc (168:180:2.23)--(-1.5,0)--cycle;
\fill[color=gray!10] (-1.5,0)+(60:2.23) arc (60:48:2.23) arc (132:120:2.23)--(0,2.598)--cycle;
\fill[color=gray!10] (-1.5,0)+(0:0.95) arc (0:30:0.95) arc (252:288:2.23) arc (150:180:0.95)--cycle;
\fill[color=gray!10] (-1.5,0)+(60:0.95) arc (60:30:0.95) arc (168:132:2.23) arc (270:240:0.95)--cycle;
\fill[color=gray!10] (1.5,0)+(120:0.95) arc (120:150:0.95) arc (12:48:2.23) arc (270:300:0.95)--cycle;
\draw[color=gray!20] (-1.5,0)--(0,2.598)--(1.5,0)--cycle;
\draw (-1.5,0)+(0:2.23) arc (0:60:2.23);
\draw (-1.5,0)+(0:0.95) arc (0:60:0.95);
\draw (1.5,0)+(180:2.23) arc (180:120:2.23);
\draw (1.5,0)+(180:0.95) arc (180:120:0.95);
\draw (0,2.598)+(240:2.23) arc (240:300:2.23);
\draw (0,2.598)+(240:0.95) arc (240:300:0.95);
\node at (0,-0.5) {B};
\end{tikzpicture}

\caption{Thurston diagrams, corresponding to a triangle for the cases:
(A) $N=2$ or $SL(2)$; (B) $N=3$ or $SL(3)$.}
\label{fi:thurtriangle}
\end{center}
\end{figure}

Let surface $\Sigma$ be triangulated with edges either entirely belonging to the boundary or to the interior. Replace every triangle by $3(N-1)$ of curves as shown of fig.~\ref{fi:thurtriangle}A for $N=2$, fig.~\ref{fi:thurtriangle}B for $N=3$, and analogously for larger values of $N$.

The constructed correspondence has the following properties:
\begin{itemize}
\item A flip of a triangulation corresponds to a sequence of Thurston moves. For $N=2$ every flip corresponds to a single Thurston move.
\item The corresponding cluster variety is the (framed) space of $SL(N)$ local systems on $\Sigma$ \cite{Higher}.
\item Every closed curve of the Thurston diagram contracts to a curve surrounding one puncture.
\item Every puncture is surrounded by exactly $N-1$ curves.
\item The dual surface $\tilde{\Sigma}$ is an $N$-fold cover of $\Sigma$ ramified in the triple intersection points.
\end{itemize}

\paragraph{Triality of Thurston diagrams.}
Thurston diagrams come in triples such that each one of the triple defines two others. Given a Thurston diagram glue a disc by its boundary to every closed curve and by a half of its boundary for nonclosed ones. We get a $CW$-complex with every 1-cell belonging to three 2-cells --- white, grey and the new discs which we will paint in yellow. Removing all discs of one type one gets a smooth 2-dimensional surface with a Thurston diagram on it. It is easy to see that the yellow-grey surface with reversed orientation is the dual to the initial (white-grey) one.

\paragraph{Thurston diagrams and double permutation group.}

The correspondence between the Thurston diagrams and the words in double permutation group has the following properties:
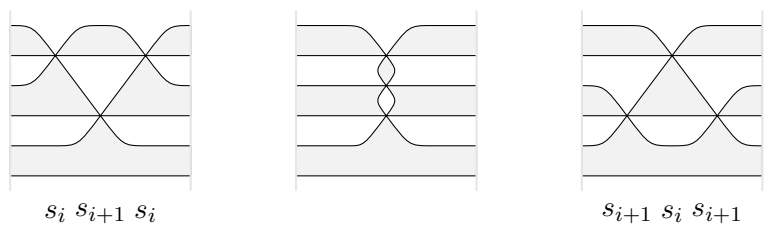
\begin{figure}[ht]
\begin{center}
\begin{tikzpicture}
\fill[color=gray!10] (-0.3,2)..controls(0,2)..(0.3,1.6)--(-0.3,1.6)--cycle;
\draw                (-0.3,2)..controls(0,2)..(0.3,1.6)--(-0.3,1.6)--cycle;
\fill[color=gray!10] (2.1,2)..controls(1.8,2)..(1.5,1.6)--(2.1,1.6)--cycle;
\draw                (2.1,2)..controls(1.8,2)..(1.5,1.6)--(2.1,1.6)--cycle;
\fill[color=gray!10] (-0.3,0.8)--(0.9,0.8)--(0.3,1.6)..controls(0,1.2)..(-0.3,1.2)--cycle;
\draw                (-0.3,0.8)--(0.9,0.8)--(0.3,1.6)..controls(0,1.2)..(-0.3,1.2)--cycle;
\fill[color=gray!10] (2.1,0.8)--(0.9,0.8)--(1.5,1.6)..controls(1.8,1.2)..(2.1,1.2)--cycle;
\draw                (2.1,0.8)--(0.9,0.8)--(1.5,1.6)..controls(1.8,1.2)..(2.1,1.2)--cycle;
\fill[color=gray!10] (-0.3,0)--(2.1,0)--(2.1,0.4)--(1.5,0.4)..controls(1.2,0.4)..(0.9,0.8)..controls(0.6,0.4)..(0.3,0.4)--(-0.3,0.4)--cycle;
\draw (-0.3,0)--(2.1,0)--(2.1,0.4)--(1.5,0.4)..controls(1.2,0.4)..(0.9,0.8)..controls(0.6,0.4)..(0.3,0.4)--(-0.3,0.4)--cycle;
\fill[color=gray!10] (0.3,1.6)--(1.5,1.6)..controls(1.2,2)..(0.9,2)..controls(0.6,2)..(0.3,1.6);
\draw                (0.3,1.6)--(1.5,1.6)..controls(1.2,2)..(0.9,2)..controls(0.6,2)..(0.3,1.6);
\draw[thick,color=gray!20] (2.1,-0.2)--(2.1,2.2);
\draw[thick,color=gray!20] (-0.3,-0.2)--(-0.3,2.2);
\node at (0.3,-0.5) {$s_i$};
\node at (0.9,-0.5) {$s_{i+1}$};
\node at (1.5,-0.5) {$s_i$};
\fill[color=gray!10] (3.5,2)--(4.1,2)..controls(4.4,2)..(4.7,1.6)--(3.5,1.6)--cycle;
\draw                (3.5,2)--(4.1,2)..controls(4.4,2)..(4.7,1.6)--(3.5,1.6)--cycle;
\fill[color=gray!10] (5.9,2)--(5.3,2)..controls(5,2)..(4.7,1.6)--(5.9,1.6)--cycle;
\draw                (5.9,2)--(5.3,2)..controls(5,2)..(4.7,1.6)--(5.9,1.6)--cycle;
\fill[color=gray!10] (4.7,1.6)..controls(4.55,1.4)..(4.7,1.2)..controls(4.85,1.4)..(4.7,1.6);
\draw                (4.7,1.6)..controls(4.55,1.4)..(4.7,1.2)..controls(4.85,1.4)..(4.7,1.6);
\fill[color=gray!10] (3.5,1.2)--(4.7,1.2)..controls(4.55,1)..(4.7,0.8)--(3.5,0.8)--cycle;
\draw                (3.5,1.2)--(4.7,1.2)..controls(4.55,1)..(4.7,0.8)--(3.5,0.8)--cycle;
\fill[color=gray!10] (5.9,1.2)--(4.7,1.2)..controls(4.85,1)..(4.7,0.8)--(5.9,0.8)--cycle;
\draw                (5.9,1.2)--(4.7,1.2)..controls(4.85,1)..(4.7,0.8)--(5.9,0.8)--cycle;
\fill[color=gray!10] (3.5,0)--(5.9,0)--(5.9,0.4)--(5.3,0.4)..controls(5,0.4)..(4.7,0.8)..controls(4.4,0.4)..(4.2,0.4)--(3.5,0.4)--cycle;
\draw (3.5,0)--(5.9,0)--(5.9,0.4)--(5.3,0.4)..controls(5,0.4)..(4.7,0.8)..controls(4.4,0.4)..(4.2,0.4)--(3.5,0.4)--cycle;
\draw[thick,color=gray!20] (3.5,-0.2)--(3.5,2.2);
\draw[thick,color=gray!20] (5.9,-0.2)--(5.9,2.2);
\fill[color=gray!10] (7.3,2)--(7.9,2)..controls(8.2,2)..(8.5,1.6)--(7.3,1.6)--cycle;
\draw                (7.3,2)--(7.9,2)..controls(8.2,2)..(8.5,1.6)--(7.3,1.6)--cycle;
\fill[color=gray!10] (9.7,2)--(9.1,2)..controls(8.8,2)..(8.5,1.6)--(9.7,1.6)--cycle;
\draw                (9.7,2)--(9.1,2)..controls(8.8,2)..(8.5,1.6)--(9.7,1.6)--cycle;
\fill[color=gray!10] (7.3,1.2)..controls(7.6,1.2)..(7.9,0.8)--(7.3,0.8)--cycle;
\draw                (7.3,1.2)..controls(7.6,1.2)..(7.9,0.8)--(7.3,0.8)--cycle;
\fill[color=gray!10] (9.7,1.2)..controls(9.4,1.2)..(9.1,0.8)--(9.7,0.8)--cycle;
\draw                (9.7,1.2)..controls(9.4,1.2)..(9.1,0.8)--(9.7,0.8)--cycle;
\fill[color=gray!10] (7.9,0.8)--(9.1,0.8)--(8.5,1.6)--cycle;
\draw (7.9,0.8)--(9.1,0.8)--(8.5,1.6)--cycle;
\fill[color=gray!10] (7.3,0)--(9.7,0)--(9.7,0.4)..controls(9.4,0.4)..(9.1,0.8)..controls(8.8,0.4)..(8.5,0.4)..controls(8.2,0.4)..(7.9,0.8)..controls(7.6,0.4)..(7.3,0.4)--cycle;
\draw (7.3,0)--(9.7,0)--(9.7,0.4)..controls(9.4,0.4)..(9.1,0.8)..controls(8.8,0.4)..(8.5,0.4)..controls(8.2,0.4)..(7.9,0.8)..controls(7.6,0.4)..(7.3,0.4)--cycle;
\draw[thick,color=gray!20] (7.3,-0.2)--(7.3,2.2);
\draw[thick,color=gray!20] (9.7,-0.2)--(9.7,2.2);
\node at (7.9,-0.5) {$s_{i+1}$};
\node at (8.5,-0.5) {$s_i$};
\node at (9.1,-0.5) {$s_{i+1}$};

\end{tikzpicture}
\end{center}
\caption{Thurston moves corresponding to the relation $s_is_{i+1}s_i=s_{i+1}s_is_{i+1}$.\label{fi:121=212}}
\end{figure}
\begin{itemize}
 \item The relations $s_is_{i+1}s_i=s_{i+1}s_is_{i+1}$ as well as $s_{\bar{i}}s_{\overline{i+1}}s_{\bar{i}}=s_{\overline{i+1}}s_{\bar{i}}s_{\overline{i+1}}$ correspond to a composition of two Thurston moves: one grey and one white (see fig. \ref{fi:121=212}).
\item The relation $s_is_{\bar{i}}=s_{\bar{i}}s_i$ corresponds to one white Thurston move.
\item The relations $s_is_{\overline{i-1}}=s_{\overline{i-1}}s_i$ corresponds to one grey Thurston move.
\item All other Weyl group relations correspond just to isotopies of the diagrams.
\item Product of words corresponds to gluing the strips.
\item Irreducible words correspond to minimal diagrams.
\item Cluster seed corresponding to a minimal Thurston diagram for a given word $s_{i_1}\cdots s_{i_n}$ is isomorphic to the cluster seed corresponding to this word
\end{itemize}

\paragraph{Thurston diagrams and the coextended double affine Weyl group.}
For the group $\widehat{PGL}(N)$ the coextended double of the Weyl group $(W\times W)^\sharp$, presented in the subsection \ref{ss:coextended-Weyl}, corresponds to a Thurston diagram on a cylinder which we shall draw horizontally. Every diagram as in the previous case consists of $2N$ curves without vertical tangents going from left to right. The diagrams representing $s_i$ and $s_{\bar{i}}$ for $i\neq 0$ are constructed in the same way as in the finite case. The diagrams for $s_0$ and $s_{\bar{0}}$ correspond to triple intersections of the curves $2N-1,2N,1$ and $2N,1,2$, respectively. The generator $\Lambda$ is represented by the diagram without intersection but connects the point number $i$ on the left of the cylinder to the point number $i-2$ modulo $N$ on the right. The list of properties of this correspondence reproduces that for the finite case.

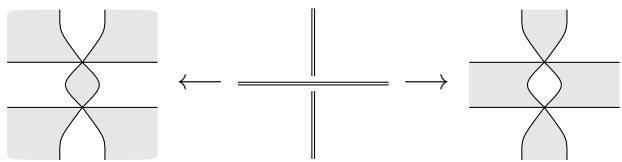
\begin{figure}[ht]
\begin{center}
\begin{tikzpicture}
\fill[color=gray!20] (-0.3,1)--(-0.3,0.9)..controls(-0.3,0.7)..(0,0.3)--(-1,0.3)--(-1,0.9)..controls(-1,1)..(-0.3,1);
\draw (-0.3,1)--(-0.3,0.9)..controls(-0.3,0.7)..(0,0.3)--(-1,0.3);
\fill[color=gray!20] (0.3,1)--(0.3,0.9)..controls(0.3,0.7)..(0,0.3)--(1,0.3)--(1,0.9)..controls(1,1)..(0.3,1);
\draw (0.3,1)--(0.3,0.9)..controls(0.3,0.7)..(0,0.3)--(1,0.3);
\fill[color=gray!20] (-0.3,-1)--(-0.3,-0.9)..controls(-0.3,-0.7)..(0,-0.3)--(-1,-0.3)--(-1,-0.9)..controls(-1,-1)..(-0.3,-1);
\draw (-0.3,-1)--(-0.3,-0.9)..controls(-0.3,-0.7)..(0,-0.3)--(-1,-0.3);
\fill[color=gray!20] (0.3,-1)--(0.3,-0.9)..controls(0.3,-0.7)..(0,-0.3)--(1,-0.3)--(1,-0.9)..controls(1,-1)..(0.3,-1);
\draw (0.3,-1)--(0.3,-0.9)..controls(0.3,-0.7)..(0,-0.3)--(1,-0.3);
\draw[black,fill=gray!20] (0,0.3)..controls(-0.3,0)..(0,-0.3)..controls(0.3,0)..(0,0.3);
\end{tikzpicture}
\begin{tikzpicture}
\draw[double](-1,0)--(1,0);
\draw[double](0,1)--(0,0.1);
\draw[double](0,-0.1)--(0,-1);
\node at (1.5,0) {$\longrightarrow$};
\node at (-1.5,0) {$\longleftarrow$};
\end{tikzpicture}
\begin{tikzpicture}
\draw[black,fill=gray!20] (-0.3,1)--(-0.3,0.9)..controls(-0.3,0.7)..(0,0.3)..controls(0.3,0.7)..(0.3,0.9)--(0.3,1);
\draw[black,fill=gray!20] (-0.3,-1)--(-0.3,-0.9)..controls(-0.3,-0.7)..(0,-0.3)..controls(0.3,-0.7)..(0.3,-0.9)--(0.3,-1);
\draw[black,fill=gray!20] (-1,0.3)--(0,0.3)..controls(-0.3,0)..(0,-0.3)--(-1,-0.3);
\draw[black,fill=gray!20] (1,0.3)--(0,0.3)..controls(0.3,0)..(0,-0.3)--(1,-0.3);
\end{tikzpicture}
\end{center}
\caption{Thurston diagrams corresponding to a link diagram.\label{fi:link}}
\end{figure}

\paragraph{Thurston diagrams and link diagrams.}
For every knot or link diagram on a surface $\Sigma$ one can associate Thurston diagrams according to the rules shown on Figure \ref{fi:link}. It can be done in two ways different by orientation.
The correspondence has the following properties:
\begin{itemize}
\item Changing an over-crossing to an under-crossing on a diagram corresponds to a Thurston move.
\item Thurston diagram of the first type has white faces corresponding to segments between crossings of the link diagram.
\item Thurston diagram of the second type has white faces corresponding to faces and to crossings of the knot diagram.
\item The Reidemeister-III move of the link diagram corresponds to a composition of 4 grey and 4 white Thurston moves. This move corresponds to the move $s_is_{\bar{i}}s_{i+1}s_{\overline{i+1}}s_is_{\bar{i}}\to s_{i+1}s_{\overline{i+1}}s_is_{\bar{i}}s_{i+1}s_{\overline{i+1}}$.
\item The Reidemeister-II move corresponds to one gray and one white Thurston move and two grey and two white Thurston reductions.
\item The Reidemeister-I move corresponds to one grey and one white Thurston reduction.
\item The bipartite graph corresponding to a link diagram on a surface of genus $g$ has $2g-2$ more white than black vertices.

\end{itemize}

An essential part of this paragraph is a translation to the language of Thurston diagrams of the constructions by M.~Cohen, O.T.~Dasbach and H.M.~Russell \cite{CoDaRu}.

\setcounter{equation}0
\subsection{Proofs of the properties of the minors generating functions
\label{ap:super}}
First, let us prove that convolution of generating functions gives a generating function for the product of matrices. Using matrix notation and representing $\bxi$ and $\bxi'$ (respectively $\Beta$ and $\Beta'$) as a row (respectively column) we get

$$\int S(M_1,\bxi',\Beta)S(M_2,\bxi,\Beta')e^{\displaystyle -\bxi\Beta}d\bxi d\Beta=
\int e^{\displaystyle \bxi' M_1\Beta +  \bxi M_2\Beta'-\bxi\Beta}d\bxi d\Beta=
$$
$$=\int e^{\displaystyle -(\bxi-\bxi' M_1)(\Beta-M_2\Beta') + \bxi' M_1M_2\Beta' }d\bxi d\Beta=
$$
$$=e^{\displaystyle \bxi' M_1M_2\Beta' }\int e^{-\displaystyle \bxi\Beta}d\bxi d\Beta= e^{\displaystyle \bxi' M_1M_2\Beta' },
$$
where $d\bxi d\Beta= \prod_{j=1}^nd\xi_jd\eta^j = d\xi_1\cdots d\xi_nd\eta^n\cdots d\eta^1$.

Second, we prove that the expression $e^{\displaystyle \bxi M\Beta}$ is indeed a generating function for the minors of the matrix
$$\exp\left(\sum_{i,j} M^i_j\xi_i\eta^j\right)
=\sum_{k=0}^n\frac{1}{k!}\mathop{\sum_{i_1,\ldots,i_k}}_{j_1,\ldots,j_k}M^{i_1}_{j_1}\cdots M^{i_k}_{j_k}\xi_{i_1}\eta^{j_1}\cdots \xi_{i_k}\eta^{j_k}=$$
$$= \sum_{k=0}^n\mathop{\sum_{i_1<\cdots<i_k}}_{j_1,\ldots,j_k}M^{i_1}_{j_1}\cdots M^{i_k}_{j_k}\xi_{i_1}\eta^{j_1}\cdots \xi_{i_k}\eta^{j_k}=\sum_{k=0}^n(-1)^\frac{k(k-1)}{2}\!\!\!\!\mathop{\sum_{i_1<\cdots<i_k}}_{j_1,\ldots,j_k}M^{i_1}_{j_1}\cdots M^{i_k}_{j_k}\xi_{i_1}\cdots \xi_{i_k}\eta^{j_1}\cdots\eta^{j_k}=
$$
$$=\sum_{k=0}^n(-1)^\frac{k(k-1)}{2}\mathop{\sum_{i_1<\cdots<i_k}}_{j_1<\cdots<j_k}\sum_{\sigma \in \mathcal{S}_k}M^{i_1}_{j_{\sigma(1)}}\cdots M^{i_k}_{j_{\sigma(k)}}\xi_{i_1}\cdots \xi_{i_k}\eta^{j_{\sigma(1)}}\cdots\eta^{j_{\sigma(k)}}=
$$
$$=\sum_{k=0}^n(-1)^\frac{k(k-1)}{2}\mathop{\sum_{i_1<\cdots<i_k}}_{j_1<\cdots<j_k}\sum_{\sigma \in \mathcal{S}_k}(-1)^{\mathop{\mbox{\small sign}}(\sigma)}M^{i_1}_{j_{\sigma(1)}}\cdots M^{i_k}_{j_{\sigma(k)}}\xi_{i_1}\cdots \xi_{i_k}\eta^{j_1}\cdots\eta^{j_k}=
$$
$$=\sum_{k=0}^n(-1)^\frac{k(k-1)}{2}\mathop{\sum_{i_1<\cdots<i_k}}_{j_1<\cdots<j_k}M^{i_1,\ldots,i_k}_{j_1,\ldots,j_k}\xi_{i_1}\cdots \xi_{i_k}\eta^{j_1}\cdots\eta^{j_k}=
$$

$$=\sum_{k=0}^n\mathop{\sum_{i_1<\cdots<i_k}}_{j_1<\cdots<j_k}M^{i_1,\ldots,i_k}_{j_1,\ldots,j_k}\xi_{i_1}\cdots \xi_{i_k}\eta^{j_k}\cdots\eta^{j_1}
$$
exactly as stated in property~\ref{generating} of Lemma~\ref{le:matrices}.

\setcounter{equation}0
\subsection{Schwartz coordinates and the pentagram map
\label{ap:schwartz}}

To simplify notations replace the indices $i+k$ by $k$ primes.
In these new notations the normalization condition (\ref{normalization}) reads as:
\begin{equation}
e'''=a(e'+e)+be''
\end{equation}
Observe first the following identity:
\begin{equation}
a(e'+e)-be'=\frac{a'+bb'}{a'}e''' - \frac{b}{a'}e''''
\end{equation}
Indeed, we have, that $a(e'+e) = e'''-be''$ and $a'(e''+e') = e''''-b'e'''$. The second equality implies that $e'=\frac{1}{a'}(e''''-b'e'''-a'e'')$. Expressing thus $e$ and $e'$ in the r.h.s. in terms of $e''$, $e'''$ and $e''''$ one obtains the identity.

This identity implies that $a(e'+e)-be'$ corresponds to the intersection point $(p,p')\cap(p''',p'''')$. Slightly abusing notations by confusing vectors and the corresponding points on the projective plane one has
\begin{equation}\label{cross-ratios}
\begin{array}{l}\displaystyle
x=[e,\ e',\ e'+e,\ a(e'+e)-be'] = [0,\infty,1,1-\frac{b}{a}] = \frac{b}{a},
\\ \displaystyle
y=[e'''',\ e''',\ e''''-b'e''',\ (a'+bb')e'''-be'''']=[0,\infty,-b',-b'-\frac{a'}{b}]=-\frac{a'}{bb'}.
\end{array}
\end{equation}
Conversely $a=-(xx'y)^{-1}$ and $b=-(x'y)^{-1}$ or, restoring the indices one finally gets the formulae (\ref{xyab}), (\ref{abxy}).

Consider now the pentagram map. The expression for the evolved vectors and the normalisation condition for them in simplified notations look as
$$\begin{array}{l}
\tilde{e}=\alpha(e+e').\\[8pt]
\tilde{e}'''=\tilde{a}(\tilde{e}+\tilde{e}')+\tilde{b}\tilde{e}''.
\end{array}
$$
Substituting the first equation into the second one gets:
\begin{equation}
\alpha'''(e''''+e''')=\alpha\tilde{a}(e'+e)+\alpha'\tilde{a}(e''+e')+\alpha''\tilde{b}(e'''+e'').
\end{equation}
Taking into account the relation (\ref{normalization}): $e'''=a(e+e')+be''$,
one can eliminate $e$ and $e'$:
\begin{equation}
\alpha'''(e''''+e''')=\frac{\alpha\tilde{a}}{a}(e'''-be'')+\frac{\alpha'\tilde{a}}{a'}(e''''-b'e''')+\alpha''\tilde{b}(e'''+e'').
\end{equation}
Assuming that $e''''$, $e'''$ and $e''$ are linearly independent one gets the system of three equations:
\begin{equation}
\begin{array}{l}
\displaystyle 0=-\frac{\alpha b \tilde{a}}{a}+\alpha''\tilde{b}\\[8pt]
\displaystyle\alpha'''=\frac{\alpha\tilde{a}}{a}-\frac{\alpha'b'\tilde{a}}{a'}+\alpha''\tilde{b}\\[8pt]
\displaystyle\alpha'''=\frac{\alpha'\tilde{a}}{a'}.
\end{array}
\end{equation}
Eliminating $\tilde{a}$ and $\tilde{b}$ one gets
$$\alpha\frac{1+b}{a}=\alpha'\frac{1+b'}{a'}$$
implying that
$$\alpha= c\frac{a}{1+b},$$
where $c$ is a nonvanishing constant.

Finally we get
\begin{equation}
\displaystyle \tilde{a}=\frac{a'''(b'+1)}{b'''+1},\quad
\displaystyle \tilde{b}=\frac{a'''b(b'+1)(b''+1)}{a''(b+1)(b'''+1)}.
\end{equation}
Changing the variables to the Schwartz cross-ratio coordinates $x$ and $y$ given by (\ref{cross-ratios}) we get
$$
\tilde{x}=x''\frac{1-x'''y''}{1-x'y},\quad
\tilde{y}=y'\frac{1-x'y}{1-x'''y''}
$$
Changing to cluster coordinates $u$ and $v$ given by
$$ u'=-(xy)^{-1}\quad v=-x'y$$ we get
\begin{equation}
\tilde{u}'=v^{-1},\quad \tilde{v}=u''\frac{(1+v)(1+v''')}{(1+v'^{-1})(1+v''^{-1})}.
\end{equation}


\end{document}